\documentclass[12pt]{article}
\usepackage{amsfonts,amsmath,amssymb,amsthm,mathrsfs,yhmath,graphicx}
\usepackage{xy}
\xyoption{all}

\usepackage{ifpdf}
\ifpdf
\usepackage{epstopdf,hyperref}
\else
\usepackage[hypertex]{hyperref}
\fi

\newtheorem{Thm}{Theorem}
\newtheorem{Prop}[Thm]{Proposition}
\newtheorem{Lem}[Thm]{Lemma}

\theoremstyle{definition}
\newtheorem{Rem}[Thm]{Remark}
\newtheorem{Def}[Thm]{Definition}

\numberwithin{equation}{section}
\numberwithin{Thm}{section}

\newcommand{\bef}{\begin{Def}}
\newcommand{\eef}{\end{Def}}
\newcommand{\bethe}{\begin{Thm}}
\newcommand{\eethe}{\end{Thm}}
\newcommand{\berop}{\begin{Prop}}
\newcommand{\eerop}{\end{Prop}}
\newcommand{\belem}{\begin{Lem}}
\newcommand{\elem}{\end{Lem}}

\newcommand\be            {\begin{equation}}
\newcommand\bea           {\begin{equation}\begin{array}l\displaystyle}
\newcommand\bearll        {\begin{array}{ll}\displaystyle}
\newcommand\ee            {\end{equation}}
\newcommand\eear          {\end{array}}
\newcommand\enl           {\\[1em]\displaystyle}
\newcommand\etb           {&\!\! \displaystyle}
\newcommand\labl[1]       {\label{#1}\ee}
\newcommand\nxt{\noindent\raisebox{.08em}{\rule{.44em}{.44em}}\hspace{.4em}}

\newcommand\arxiv[2]      {\href{http://arXiv.org/abs/#1}{#2}}
\newcommand\doi[2]        {\href{http://dx.doi.org/#1}{#2}}
\newcommand\httpurl[2]    {\href{http://#1}{#2}}

\newcommand\eps           {\varepsilon}

\newcommand\id            {{\rm id}}

\newcommand\Cb            {\mathbb{C}}

\newcommand\Rb            {\mathbb{R}}
\newcommand\Zb            {\mathbb{Z}}

\newcommand\Cc            {\mathcal{C}}

\newcommand\Gc            {\mathcal{G}}

\newcommand\Ic            {\mathcal{I}}

\newcommand\Mc            {\mathcal{M}}

\newcommand\Oc            {\mathcal{O}}

\newcount\hour\newcount\minute
        \hour=\time \divide\hour by60 \minute=\time
        {\multiply\hour by60 \global\advance\minute by-\hour}
        \edef\militarytime{\number\hour:\ifnum\minute<10 0\fi\number\minute}

\newcommand\gt[1]         {\mathfrak{#1}}
\newcommand\sfk           {{\mathsf k}}
\newcommand\sfi           {{\mathsf i}}
\newcommand\sfd           {{\mathsf d}}
\newcommand\sfT           {{\mathsf T}}

\newcommand\qq            {\begin{eqnarray}}
\newcommand\qqq           {\end{eqnarray}}
\newcommand\tx[1]         {\textrm{#1}}

\newcommand\erm           {\textrm{e}}

\def\Gx{\textrm{G}}

\def\a{\alpha}

\def\vep{\varepsilon}
\def\D{\Delta}

\def\La{\Lambda}
\def\la{\lambda}
\def\Si{\Sigma}
\def\si{\sigma}

\def\Om{\Omega}
\def\om{\omega}

\def\th{\theta}

\def\txH{{\rm H}}

\def\x{\times}
\def\Ad{{\rm Ad}}

\def\emb{\hookrightarrow}

\def\p{\partial}

\def\con{\lrcorner\,}

\def\unl{\underline}

\def\tr{{\rm tr}}
\def\tgt{\gt{t}}

\def\xcV{\mathscr{V}}

\def\pr{{\rm pr}}
\def\txm{{\rm m}}

\def\cA{\mathcal{A}}

\def\cO{\mathcal{O}}

\def\cS{\mathcal{S}}

\def\ox{\otimes}

\def\Lie{{\rm Lie}}

\def\bC{{\mathbb{C}}}

\def\bR{{\mathbb{R}}}
\def\bZ{{\mathbb{Z}}}

\newcommand{\sug}{{\rm SU}(2)}
\newcommand{\sua}{\mathfrak{su}(2)}
\newcommand{\sugt}{\widetilde{\rm SU}(2)}
\newcommand{\faff}{P^\sfk_+}
\newcommand{\suk}{\widehat{\mathfrak{su}}(2)_\sfk}
\newcommand{\uj}{{\rm U}(1)}
\newcommand{\cG}{\mathcal{G}}
\newcommand{\curv}{{\rm curv}}
\newcommand{\triv}{{\boldsymbol{1}}}

\newcommand{\cT}{\mathcal{T}}

\newcommand{\cC}{\mathcal{C}}
\newcommand{\cM}{\mathcal{M}}
\newcommand{\cI}{\mathcal{I}}

\newcommand{\txp}{{\rm p}}
\newcommand{\txq}{{\rm q}}
\newcommand{\xcD}{\mathscr{D}}
\newcommand{\xcF}{\mathscr{F}}
\newcommand{\z}{\zeta}

\newcommand{\Tlmn}{\cT_{\la,\mu}^{\ \nu}}
\newcommand{\filmn}{\varphi_{\la,\mu}^{\ \nu}}
\newcommand{\adj}{\tx{ad}}
\newcommand{\Adj}{\tx{Ad}}

\newcommand{\barr}{\begin{array}}
\newcommand{\earr}{\end{array}}

\newcommand{\pLie}[1]{{-\hspace*{-11pt}\mathscr{L}\!}_{#1}}

\newcommand{\alxydim}[2]{\begin{aligned}\xymatrix#1{#2}\end{aligned}}

\begin{document}

\thispagestyle{empty}
\def\thefootnote{\fnsymbol{footnote}}
\begin{flushright}
KCL-MTH-09-08\\
\end{flushright}
\vskip 5.0em
\begin{center}\LARGE
Affine su(2) fusion rules from gerbe 2-isomorphisms
\end{center} \vskip 4em
\begin{center}\large
  Ingo Runkel\footnote{Present address: Department Mathematik,
Bereich Algebra und Zahlentheorie, Universit\"at Hamburg,
Bundesstra\ss e 55, 20146 Hamburg, Germany; Email: {\tt
ingo.runkel@uni-hamburg.de}}
  and
  Rafa\l ~R.\ Suszek\footnote{Present address: Katedra Metod
  Matematycznych Fizyki, Wydzia\l ~Fizyki Uniwersytetu
  Warszawskiego, ul.\,Ho\.za 74, PL-00-682 Warszawa, Poland;
  Email: {\tt suszek@fuw.edu.pl}}
\end{center}
\begin{center}
  Department of Mathematics, King's College London \\
  Strand, London WC2R 2LS, United Kingdom
\end{center}
\vskip 1em
\begin{center}
September 2009
\end{center}

\vskip 4em
\begin{abstract}
We give a geometric description of the fusion rules of the affine
Lie algebra $\,\suk\,$ at a positive integer level $\,\sfk\,$ in
terms of the $\sfk$-th power of the basic gerbe over the Lie group
$\,\sug$.\ The gerbe can be trivialised over conjugacy classes
corresponding to dominant weights of $\,\suk\,$ via a 1-isomorphism.
The fusion-rule coefficients are related to the existence of a
2-isomorphism between pullbacks of these 1-isomorphisms to a
submanifold of $\,\sug\x\sug\,$ determined by the corresponding
three conjugacy classes. This construction is motivated by its
application in the description of junctions of maximally symmetric
defect lines in the Wess--Zumino--Witten model.
\end{abstract}

\setcounter{footnote}{0}
\def\thefootnote{\arabic{footnote}}

\newpage

\tableofcontents

\section{Introduction}\label{sec:intro}

Consider the affine Lie algebra $\,\suk\,$ at a positive integer
level $\,\sfk$.\ The integrable highest-weight representations of
$\,\suk\,$ are labelled by elements of the set of dominant affine
weights which we identify with the subset
$\,\faff=\{0,1,\dots,\sfk\}\,$ of the integers. The fusion ring of
these representations is given by
\be
  [\la] \ast [\mu] = \sideset{}{^{{}_{(+2)}}}\sum_{\nu = |\la-\mu|}^{\min(\la+\mu,
  2\sfk-\la-\mu)}\, [\nu] \qquad
  \text{for $\,\la,\mu\in\faff$}\,,
\labl{eq:affine-fusion}
where the superscript $\,(+2)\,$ means that the sum is carried out
in steps of two. We can collect the fusion rules in a discrete
subset of $\,[0,\sfk]^{\times 3}$,
\be
  V = \{\ (\la,\mu,\nu)\in (\faff)^{\times 3} \ | \
  [\nu]~\text{appears in}~[\la] \ast [\mu] \ \}\,.
\labl{eq:setV-def}
We should, in principle, introduce a conjugation in this definition,
but for $\,\suk\,$ this does not make a difference.

The fusion rules appear in numerous places. The most interesting one
for us in the present context is the application in two-dimensional
conformal quantum field theory. There, they can be used to compute
the dimension of the spaces of conformal blocks in the quantum
Wess--Zumino--Witten (WZW) model at level $\,\sfk$,\ cf., e.g.,
\cite{Beauville:1994}, and they are famously related to the modular
properties of affine characters by the Verlinde formula
\cite{Verlinde:1988sn}.

The application we have in mind is the study of the WZW model on
surfaces with defect lines and defect junctions. Without going into
any detail, we merely mention that the elementary maximally
symmetric defects in the quantum WZW model are labelled by weights
in $\,\faff$.\ Defect lines can meet at points on the world-sheet
which we call defect junctions. A nonzero defect junction (of
conformal dimension zero) of three defect lines exists if and only
if the fusion rule for the three weights labelling the defect lines
is non-zero \cite{Frohlich:2006ch}.

Defect lines and defect junctions also have a description in
classical $\si$-models \cite{Fuchs:2007fw,Runkel:2008gr}. The
analysis of maximally symmetric defect lines and junctions in the
WZW model motivates the geometric construction of the present paper.
The relation to the WZW model will be elaborated in
\cite{Runkel:2009}.

We shall identify the fundamental affine Weyl alcove of $\,\sug\,$
at level $\,\sfk\,$ with the interval $\,[0, \sfk]\,$ on the real
axis, and assign to $\,\la\in [0,\sfk]\,$ the conjugacy class
\be
  \Cc_\la = \bigl\{\ g\cdot\bigl(\begin{smallmatrix}
  \erm^{\pi\sfi\la/\sfk} & 0 \\ 0 & \erm^{-\pi\sfi\la/\sfk}
  \end{smallmatrix}\bigr)\cdot g^{-1} \ \big\vert \ g\in\sug \
  \bigr\}\,.
\labl{eq:Clam-def}
We introduce the factor $\,\sfk\,$ at this stage in order to reduce
the number of its appearances later on. Define three maps $\,\txp_1,
\txp_2, \txm : \sug \times \sug \rightarrow \sug\,$ via
\be
  \txp_1(g,h) = g\,,\qquad\qquad\txp_2(g,h) = h\,,\qquad\qquad
  \txm(g,h) = g\cdot h\,.
\labl{eq:p1p2m-def}
We shall be interested in the following submanifolds of $\,\sug
\times \sug$:
\be
  \Tlmn
  = \txp_1^{-1}(\Cc_\la) \cap \txp_2^{-1}(\Cc_\mu) \cap
  \txm^{-1}(\Cc_{\nu})\,,\quad\text{where}\quad\la,\mu,\nu
  \in [0,\sfk]\,.
\labl{eq:setT-def}
It is easy to see that $\,\Tlmn \neq \emptyset\,$ if and only if
$\,(\Cc_\la \cdot \Cc_\mu) \cap \Cc_\nu \neq \emptyset$,\ and we show in
section \ref{sec:el-mult} that the $\,\Tlmn\,$ are in fact homogeneous spaces. This
defines the fusion polytope of $\,\sug$,
\be
  \xcF = \{\ (\la,\mu,\nu) \in [0,\sfk]^{\times 3} \ | \ \Tlmn \neq
  \emptyset \ \}\,,
\labl{eq:setF-def}
see Fig.\,\ref{fig:SU2-fus}. One finds by inspection (see Section
\ref{sec:el-mult} below) that if $\,[\nu]\,$ appears in the fusion
product $\,[\la]\ast[\mu]$,\ for $\,\la,\mu \in \faff$,\ then also
$\,(\la,\mu,\nu) \in\xcF$,\ i.e.\ $\,V \subset\xcF$.\ In fact, this
relation between intersections of conjugacy classes and affine
fusion rules holds more generally \cite{Hayashi:1999,Teleman:2000},
but we shall only consider $\,\sug\,$ here.

\begin{figure}[bt]

\begin{center}
\begin{picture}(180,220)
\put(0,0){
  \put(0,0){\includegraphics[height=18em]{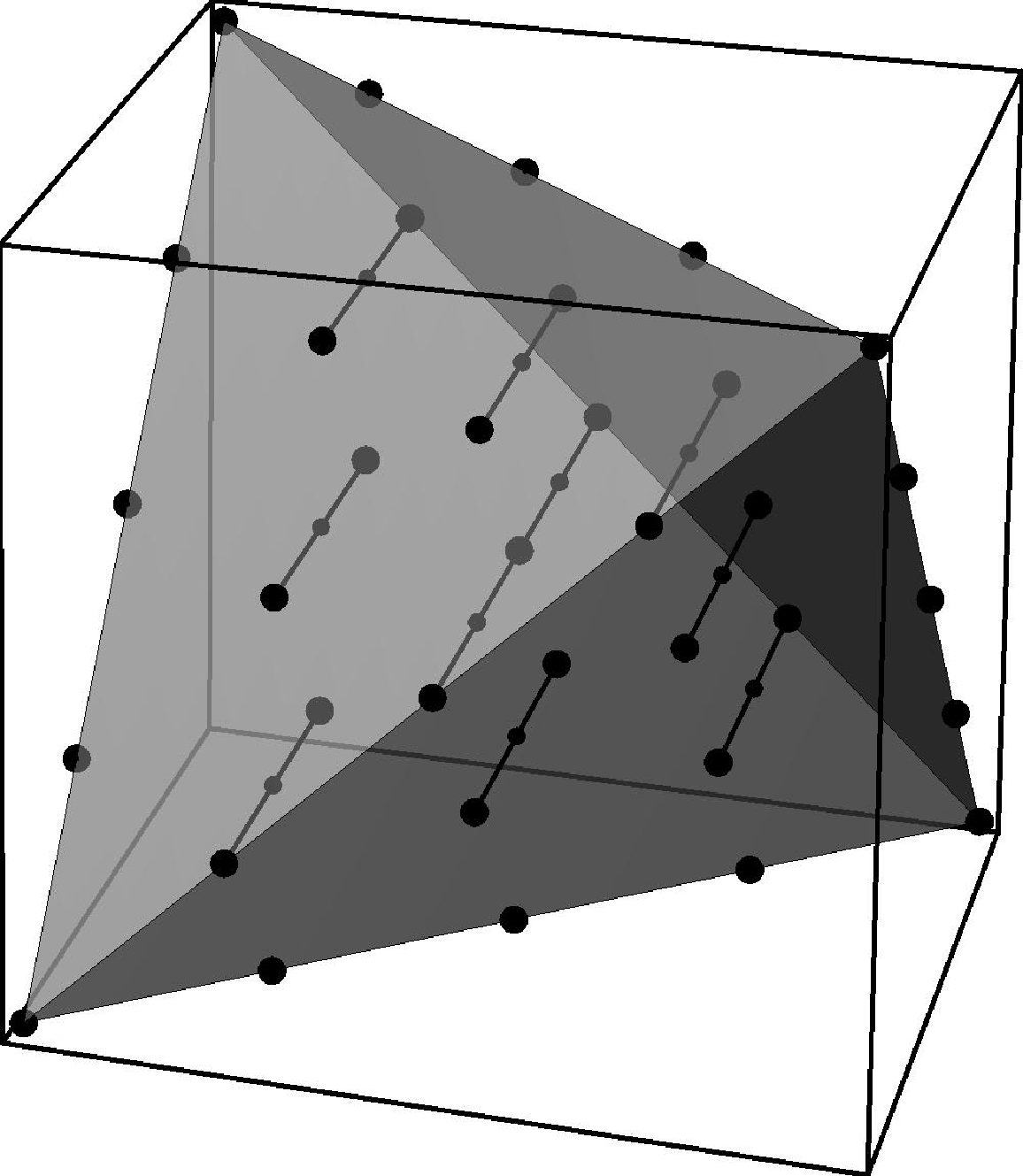}}
  \put(-12,16)     { $0$ }
  \put(135,-12)   { $\lambda=\sfk$ }
  \put(-35,161)   { $\mu=\sfk$ }
  \put(18,72)      { \scriptsize $\nu=\sfk$ }
}
\end{picture}
\end{center}

\caption{The fusion polytope of $\,\sug\,$ is a tetrahedron. Shown
here is the example of $\,\sfk=4$; the lines inside $\,\xcF\,$ are
the lines of constant $\,\la,\mu\in\faff$,\ the dots give the
intersection of $\,\xcF\,$ with $\,(\faff)^{\x 3}$,\ and the bold
dots mark points for which $\,[\nu] \in [\la] \ast [\mu]$,\ i.e.\
the intersection $\,\xcF \cap V$.\ We give a construction based on
the basic gerbe of $\,\sug\,$ which singles out the bold dots.}
\label{fig:SU2-fus}
\end{figure}
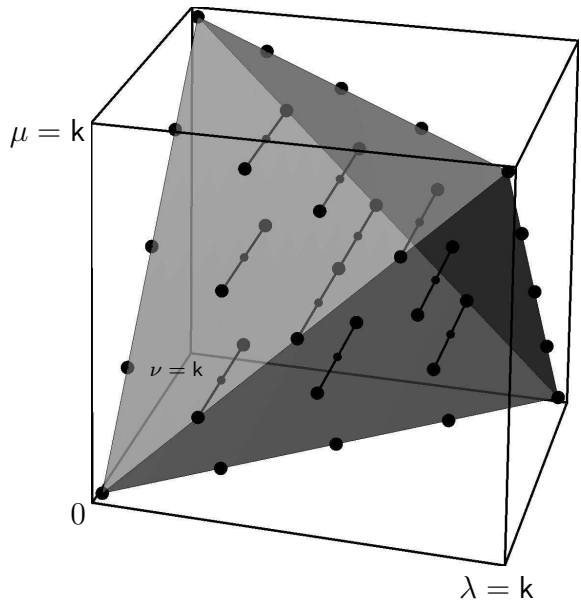

One can now ask about the converse, namely if one can determine
$\,V\,$ starting from $\,\xcF$.\ The simplest idea would be to just
intersect $\,\xcF\,$ with the set of dominant affine weights
$\,(\faff)^{\x 3}$,\ but this does not respect the
parity-conservation rule by which the sum \eqref{eq:affine-fusion}
has to be carried out in steps of two. Obtaining this rule from
geometric considerations is the main aim of this paper. To this end,
we shall need to introduce some additional structures on $\,\sug$,\
to wit, a certain gerbe $\,\Gc\,$ and related 1- and 2-isomorphisms.
The construction is summarised below, and then detailed in Sections
\ref{sec:basgerb}--\ref{sec:fusmorph}.

\medskip

The Lie group $\,\sug\,$ is equipped with the family of Cartan
3-forms
\qq
\txH=\frac{r}{12\pi}\,\tr(\th_L\wedge\th_L\wedge\th_L)\,,\qquad
\qquad r \in \Rb\,,
\qqq
defined in terms of the standard left-invariant Maurer--Cartan
1-forms $\,\th_L(g)=g^{-1}\,\sfd g\,$ on $\,\sug$.\ There is a gerbe
on $\,\sug\,$ (or on any compact simple connected and simply
connected Lie group, for that matter) with curvature $\,\txH\,$ if
and only if $\,r\in\Zb$,\ and this gerbe is unique up to a
1-isomorphism \cite{Gawedzki:1987ak,Meinrenken:2002}. We shall set
$\,r=\sfk\,$ and denote the corresponding gerbe $\,\Gc$.\ We review
the definition of bundle gerbes and the detailed construction of
$\,\Gc\,$ in Section \ref{sec:basgerb} below.

Gerbes have an important application in physics: they describe the
topological term in the $\si$-model action functional which is
necessary to preserve the conformal symmetry of the WZW model upon
quantisation
\cite{Witten:1983ar,Alvarez:1984es,Gawedzki:1987ak,Felder:1988sd}.
\medskip

If we restrict the 3-form $\,\txH\,$ to a conjugacy class
$\,\Cc_\la \subset \sug\,$ for $\,\la \in [0,\sfk]$,\ it
becomes exact, $\,H|_{\Cc_\la} = \sfd \omega_\la\,$ for
$\,\om_\la\in\Om^2(\cC_\la)$.\ One can then enquire whether also the
gerbe $\,\Gc\,$ can be trivialised when restricted to
$\,\Cc_\la$,\ i.e.\ whether there is a 1-isomorphism
\be
  \Phi_\la\ :\ \Gc\big|_{\Cc_\la} \rightarrow \Ic(\omega_\la)
\labl{eq:Philam-def}
between $\,\Gc|_{\Cc_\la}\,$ and the trivial gerbe with curving
$\,\omega_\la\,$ (consult Sections \ref{sec:basgerb} and
\ref{sec:triv} for definitions). This turns out to be possible if
and only if $\,\la \in \faff$,\ and $\,\Phi_\la\,$ is unique up to a
2-isomorphism in this case \cite{Gawedzki:2002se,Gawedzki:2004tu}
(again, a similar statement holds for other compact simple connected
and simply connected Lie groups). We review the construction of
$\,\Phi_\la\,$ in Section \ref{sec:triv}.

The trivialisation of the gerbe over $\,\Cc_\la\,$ finds its
application in the description of the WZW model on surfaces with a
non-empty boundary. There, it describes a maximally symmetric
boundary condition for the $\si$-model fields
\cite{Gawedzki:2002se}, and, more relevant to our present concerns,
it also describes the WZW model on surfaces with defect lines
\cite{Fuchs:2007fw,Runkel:2008gr}, where it defines a maximally
symmetric defect gluing condition
\cite{Waldorf:2008mult,Runkel:2009}.
\medskip

Each of the three maps $\,\txp_1,\txp_2\,$ and $\,\txm\,$ defined in
\eqref{eq:p1p2m-def} can be used to pull back the gerbe $\,\Gc\,$
from $\,\sug\,$ to $\,\sug \times \sug$.\ There is a unique
1-isomorphism $\,\cM : \txp_1^* \Gc \star \txp_2^*\Gc \rightarrow
\txm^* \Gc \star \Ic(\rho)$,\ where $\star$ is a product of gerbes
(described in Section \ref{sec:basgerb}) and $\,\rho\,$ is a 2-form
defined globally on $\,\sug \times \sug$.\ The 1-isomorphism
$\,\cM\,$ (together with a preferred 2-isomorphism subject to
certain conditions) is called a multiplicative structure on
$\,\cG\,$ \cite{Carey:2004xt,Waldorf:2008mult}. It enters the final
step of our construction: Given $\,\la,\mu,\nu \in \faff$,\ we use
the maps $\,\txp_1,\txp_2\,$ and $\,\txm\,$ to pull back the
1-isomorphisms $\,\Phi_\la,\Phi_\mu\,$ and $\,\Phi_\nu\,$ from their
respective conjugacy classes to $\,\Tlmn$.\ We then ask if there
exists a 2-isomorphism
\be
\filmn\ :\ \txp_1^*\Phi_\la\star\txp_2^*\Phi_\mu
\Longrightarrow\bigl(\txm^*\Phi_\nu\star\id_{\cI(\rho)}\bigr)\circ\cM
\labl{eq:fusion-2mor-def}
between the two 1-isomorphisms trivialising $\,\txp_1^* \Gc \star
\txp_2^*\Gc\,$ over $\,\Tlmn$,\ as detailed in Section
\ref{sec:fusmorph} (the product and composition of 1-morphisms is
defined in Section \ref{sec:triv}). The existence of the 2-isomorphism is obstructed
in general, and we find by inspection that the obstruction vanishes
if and only if $\,[\nu]\,$ appears in the fusion product $\,[\la]
\ast [\mu]$.\ Thus, if we define
\be
  V_\cG = \{\ (\la,\mu,\nu)\in (\faff)^{\times 3} \ | \
  \filmn \text{ exists}\ \}\,,
\labl{eq:setV'-def}
we can state our main result as

\bethe\label{thm:V=V}
$V=V_\cG$,\ where $\,V\,$ is given in \eqref{eq:setV-def} and
$\,V_\cG\,$ in \eqref{eq:setV'-def}. \eethe

The definition \eqref{eq:fusion-2mor-def} and -- in particular --
the appearance of the multiplicative structure may seem ad hoc, but
they will turn out to be very natural from the point of view of the
WZW model in the presence of defect lines. In this context, the
2-isomorphism \eqref{eq:fusion-2mor-def} is used to define a
junction of maximally symmetric defect lines \cite{Runkel:2008gr},
as we shall describe in detail in \cite{Runkel:2009}.

\medskip

Let us finish this introduction with a general remark in which we
try to put our result into context. It has always been surprising to
what extent it is possible to gain information about the quantum WZW
model by considering the underlying classical sigma model.
Conversely, knowing the exact solution of the quantum WZW model from
methods in representation theory has inspired the search for new
geometric concepts which give geometrical explanations for the
observed quantum effects. Indeed, this is the very origin of the
basic gerbe on a Lie group \cite{Gawedzki:2002se,Meinrenken:2002},
as well as gerbe modules \cite{Carey:2002} and the bi-categorial
structure on gerbes considered in
\cite{Waldorf:2007mm,Fuchs:2007fw}. In this paper, we extend this
dictionary -- for the first time -- to fusion rules. Previous
partially geometric constructions identifying the fusion rules
always used quantisation methods at some point (see e.g.\
\cite{Bismut:1999,Alekseev:2000} and references therein), while our
approach is purely geometrical. As our result is obtained only in
the special case of $\,\sug\,$ and since we require over forty pages
of machinery and, in some places, very technical calculations to
arrive at a $\,\Zb_2\,$ selection rule, we regard the present result
as a proof of concept rather than as the final answer. Currently, we
have not much to say about how Theorem \ref{thm:V=V} generalises to
other Lie groups, in particular if or how fusion-rule multiplicities
could appear in the context of gerbes. But we expect that such a
generalisation exists and that a more elegant and conceptual proof
than ours -- which boils down to the comparison of two lists
computed via different methods -- can be found.

\medskip

This paper is organised as follows: In Sections \ref{sec:basgerb}
and \ref{sec:triv}, we review the construction of the gerbe
$\,\cG\,$ and of the 1-isomorphisms trivialising $\,\cG\,$ upon
restriction to conjugacy classes, following \cite{Gawedzki:2002se}.
Our first new result is the explicit construction of the
1-isomorphism $\,\cM\,$ restricted to $\,\Tlmn$ in Section
\ref{sec:el-mult};\ previously, only the existence of $\,\cM\,$ had
been proved. In Section \ref{sec:fusmorph}, we analyse when
$\,\filmn\,$ exists and prove Theorem \ref{thm:V=V}.

\bigskip\noindent
{\bf Acknowledgements:} We thank J.~Fuchs,
K.~Gaw{\c{e}}dzki and C.~Schweigert for discussions and sustained
interest in our work. This research was partially supported by the
EPSRC First Grant EP/E005047/1 and the STFC Rolling Grant
ST/G000395/1.
\newpage

\section{The basic gerbe over SU(2)}\label{sec:basgerb}

Abelian gerbes \cite{Giraud:1971} over a given manifold $\,M\,$
provide a geometric realisation of elements of the integral
cohomology group $\,H^3(M,\bZ)$.\ They admit a number of
descriptions, such as the local description
\cite{Gawedzki:1987ak,Brylinski:1993ab} in terms of Deligne
hypercohomology, and the geometric one
\cite{Murray:1994db,Murray:1999ew} formulated in terms of complex
line bundles over a surjective submersion over $\,M$.\ In the
context of two-dimensional non-linear $\si$-models, the former
description was the first one to be used as it is more intuitive and
can be abstracted from the construction of the action functional
\cite{Gawedzki:1987ak}. Nonetheless, the geometric description is
more convenient for our explicit calculations, and we shall
use it throughout the main text.

Let us begin our lightning review of the necessary elements of the
geometric description with the concept of an Abelian bundle gerbe
with curving, connection and fixed curvature, as represented
concisely by the following diagram \cite{Stevenson:2000wj}
\qq
\alxydim{@C=1cm@R=1cm}{\hspace{10pt}(L,\nabla_L,\mu_L)
\ar[d]_{\pi_L} & \\ Y^{[2]}M \ar@<.5ex>[r]^{\pr_1}
\ar@<-.5ex>[r]_{\pr_2} & (YM,B) \hspace{-14pt}\ar[d]^{\pi_{YM}} \\ &
(M,\txH)}\quad.
\qqq
Here, $\,Y^{[2]}M\,$ is the fibred product
\qq
Y^{[2]}M\equiv YM\x_M YM=\{\ (y_1,y_2)\in YM\x YM \ \vert \
\pi_{YM}(y_1)=\pi_{YM}(y_2) \ \}\,,
\qqq
with its obvious generalisations $\,Y^{[n]}M\,$ for $\,n>2$.\ We
shall also need the various canonical projections
\be
  \pr_{i_1,\dots,i_n}\ :\ Y^{[m]}M \rightarrow Y^{[n]}M
  \ :\ (y_1,\dots,y_m)\mapsto (y_{i_1},\dots,y_{i_n}) \,.
\ee
For example, $\,\pr_i : Y^{[2]}M\to YM\,$ maps $\,(y_1,y_2)\,$ to
$\,y_i$.\ The notation does not keep track of the index $m$ or of
the manifold $\,M$,\ and so we shall mention the source and the
target explicitly if it is not clear from the context.

\bef\cite{Murray:1994db}\label{def:gerbe}
A {\it Hermitian Abelian bundle gerbe with curving and connection}
(or {\it gerbe} for short) of curvature $\,\txH\in\Om^3(M)\,$ over a
smooth base $\,M\,$ is a quadruple $\,\cG=(YM,B,L,\mu_L)$ with the
following entries:
\begin{itemize}
\item a surjective submersion\footnote{A differentiable map $\,f:X
\to Y\,$ between a pair of smooth manifolds $\,X,Y\,$ is called a
{\em surjective submersion} if both $\,f\,$ and its tangent map $\,f_*:
\sfT X\to\sfT Y\,$ are surjective, and so $\,f\,$ admits smooth
local sections.} $\,YM\xrightarrow{\pi_{YM}}M$,\ together with a
global 2-form $\,B\in\Om^2(YM)\,$ on $\,YM$,\ termed the curving of
the gerbe, satisfying the relation
\qq
\pi_{YM}^*\txH=\sfd B\,;
\qqq
\item a Hermitian line bundle $\,L\xrightarrow{\pi_L}Y^{[2]}M\,$
with unitary connection $\,\nabla_L\,$ of curvature
\qq
\curv(\nabla_L)=\pr_2^*B-\pr_1^*B\,,
\qqq
determined by the pullbacks $\,\pr_i^*B\,$ of the curving along
$\,\pr_i : Y^{[2]}M\to YM$;
\item a groupoid structure $\,\mu_L\,$ on
$\,L\xrightarrow{\pi_L}Y^{[2]}M$,\ i.e.\ a unitary
connection-preserving bundle isomorphism
\qq
\mu_L\ :\ \pr_{1,2}^*L\ox\pr_{2,3}^*L \xrightarrow{\sim} \pr_{1,3}^*
L
\qqq
between pullbacks of $\,L\,$ to $\,Y^{[3]}M\,$ along $\,\pr_{i,j} :
Y^{[3]}M\to Y^{[2]}M$,\ which is associative in the sense specified
by the commutativity of the diagram
\qq
\alxydim{@C=8em@R=5em}{\pr_{1,2}^*L\ox\pr_{2,3}^*L\ox\pr_{3,4}^*L
\ar[d]_{\pr_{1,2,3}^*\mu_L\ox\id_{\pr_{3,4}^*L}}\ar[r]^{\hspace{30pt}
\id_{\pr_{1,2}^*L}\ox\pr_{2,3,4}^*\mu_L} &
\pr_{1,2}^*L\ox\pr_{2,4}^*L \ar[d]^{\pr_{1,2,4}^*\mu_L} \\
\pr_{1,3}^*L\ox\pr_{3,4}^*L \ar[r]_{\hspace{12pt}\pr_{1,3,4}^*\mu_L}
& \pr_{1,4}^*L}
\qqq
of bundle isomorphisms over $\,Y^{[4]}M$,\ the latter coming with
the canonical projections $\,\pr_{i,j}: Y^{[4]}M\to Y^{[2]}M\,$ and
$\,\pr_{i,j,k}:Y^{[4]}M\to Y^{[3]}M$.
\end{itemize}
\eef

Amidst all gerbes over a given base $\,M$,\ there is a distinguished
class of trivial gerbes, characterised by the existence of a globally
defined curving $\,\om\in\Om^2(M)$.\ They are represented by
quadruples
\qq
\cI(\om)=(M,\om,\triv_M,\mu_{\triv_M})\,,
\qqq
with the surjective submersion $\,\id_M:M\to M$,\ the trivial bundle
$\,\triv_M=M\x\bC\to M\,$ with the trivial connection
$\,\nabla_{\triv_M}=\sfd$,\ and the canonical groupoid structure
\qq
\mu_{\triv_M}\ :\ \triv_M\ox\triv_M\to\triv_M\ :\ (x,z)\ox(x,z')
\mapsto(x,z\cdot z')\,.
\qqq

The following two natural operations defined on gerbes will be
useful for our purposes.

\medskip\noindent{\bf Pullback:}
Given a smooth map $\,f:N\to M\,$ from a manifold $\,N\,$ to the
base $\,M\,$ of the gerbe $\,\cG$,\ and a surjective submersion
$\,\pi_{Y N}:Y N \to N\,$ together with a map $\,\widehat f:Y N\to
YM\,$ that covers $\,f\,$ in the sense that it renders the diagram
\qq\label{diag:pullback}
\alxydim{@C=5em@R=3em}{Y N \ar[d]_{\pi_{Y N}} \ar[r]^{\widehat
f} & YM \ar[d]^{\pi_{YM}} \\ N \ar[r]_{f} & M}
\qqq
commutative (and so it also induces maps $\,\widehat f^{[n]}:
Y^{[n]}N\to Y^{[n]}M$), we define the {\em pullback of} $\,\cG\,$
{\em to} $\,N\,$ {\em along} $\,f\,$ as the gerbe
\qq
f^*\cG=\bigl(Y N,\widehat f^*B,\widehat f^{[2]*}L,\widehat
f^{[3]*}\mu_L\bigr)\,.
\qqq

\medskip\noindent{\bf Product:}
Consider a pair of gerbes $\,\cG_i=(Y_i M, B_i,L_i,\mu_{L_i}),\
i\in\{1,2\}\,$ over a common base $\,M$.\ We define their {\em
product} as the gerbe
\qq\label{eq:gerbe-prod}
\cG_1\star\cG_2=\bigl(Y_{1,2}M,\pr_1^*B_1+\pr_2^*B_2,\pr_{1,3}^*L_1
\ox\pr_{2,4}^*L_2,\pr_{1,3,5}^*\mu_{L_1}\ox\pr_{2,4,6}^*\mu_{L_2}
\bigr)
\qqq
with the surjective submersion
\qq
\pi_{Y_i M}\circ\pr_i\ :\ Y_{1,2}M=Y_1 M\x_M Y_2 M\to M
\qqq
expressed in terms of either of the canonical projections
\qq
\pr_i\ :\ Y_{1,2}M\to Y_i M\ :\ (y_1,y_2)\mapsto y_i
\qqq
and their obvious generalisations to higher fibred products
$\,Y_{1,2}^{[n]}M$.

\medskip

The basic gerbe $\,\cG_1\,$ over the group manifold of $\,\sug\,$ is
the unique (up to a 1-isomorphism, see Section \ref{sec:triv}) gerbe
with the curvature given by the Cartan 3-form
\qq
\curv(\cG_1)=\frac{1}{12\pi}\,\tr\bigl(\th_L\wedge\th_L\wedge\th_L
\bigr)\,,\qquad\qquad\th_L(g)=g^{-1}\,\sfd g\,,\qquad g\in\sug\ .
\qqq
The class of the rescaled 3-form $\,\tfrac{1}{2\pi}
\,\curv(\cG_1)\,$ is a generator of the integral cohomology group
$\,H^3\bigl(\sug,\bZ\bigr)\cong\bZ\,$ embedded in
$\,H^3\bigl(\sug,\bR \bigr)$.\

Let us fix a positive integer $\,\sfk\,$ for the rest of the paper.
We shall be concerned with a gerbe $\,\cG\,$ of $\,\sfk\,$ times the
curvature of $\,\cG_1$,\ which is otherwise constructed in the same
manner as $\,\cG_1$.\ We shall refer to $\,\cG\,$ as the {\em
canonical gerbe}.
It is 1-isomorphic to the $\sfk$-th $\star$-power of the basic gerbe.
An explicit construction of the canonical gerbe
over $\,\sug\,$ in terms of its local data (i.e.\ in the language of
the Deligne hypercohomology) was originally given in
\cite{Gawedzki:1987ak}. In the remainder of this section, we explain
in some detail the geometric construction along the lines of
\cite{Gawedzki:2002se}, providing the level of detail that we need
for Sections \ref{sec:triv}--\ref{sec:fusmorph}.

We shall first set up our conventions and assemble the
necessary algebraic objects. The standard Pauli
matrices are
\qq
\si_1=\left(\begin{smallmatrix} 0 & 1 \\ 1 & 0
\end{smallmatrix}\right)\,,\qquad\qquad\si_2=\left(
\begin{smallmatrix} 0 & -\sfi \\ \sfi & 0 \end{smallmatrix}\right)\,,
\qquad\qquad\si_3=\left(\begin{smallmatrix} 1 & 0 \\ 0 & -1
\end{smallmatrix}\right)\ .
\qqq
In the complexified Lie algebra $\,\sua^\Cb$,\ we choose the Cartan
subalgebra $\,\tgt=\bC\,\si_3$ and fix the $\adj$-invariant bilinear
(Killing) form given by the trace as
\qq
K(X,Y)=\tr(X\cdot Y)\,,\qquad\qquad X,Y\in\sua^\Cb\,.
\qqq
We subsequently use the latter to identify the Lie algebra with its
dual. Under this identification, the weight lattice of the algebra
becomes $\,P=\tfrac{1}{2}\,\bZ\,\si_3$, with the single fundamental
weight $\,\La=\tfrac{1}{2}\,\si_3$,\ and the root lattice takes the
form $\,Q=\bZ\,\si_3$, with the single simple root $\,\a=\si_3$.\
The fundamental affine Weyl alcove at level $\,\sfk\,$ is the
1-simplex
\qq
\cA_{\rm W}^\sfk=\{\ \la\,\La \ \vert \ \la\in[0,\sfk] \ \}\,,
\qqq
whence also its identification with the closed segment
$\,[0,\sfk]\,$ mentioned in the introduction and used throughout the
paper. Upon exponentiation
\qq
\erm^{2\pi\sfi\la\,\La/\sfk}=\Big(\begin{smallmatrix}
\erm^{\pi\sfi\la/\sfk} & 0 \\ 0 &
\erm^{-\pi\sfi\la/\sfk}\end{smallmatrix}\Big)=:t_\la\,,
\qqq
it produces an interval in the Cartan subgroup $\,T=\uj\,$ which
intersects each $\Adj$-orbit at exactly one point.

\subsection*{Parametrisations of $\,\boldsymbol\sug$}

In what follows, we shall make use of two particularly convenient
parametrisations of the group manifold. The first of them is a
redundant parametrisation given by the surjective map
\qq\label{eq:su2-param}
c\ :\ \sugt=[0,\sfk]\x\sug\to\sug\ :\
(\la,h)\mapsto\Ad_h(t_\la)\,.
\qqq
Since $\,c(\la,h) = c(\la',h')\,$ implies $\,\la=\la'$,\ we can
consider the isotropy group for each value of $\,\la\,$ separately,
\qq
\sug_{t_\la}=\{\ g\in\sug \ \vert \ \Ad_g(t_\la)=t_\la \ \}\,.
\qqq
It acts on $\,\sug\,$ by right regular translations
\qq\label{eq:iso-act}
\sug\x\sug_{t_\la}\to\sug\ :\ (h,t)\mapsto h\cdot t\,,
\qqq
so that we have $\,c(\la,h \cdot t) = c(\la,h)$.\ Many of the
expressions below will simplify once pulled back to $\,\sugt\,$.
However, to ensure that functions, forms and bundles on $\,\sugt\,$
arise as pullbacks from $\,\sug\,$, one has to impose appropriate
equivariance conditions on these objects over each of the
submanifolds $\,\{\la\}\x\sug\,$ with respect to the isotropy groups
$\,\sug_{t_\la}$.\ We shall denote objects defined on $\,\sugt\,$
with a tilde $\,\widetilde{~~}$.

\medskip

The other useful parametrisation is in terms of the Euler angles
$\,\varphi\in[0,\pi[,\ \th\in[0,\tfrac\pi2],\ \psi\in [0,2\pi[$.\
Here, an arbitrary group element $\,g\in\sug\,$ is expressed as
\qq\label{eq:Euler-param}
g(\varphi,\th,\psi)=\erm^{\sfi\varphi\si_3}\cdot\erm^{\sfi\th\si_2}
\cdot\erm^{\sfi\psi\si_3}\,.
\qqq

\subsection*{Surjective submersion and curvings}

The construction of the canonical gerbe $\,\cG\,$ begins with the
choice of a cover of the group manifold.
We take a pair of contractible open subsets ($e\,$ is the
group unit)
\qq
\cO_0\cong\sug\setminus\{-e\}\,,\qquad\qquad\cO_1\cong\sug\setminus
\{e\}
\qqq
within $\,\sug$.\ The surjective submersion is defined as the
disjoint union of the two elements of the open cover,
\qq\label{eq:cover}
Y\sug=\cO_0\sqcup\cO_1\xrightarrow{\pi_{Y\sug}}\sug\,.
\qqq
This choice is related to the
shape of the fundamental affine Weyl alcove via the pullback to
$\,\sugt$.\ Namely, $\,\cO_0\,$ and $\,\cO_1\,$ are
associated with the respective vertices
\qq
\la_0=0\,,\qquad\qquad\la_1=\sfk
\qqq
of $\,[0,\sfk]\cong\cA_{\rm W}^\sfk\,$ as per
\qq\label{eq:tilde-O}
\widetilde\cO_i:=c^{-1}(\cO_i)=\bigl([0,\sfk]\setminus\{\la_{1-i}\}
\bigr)\x\sug\,.
\qqq

\medskip

Next, we need to find a global primitive for the curvature
\be
  \txH = \curv(\cG) = \frac{\sfk}{12\pi}\,\tr\bigl(\th_L\wedge\th_L
  \wedge\th_L\bigr)
\ee
of $\,\cG\,$ when pulled back to $\,Y\sug$,\ i.e.\ we need 2-forms
$\,B_i\,$ on $\,\cO_i\,$ such that $\,\sfd B_i = \txH|_{\cO_i}$.\ We
shall construct these by giving a primitive of the pullback of
$\,\txH\,$ by $\,c\,$ and then checking the necessary equivariance
conditions.

Using the Maurer--Cartan equation $\,\sfd\th_L+\th_L\wedge\th_L=0$,\
the pullback $\,c^*\txH\,$ is readily verified to trivialise
globally on $\,\sugt\,$ as
\qq
c^*\txH(\la,h)=\sfd\bigl(\widetilde Q(\la,h)+\widetilde
F(\la-\la_c,h)\bigr)\,,\label{HQt}
\qqq
with
\qq
\widetilde
Q(\la,h)=\tfrac{\sfk}{4\pi}\,\tr\bigl(\th_L(h)\wedge\Ad_{t_\la}
\th_L(h)\bigr)\,,\qquad\qquad \widetilde
F(\la,h)=-\sfi\,\la\,\tr\bigl(\La\,\sfd\th_L(h)\bigr)\,,
\qqq
and $\,\la_c\,$ an arbitrary constant. For a form $\,\eta\,$ on
$\,\sugt\,$ to be the pullback of a form on $\,\sug$,\ it
necessarily has to be basic with respect to the action of the local isotropy
groups $\,\sug_{t_\la}$.\ This means that it has to be horizontal and
invariant for each $\,\la\in [0,\sfk]$,\ that is, for all vector
fields on $\,\{\la\}\x\sug\,$ of the form
\qq\label{eq:XL-vec-flds}
X(h)=X^A\,L_A(h)\,,\qquad\qquad X^A\,\sfi\,\si_A\in
\Lie\sug_{t_\la}\,,
\qqq
it has to satisfy
\qq\label{eq:B-horinv}
X \con \eta=0
\qquad \text{and} \qquad
\pLie{X}\eta=0\,.
\qqq
Here, $\,L_A\,$ are the standard left-invariant vector fields on
$\,\sug\,$ dual to the Maurer--Cartan 1-forms,
\qq
L_A\con\th_L=\sfi\,\si_A\,,
\qqq
and $\,\pLie{X}\,$ is the Lie derivative along $\,X$.\ The vector
fields in \eqref{eq:XL-vec-flds} generate the action
\eqref{eq:iso-act} of the isotropy groups $\,\sug_{t_\la}$.\ For
$\,\la \in ]0,\sfk[$,\ the Lie algebra $\,\Lie\sug_{t_\la}\,$ is
spanned by $\,\sfi\,\si_3\,$ and we find that
\qq
L_3(h)\con \widetilde
Q(\la,h)=\tfrac{\sfk}{4\pi}\,\tr\bigl(\sfi\,\si_3\,\Ad_{t_\la}
\th_L(h)\bigr)-\tfrac{\sfk}{4\pi}\,\tr\bigl(\th_L(h)\,\Ad_{t_\la}
\sfi\,\si_3\bigr)=0\,,
\qqq
and, similarly, $\,L_3(h)\con \widetilde F(\la,h)=0$.\ Evaluating
the Lie derivatives as $\,\pLie{X}\eta = X \con \sfd\eta + \sfd(X
\con \eta)$,\ we find
\qq
\pLie{L_3}(\widetilde Q+\widetilde F)(\la,h) = L_3(h) \con c^*
\txH(\la,h)=0\,,
\qqq
and hence, altogether, $\,\widetilde Q(\la,h)+\widetilde F(\la,h)\,$ is basic on
$\,]0,\sfk[ \x\sug$.

Furthermore, $\,\sfd \widetilde
F(\la,h)=-\sfi\,\sfd\la\wedge\tr\bigl(\La\, \sfd\th_L(h)\bigr)\,$
gives $\,\pLie{L_3} \widetilde F(\la,h) = 0$.\ This implies that
also $\,\pLie{L_3} \widetilde Q(\la,h) = 0$.\ Since $\,\widetilde
Q(\lambda_i,h)=0\,$ for $\,i\in\{0,1\}$,\ it is automatically basic
for these two values of $\,\la$.\ On the other
hand, $\,\widetilde F(\la_i-\la_c,h)\,$ is horizontal and invariant
with respect to $\,\sug_{t_{\la_i}} = \sug\,$ only for
$\,\la_c=\la_i\,$ (in which case $\,\widetilde F(\la_i-\la_c,h)=0$).
This prompts us to define
\qq\label{eq:Bi-result}
\widetilde B_i(\la,h)=\widetilde Q(\la,h)+\widetilde F(\la-\la_i,h)\,,
\qqq
such that $\,\widetilde B_i(\la,h)\,$ is basic on
$\,\widetilde\cO_i\,$ and satisfies $\,\widetilde B_i(\la_i,h)=0$.

If the stabiliser is constant, a basic form is automatically the
pullback of a smooth form on the quotient. If the stabiliser jumps,
as it does here at $\,\lambda=0\,$ and $\,\lambda=\sfk$,\ one has
to verify additionally smoothness of the descended form.
For $\,\widetilde B_i$,\ this was done in \cite{Gawedzki:2002se}
using the Poincar\'e lemma. The $\,\widetilde B_i\,$ are hence pullbacks to
\qq
Y\sugt=\widetilde\cO_0\sqcup\widetilde\cO_1
\qqq
of the sought-after curvings $\,B_i\,$ on $\,\Oc_i\,$ along the map
\qq \label{eq:c-hat-def}
\widehat c\ :\ Y\sugt\to Y\sug\ :\ (i,\la,h)\mapsto\bigl(i,
\Ad_h(t_\la)\bigr)
\qqq
that covers $\,c$.

\subsection*{Kirillov--Kostant--Souriau bundles}

The line bundle over $\,Y^{[2]}\sug\,$ will be constructed in terms
of the Kirillov--Kostant--Souriau (KKS) bundles over co-adjoint
orbits of $\sug$.\ Denote by
\qq
\sug_{\la}=\{\ g_{\la}\in\sug \ \vert \ \Ad_{g_{\la}}
(\la\,\La)=\la\,\La \ \}
\qqq
the isotropy group of the point $\,\la\,\La\,$ with respect to the
co-adjoint action of $\,\sug$.\ The isotropy group acts on
$\,\sug\,$ through right regular translations,
\qq\label{eq:iso-act-right}
\sug\x\sug_{\la}\to\sug\ :\ (h,t)\mapsto h\cdot t\,.
\qqq
The corresponding co-adjoint orbit $\,\sug/\sug_\la\,$ is endowed
with a canonical symplectic structure, and that structure defines a
canonical line bundle $\,K_\la\to\sug/\sug_\la$,\ the KKS bundle
\cite{Kostant:1970,Souriau:1970,Kirillov:1975}, which we proceed to
describe. We start from the trivial line bundle
\qq
\pr_1\ :\ \overline K_\la=\sug\x\bC\to\sug
\qqq
with connection
\qq\label{eq:KKS-con}
\nabla=\sfd+\la\,\tr(\La\,\th_L)
\qqq
of curvature
\qq\label{eq:KKS-curv}
\curv(\nabla)= \sfi\la\,\tr\bigl(\La\,\sfd \th_L \bigr)\,.
\qqq
For $\,\la \in \Zb$,\ we can lift the action
\eqref{eq:iso-act-right} to the bundle $\,\overline K_\la\,$ as
\be
\overline K_\la \x \sug_{\la} \to \overline K_\la\ :\ (h,z,t)\mapsto\bigl(h\cdot t,
\chi_{\la}(t)\cdot z\bigr)
\labl{eq:KKS-S-action}
by means of the characters
\qq\label{eq:char-la}
\chi_\la\ :\ \sug_{\la} \to \Cb^\times\ :\
\erm^{\sfi\varphi\si_3}\mapsto \erm^{-\sfi\la\varphi} \quad
\text{for}~ \la \neq 0\,.
\qqq
For $\,\la=0$,\ we have $\,\sug_{\la} = \sug\,$ and it is convenient
to set $\,\chi_0(g) = 1\,$ for all $\,g \in \sug$.

For $\,\la \neq 0$,\ this lift of the action to the bundle is fixed
uniquely by the demand that the connection form (cf.\
\cite[Def.\,2.2.4]{Brylinski:1993ab} or
\cite[App.\,A.3]{Woodhouse:1992de}),
\qq
\widehat A(h,z)=\tfrac{\sfi\,\sfd z}{z}+\sfi\la\,\tr\bigl(\La\,
\th_L(h)\bigr)\,,
\qqq
induced by $\,\nabla\,$ on the complement of the zero section in the
total space $\,\overline K_\la\,$ be annihilated by the vector field
\qq
\widehat L_3(h,z)=L_3(h)+\tfrac{\sfd\ }{\sfd\vep}\big
\vert_{\vep=0}\chi_\la\bigl(\erm^{\sfi\vep\si_3}\bigr)\cdot z\,
\tfrac{\p\ }{\p z}
\qqq
on $\,\sug\x\bC^\x\,$ that generates the lifted action. It is also
straightforward to check that $\,\widehat A\,$ is invariant since
$\,\sug_\la\,$ is just the maximal torus. In the case $\,\la=0$,\
horizontality and invariance with respect to the action of
$\,\sug_\la\,$ hold trivially.

Horizontality and invariance of $\,\nabla\,$ with respect to the
action of $\,\sug_\la\,$ on $\,\overline K_\la\,$ show that we can
pass to the quotient line bundle $\,K_\la = \overline K_\la /
\sug_\la\,$ with base $\,\sug / \sug_\la$.\ For $\,\la=0$,\ we
obtain the trivial line bundle over a point.

\subsection*{Line bundle with connection over
$\,\boldsymbol{Y^{[2]}\sug}$}

For the surjective submersion $\,Y\sug\,$ of \eqref{eq:cover}, the
fibred product is given by
\qq
Y^{[2]}\sug=\bigsqcup_{i,j\in\{0,1\}}\,\cO_{i,j}\,,\qquad\qquad
\cO_{i,j}=\cO_i\cap\cO_j\,.
\qqq
We need to specify a Hermitian line bundle
\qq
\pi_L\ :\ L\to Y^{[2]}\sug
\qqq
equipped with a unitary connection $\,\nabla_L\,$ of curvature
\qq\label{eq:con-from-BB}
\curv(\nabla_L)\vert_{\cO_{i,j}}=(B_j-B_i)\vert_{\cO_{i,j}}\,.
\qqq
We shall give this bundle in terms of an equivariant line bundle
$\,\widetilde L_{i,j}\,$ over the base
$\,\widetilde\cO_{i,j}=\widetilde\cO_i \cap\widetilde\cO_j$.\ The
unitary connection $\,\nabla_{\widetilde L_{i,j}}\,$ on
$\,\widetilde L_{i,j}\,$ must have the curvature
\be
\curv(\nabla_{\widetilde L_{i,j}})(\la,h) = (\widetilde
B_j-\widetilde B_i)(\la,h) = \widetilde F(\la_i-\la_j,h) = \sfi \,
\la_{i,j} \, \tr(\La\,\sfd \th_L(h))\,,
\ee
where we have abbreviated $\,\la_{i,j} = \la_j-\la_i$.\ A quick
glance at \eqref{eq:KKS-curv} shows that the KKS bundle (or rather
the equivariant bundle $\,\overline K_{\la_{i,j}}\,$ on $\,\sug$)
has the desired curvature. We therefore define
\be
\widetilde L_{i,j} = \pi_{\sug}^* \overline K_{\la_{i,j}}\,,\qquad
\text{where}\quad\pi_{\sug}\ :\ \widetilde\cO_{i,j} \to \sug\ :\
(\la,h)\mapsto
  h\,.
\ee
The result is simply the trivial bundle $\,\widetilde L_{i,j} =
\widetilde\cO_{i,j} \times \Cb \to \widetilde\cO_{i,j}\,$ with
connection $\,\nabla_{\widetilde L_{i,j}} = \sfd +
\la_{i,j}\,\tr(\La\,\th_L)$.\ For $\,i\neq j$,\ it inherits
equivariance with respect to the action of $\,\sug_{\la_{i,j}}=
\uj\,$ from $\,\overline K_{\la_{i,j}}\,$ and thus yields a
well-defined quotient bundle over $\,\Oc_{i,j}$.\ It should be
stressed that there are two distinct isotropy groups entering the
above construction, to wit, $\,\sug_{t_\la}\,$ whose action has to
be divided out (over each point $\,(\la,h)\in \widetilde\cO_{i,j}$)
when passing from $\,\widetilde\cO_{i,j}\,$ to $\,\cO_{i,j}$,\ and
$\,\sug_{\la_{i,j}}\,$ whose action on $\,\widetilde L_{i,j}\,$ is
inherited from the KKS bundle. It is the relation $\,\sug_{t_\la}
\subset \sug_{\la_{i,j}}$,\ valid for all
$\,(\la,h)\in\widetilde\cO_{i, j}$,\ that enables us to descend the
$\sug_{\la_{i,j}}$-equivariant bundle $\,\widetilde L_{i,j}\,$ to
the $\sug_{t_\la}$-quotient $\,\Oc_{i,j}\,$ of
$\,\widetilde\cO_{i,j}$.\ For $\,i=j$,\ both the connection and the
action of the isotropy groups on the fibre are trivial, and so
$\,\widetilde L_{i,i}\,$ induces the trivial bundle with trivial
connection over $\,\Oc_{i,i}$.

This completes the definition of the line bundle $\,L\to Y^{[2]}
\sug\,$ in terms of appropriate equivariant bundles $\,\widetilde
L_{i,j} \to \widetilde\Oc_{i,j}$.

\subsection*{Groupoid structure}

The construction of $\,\cG\,$ is completed by specifying the
groupoid structure on the fibres of $\,L\,$ pulled back to
\qq
Y^{[3]}\sug=\bigsqcup_{i,j,k\in\{0,1\}}\,\cO_{i,j,k}\,,\qquad\qquad
\cO_{i,j,k}=\cO_i\cap\cO_j\cap\cO_k\,.
\qqq
The map $\,\mu:\pr_{1,2}^*L \ox \pr_{2,3}^*L \to \pr_{1,3}^* L\,$
will be described in terms of an equivariant map $\,\widetilde\mu\,$
on $\,\widetilde\cO_{i,j,k}=\widetilde\cO_i\cap \widetilde\cO_j\cap
\widetilde\cO_k$.\ Namely, we set
\be
\widetilde\mu\vert_{\widetilde\cO_{i,j,k}}\ :\ \widetilde L_{i,j}
\ox \widetilde L_{j,k}\big\vert_{\widetilde\cO_{i,j,k}} \to
\widetilde L_{i,k}\big\vert_{\widetilde\cO_{i,j,k}}\,,\qquad
\widetilde\mu(i,j,k,\la,h, z \otimes z') = (i,j,k,\la,h, z \cdot
z')\,.
\ee
This map is clearly unitary and it is easy to see that it is
compatible with the connections. Equivariance with respect to the
action of the maximal torus, common to all isotropy groups, amounts
to the statement that for all $\,t\in\uj$,
\be
\widetilde\mu\bigl(i,j,k,\la,h \cdot t,\chi_{\la_{i,j}}(t) \cdot z
\otimes \chi_{\la_{j,k}}(t)\cdot z'\bigr) = \bigl(i,j,k,\la,h \cdot
t, \chi_{\la_{i,k}}(t) \cdot z \cdot  z'\bigr)\,.
\ee
This holds true by virtue of the equality
\qq
\chi_{\la_{i,j}}(t)\cdot\chi_{\la_{j,k}}(t)=
\chi_{\la_{i,k}}(t)\,,
\qqq
readily verified by direct inspection.
Once again, we are using the equivariant structures with respect to
the action the isotropy groups $\,\sug_{\la_{i,j}},\sug_{\la_{j,
k}}\,$ and $\,\sug_{\la_{i,k}}$,\ defining the three bundles
involved, to divide out the action of the isotropy group
$\,\sug_{t_\la}\,$ acting on their base. In so doing, we exploit the
relation $\,\sug_{t_\la}\subset\sug_{\la_{i,j}}\cap\sug_{\la_{j,k}}
\cap\sug_{\la_{i,k}}$,\ valid for all $\,(\la,h)\in\widetilde\cO_{i,
j,k}$.\ For the non-generic isotropy groups over $\,\la=0\,$ (in
which case $\,i=j=k=0$) and $\,\la=\sfk\,$ (in which case
$\,i=j=k=1$), the argument is exactly the same, provided that we
allow $\,t \in \sug\,$ and recall that we defined $\,\chi_0(t)=1\,$
in this case.

Finally, it is clear that $\,\mu\,$ satisfies the associativity
condition.

\section{Trivialisation over conjugacy classes}\label{sec:triv}

The next piece of the general theory that we shall need in the
subsequent discussion is the definitions of a 1-isomorphism and of a
stable isomorphism,

\bef\cite{Murray:1999ew,Waldorf:2007mm}\label{def:stab-iso} Let
$\,\cG_i=(Y_i M,B_i,L_i,\mu_{L_i}),\ i\in\{1,2\}\,$ be a pair of
gerbes of the same curvature $\,\txH\,$ over a common base $\,M$.\\
(i) A {\it
1-isomorphism} $\,\Phi_{1,2}\,$ between $\,\cG_1\,$ and $\,\cG_2$,\
denoted as
\qq
\Phi_{1,2}\ :\ \cG_1\xrightarrow{\sim}\cG_2\,,
\qqq
is a triple $\,\Phi_{1,2}=(YY_{1,2}M,E_{1,2},\a_{1,2})\,$ which
consists of
\begin{itemize}
\item a surjective submersion $\,YY_{1,2}M\xrightarrow{\pi_{YY_{1,
2}M}}Y_{1,2}M=Y_1 M\x_M Y_2 M\,$ over the fibred product
\qq
Y_1 M\x_M Y_2 M=\{\ (y_1,y_2)\in Y_1 M\x Y_2 M \ \vert \ \pi_{Y_1
M}(y_1)=\pi_{Y_2 M}(y_2) \ \}\,;
\qqq
\item a Hermitian line bundle $\,E_{1,2}\xrightarrow{\pi_{E{1,
2}}}YY_{1,2}M\,$ with connection $\,\nabla_{E_{1,2}}\,$ of curvature
\be
\curv(\nabla_{E_{1,2}})=\pi_2^*B_2-\pi_1^*B_1
\labl{eq:curv-E12}
fixed by the pullbacks of the two curvings along the maps $\,\pi_i=
\pr_i\circ\pi_{YY_{1,2}M}\,$ written in terms of the canonical
projections $\,\pr_i:Y_1 M\x_M Y_2 M\to Y_i M:(y_1,y_2)\mapsto y_i$;
\item a unitary connection-preserving isomorphism
\qq\label{eq:alpha12-iso}
\a_{1,2}\ :\ \pi_{1,3}^*L_1\ox p_2^*E_{1,2}
\xrightarrow{\sim}p_1^*E_{1,2}\ox\pi_{2,4}^*L_2\,,
\qqq
of line bundles over $\,Y^{[2]}Y_{1,2}M=YY_{1,2}M\x_M YY_{1,2}M$,\
defined in terms of the obvious maps
$\,\pi_{i,i+2}=\pr_{i,i+2}\circ(\pi_{YY_{1
,2}M}\x\pi_{YY_{1,2}M}):YY_{1,2}M\x_M YY_{1,2}M\to Y_i^{[2]}M\,$ and
$\,p_i:YY_{1,2}M\x_M YY_{1,2}M\to YY_{1,2}M$;\ the isomorphism must
be compatible with the two groupoid structures $\,\mu_{L_i}\,$ in
the sense specified by the commutativity of the diagram
\qq
\tx{\tiny{$\alxydim{@C=-.2cm@R=2.cm}{
& & \pi_{1,3}^*L_1\ox\pi_{3,5}^*L_1\ox p_3^*E_{1,2}
\ar[dll]_{\hspace{-2.5cm}\id_{\pi_{1,3}^*L_1}\ox p_{2,3}^*\a_{1,2}}
\ar[drr]^{\hspace{.5cm}\pi_{1,3,5}^*\mu_{L_1}\ox\id_{p_3^*E_{1,2}}}
& & \\ \pi_{1,3}^*L_1\ox p_2^*E_{1,2}\ox\pi_{4,6}^*L_2
\hspace{-1.cm} \ar[dr]_{p_{1,2}^*\a_{1,2}\ox\id_{\pi_{4,6}^* L_2}}
& & & & \pi_{1,5}^*L_1\ox p_3^*E_{1,2} \ar[dl]^{p_{1,3}^*\a_{1,2}} \\
& p_1^*E_{1,2}\ox\pi_{2,4}^*L_2\ox\pi_{4, 6}^*L_2
\ar[rr]^{\hspace{.5cm}\id_{p_1^*E_{1,2}}\ox\pi_{2,4,6}^*\mu_{L_2}} &
& p_1^*E_{1,2}\ox\pi_{2,6}^*L_2 & }$}}\cr\cr
\label{diag:stab-iso-comp-mu}
\qqq
of bundle isomorphisms over $\,Y^{[3]}Y_{1,2}M=YY_{1,2}M\x_M
YY_{1,2}M\x_M YY_{1,2}M$,\ the latter being endowed with the
canonical projections $\,p_i:Y^{[3]}Y_{1,2}M\to YY_{1,2}M\,$ and
$\,p_{i,j}:Y^{[3]}Y_{1,2}M\to Y^{[2]}Y_{1,2}M$,\ as well as with the
maps $\,\pi_{i,j}\,$ and $\,\pi_{i,j,k,l}\,$ given by the
corresponding canonical projections $\,\pr_{i,j}\,$ and
$\,\pr_{i,j,k,l}$ precomposed with $\,(\pi_{YY_{1,2}M})^{\times 3}$.
\end{itemize}
(ii) A {\em stable isomomorphism} is a 1-isomorphism whose
surjective submersion is given by $\,YY_{1,2}M=Y_{1,2}M\,$ with
$\,\pi_{YY_{1,2}M}=\id_{Y_{1,2}M}$.

\eef\smallskip

The more general surjective submersion allowed in the definition of
a 1-isomorphism is necessary when formulating composition below.
However, the notions of 1-isomorphism classes and stable-isomorphism
classes of gerbes are the same,

\berop{\rm\cite{Waldorf:2007mm}}\label{prop:stabiso-vs-1-iso}
Let $\,\cG_i,\ i\in\{1,2\}\,$ be two gerbes. There exists a stable
isomorphism between $\,\cG_1\,$ and $\,\cG_2\,$ if and only if there
exists a 1-isomorphism between these gerbes.
\eerop

\smallskip

Trivialisations compose a distinguished class of stable
isomorphisms. Given a gerbe $\,\cG\,$
over base $\,M$,\ they take the form
\qq
\Phi\ :\ \cG\xrightarrow{\sim}\cI(\om)\,,
\qqq
for some $\,\om\in\Om^2(M)$.\ These are special examples of a larger
family of bundle-gerbe modules \cite{Bouwknegt:2001vu} whose
definition generalises that of trivialisations in that it replaces
the notion of a stable isomorphism between a given bundle gerbe and
a trivial one with the notion of a 1-morphism, where a higher-rank
vector bundle is allowed instead of a line bundle in Definition
\ref{def:stab-iso} \cite{Waldorf:2007mm}.
\smallskip

Just as for gerbes, there are a number of natural operations on
1-isomorphisms. In the present paper, we shall only need the
pullback of a 1-isomorphism, the product and composition of
1-isomorphisms, and so we confine our presentation to these
particular operations.

\medskip\noindent{\bf Pullback:}
Let $\,\cG_i=(Y_iM,B_i,L_i,\mu_{L_i}),\ i\in\{1,2\}\,$ be a pair of
gerbes over a common base $\,M\,$ and let $\,\Phi_{1,2}:\cG_1
\xrightarrow{\sim}\cG_2\,$ be a 1-isomorphism with data
$\,\Phi_{1,2}=(YY_{1,2}M,E_{1,2},\a_{1,2})$.\ Furthermore, let
$\,N\,$ be a smooth manifold with surjective submersions $\,\pi_{Y_i
N}:Y_i N\to N$,\ and suppose that we are given a smooth map
$\,f:N\to M\,$ together with maps $\,\widehat f_i:Y_i N\to Y_i M\,$
that cover $\,f$,\ so that we obtain the pullback gerbes
$\,f^*\cG_i$.\ To define the pullback of $\,\Phi_{1,2}$,\ we need,
in addition, a surjective submersion $\,\pi_{YY_{1,2}N}:Y
Y_{1,2}N\to Y_{1,2}N\equiv Y_1 N\x_N Y_2 N\,$ together with a map
$\,\widetriangle f_{1,2}:YY_{1,2}N\to YY_{1,2} M\,$ that covers
$\,\widehat f_{1,2}=\widehat f_1\x\widehat f_2$,\ as expressed by
the commutative diagram
\qq
\alxydim{@C=5em@R=3em}{YY_{1,2}N \ar[d]_{\pi_{YY_{1,2}N}}
\ar[r]^{\widetriangle f_{1,2}} & YY_{1,2}M \ar[d]^{\pi_{YY_{1,2}M}}
\\ Y_{1,2}N \ar[r]_{\widehat f_{1,2}} & Y_{1,2}M}\,.
\qqq
Given these, the {\em pullback of} $\,\Phi_{1,2}\,$ {\em to} $\,N\,$
{\em along} $\,f\,$ is the 1-isomorphism
\qq
f^*\Phi_{1,2}=(YY_{1,2}N,\widetriangle f_{1,2}^*E_{1,2},
\widetriangle f_{1,2}^{[2]*}\a_{1,2})\,,\qquad\qquad f^*\Phi_{1,2}\
:\ f^*\cG_1 \xrightarrow{\sim}f^*\cG_2\,,
\qqq
where we use the maps $\,\widetriangle f_{1,2}^{[n]}:Y^{[n]}
Y_{1,2}N\to Y^{[n]}Y_{1,2}M\,$ induced by $\,\widetriangle f_{1,
2}$.

\medskip\noindent{\bf Product:}
Consider a quadruple of gerbes $\,\cG_i=(Y_i M,B_i,L_i,\mu_{L_i}),\
i\in\{1,2,3,4\}\,$ over a common base $\,M$,\ with products
$\,\cG_1\star\cG_3\,$ and $\,\cG_2\star\cG_4\,$ as in
\eqref{eq:gerbe-prod}, and a pair of 1-isomorphisms $\,\Phi_{1,2}:
\cG_1\xrightarrow{\sim}\cG_2,\ \Phi_{3,4}:
\cG_3\xrightarrow{\sim}\cG_4\,$ with data
$\,\Phi_{i,j}=(YY_{i,j}M,E_{i,j},\a_{i,j})$.\ The latter permit to
define a {\em product 1-isomorphism}
$\,\Phi_{1,2}\star\Phi_{3,4}:\cG_1
\star\cG_3\xrightarrow{\sim}\cG_2\star\cG_4\,$ as
\qq
\Phi_{1,2}\star\Phi_{3,4}=(YY_{1,2}M\x_M YY_{3,4}M,\pr_1^*E_{1,2}\ox
\pr_2^*E_{3,4},\pr_{1,3}^*\a_{1,2}\ox\pr_{2,4}^*\a_{3,4})\,,
\qqq
where the surjective submersion $\,\pi_{Y(Y_{1,3}M\x_M Y_{2,4}M)}:
Y(Y_{1,3}M\x_M Y_{2,4}M)\to Y_{1,3}M\x_M Y_{2,4}M\,$ for $\,Y(Y_{1,
3}M\x_M Y_{2,4}M)=YY_{1,2}M\x_M YY_{3,4}M\,$ is given by the map
$\,\pi_{Y(Y_{1,3}M\x_M Y_{2,4}M)}=\tau_{2,3}\circ(\pi_{YY_{1,2}M}\x
\pi_{YY_{3,4}M})\,$ expressed in terms of the transposition map
$\,\tau_{2,3}:Y_1 M\x_M Y_2 M\x_M Y_3 M\x_M Y_4 M\to Y_1 M\x_M Y_3
M\x_M Y_2 M\x_M Y_4 M:(y_1,y_2,y_3,y_4)\mapsto(y_1,y_3,y_2,y_4)$,\
and where $\,\pr_i:YY_{1,2}M\x_M YY_{3,4}M\to YY_{2i-1,2i}M\,$ and
$\,\pr_{i,i+2}:YY_{1,2}M\x_M YY_{3,4}M\x_M YY_{1,2}M\x_M YY_{3,4}M
\to Y^{[2]}Y_{2i-1,2i}M\,$ are the canonical projections.

\medskip\noindent{\bf Composition:}
For a triple $\,\cG_i=(Y_i M,B_i,L_i,\mu_{L_i}),\ i\in\{1,2 ,3\}\,$
of gerbes over a common base $\,M$,\ and a pair $\,\Phi_{i,j}=
(YY_{i,j}M,E_{i,j},\a_{i,j}),\ (i,j)\in\{(1,2),(2,3)\}\,$ of
1-isomorphisms $\,\Phi_{i,j}:\cG_i\xrightarrow{\sim}\cG_j$,\ the
{\it composition of 1-isomorphisms} $\,\Phi_{2,3}\,$ {\it and}
$\,\Phi_{1,2}\,$ is the 1-isomorphism
\qq
\Phi_{2,3}\circ\Phi_{1,2}=(YY_{1,3}M,E_{1,2,3},\a_{1,2,3})\ :\
\cG_1\xrightarrow{\sim}\cG_3
\qqq
with the surjective submersion
\qq\label{eq:can-surj-sub-comp}
(\pr_1,\pr_3)\circ(\pi_{YY_{1,2}M}\x\pi_{YY_{2,3}M})\ :\ YY_{1,3}
M=YY_{1,2}M\x_{Y_2 M}YY_{2,3}M\to Y_1 M\x_M Y_3 M\,,
\qqq
written in terms of the canonical projections $\,\pr_i:Y_1 M\x_M Y_2
M\x_M Y_3 M\to Y_i M$,\ with the Hermitian line bundle
\qq
E_{1,2,3}=p_1^*E_{1,2}\ox p_2^*E_{2,3}\to YY_{1,3}M
\qqq
with the connection
\qq
\nabla_{E_{1,2,3}}=p_1^*\nabla_{E_{1,2}}+p_2^*\nabla_{E_{2,3}}\,,
\qqq
both written in terms of the canonical projections $\,p_i:YY_{1,3}M
\to YY_{i,i+1}M$,\ and with the composite bundle isomorphism
\qq
\a_{1,2,3}:=\bigl(\id_{p_1^*E_{1,2}}\ox p_{2,4}^*\a_{2,3}\bigr)
\circ\bigl(p_{1,3}^*\a_{1,2}\ox\id_{p_4^*E_{2,3}}\bigr)\,,
\qqq
defined in terms of the canonical projections $\,p_{i,i+2}:
Y^{[2]}Y_{1,3}M\to Y^{[2]}Y_{i,i+1}M\,$ and $\,p_1=\pr_1\circ
p_{1,3},\ p_4=\pr_2\circ p_{2,4}$,\ where $\,\pr_i:Y^{[2]}Y_{i,i+1}M
\to Y_{i,i+1}M\,$ are, again, the canonical projections.\smallskip

A useful description of gerbes with a given curvature is provided by
the following statement \cite{Gawedzki:1987ak}, which derives from
the relation between bundle gerbes and the Deligne hypercohomology,
and from the relation of the latter to sheaf cohomology
\cite{Brylinski:1993ab}, cf., e.g.,
\cite[Sect.\,2.3]{Gawedzki:2002se} and \cite[Prop.\,2.2]{Gomi:2003}.

\berop\label{prop:stab-iso-cl}
The set of stable-isomorphism classes of gerbes with a given
curvature over base $\,M\,$ is a torsor under the action of the
cohomology group $\,H^2\bigl(M,\uj\bigr)$.
\eerop

In particular, since the cohomology group
$\,H^2\bigl(\sug,\uj\bigr)\,$ is trivial, this proposition shows
that the canonical gerbe is unique up to a 1-isomorphism.

\begin{Rem}
Stable isomorphisms play a central r\^ole in the construction of
two-dimensional non-linear $\si$-models. The topological term of the
action functional is given as (the logarithm of) the surface
holonomy of a trivialisation of the pullback, along the embedding
map $\,X:\Si\to M$,\ of the gerbe over the target space of the
$\si$-model to the two-dimensional world-sheet $\,\Si\,$
\cite{Gawedzki:1987ak,Gawedzki:2002se}. From the point of view of
this paper, it is also important that a stable isomorphism together
with a definition of the world-volume and of the curvature 2-form of
the so-called bi-brane describe conformal defects separating phases
of the two-dimensional field theory
\cite{Fuchs:2007fw,Runkel:2008gr,Sarkissian:2008dq} (cf.\ also
\cite{Carey:2002,Gawedzki:2002se} for the boundary case).
\\
A particular instantiation of this latter fact occurs on
world-sheets with a boundary, capturing the dynamics of the open
string. Here, the $\sigma$-model action functional is defined for a
submanifold $\,\xcD\,$ of $\,M$,\ embedded in the target space as
$\,\iota:\xcD\emb M\,$ and termed the D-brane (or $\cG$-brane)
world-volume, alongside a trivialisation of the bulk\footnote{In the
context of two-dimensional non-linear $\si$-models, the term `bulk'
refers to the conformally invariant field theory defined away from
the boundary of the world-sheet, describing the dynamics of the
closed string.} gerbe over $\,\xcD$,\
$\,\Phi:\iota^*\cG\xrightarrow{\sim} \cI(\om)$.\ The trivialisation
is defined in terms of a trivial gerbe $\,\cI(\om)\,$ with a
globally defined curving $\,\om\in\Om^2(\xcD)\,$ satisfying
$\,\iota^*\curv(\cG)=\sfd\om\,$ and called the D-brane curvature.
\\
In the WZW model,\ there is a distinguished class of D-branes -- the
so-called (untwisted) maximally symmetric D-branes -- whose
world-volume is given by conjugacy classes. Below, we review the
construction of the associated 1-isomorphism in the case of
$\,\sug$.

\end{Rem}

We shall denote the conjugacy class of the element $\,t_\la\,$ of
the maximal torus as
\qq
\cC_\la=\{\ \Ad_g(t_\la) \ \vert \ g\in\sug \ \}\,.
\qqq
The restriction of the Cartan 3-form to $\,\cC_\la\,$ admits a
global primitive $\,\txH\vert_{\cC_\la}=\sfd\om_\la$.\ In the equivariant
formulation on $\,[0,\sfk]\x\sug$,\ these objects are given by
\qq\label{eq:om_la-def}
\widetilde\cC_\la = c^{-1}(\cC_\la) = \{ \la \} \times \sug\,,\qquad
\qquad \widetilde\omega_\la(h) = \widetilde Q(\la,h)|_{\widetilde\cC_\la}\,.
\qqq
In the remainder of this section, we describe -- after
\cite{Gawedzki:2002se} -- the conditions under which there exists a
trivialisation of the restricted gerbe
\qq
\Phi_\la\ :\ \cG\vert_{\cC_\la}\xrightarrow{\sim}\cI(\om_\la)
\qqq
and give the details of its construction.

\subsection*{Surjective submersion and curvature 2-forms}

We need to give a surjective submersion
$\,YY_{1,2}\cC_\la\xrightarrow{\pi_{YY_{1,2}\cC_\la}} Y_1\cC_\la\x_M
Y_2\cC_\la$,\ where $\,Y_1\cC_\la\,$ and $\,Y_2\cC_\la\,$ are the
surjective submersions of the gerbes $\,\cG|_{\cC_\la}\,$ and
$\,\cI(\om_\la)$,\ respectively.

For the restricted canonical gerbe $\,\cG\vert_{\cC_\la}$,\ we
have the surjective submersion
\qq
Y\cC_\la=\cC_{\la;0}\sqcup\cC_{\la;1}\xrightarrow{\pi_{Y\cC_\la}}
\cC_\la \qquad \text{with} \quad \cC_{\la;l}=\cC_\la\cap\cO_l\,.
\qqq
Thus, $\,\cC_{\la;l}=\cC_\la\,$ except for the special cases
$\,\cC_{0;1}=\emptyset=\cC_{\sfk;0}$,\ which are simply dropped,
leaving us with $\,Y\cC_0=\cC_{0,0}=\{e\}\,$ and $\,Y\cC_\sfk=
\cC_{\sfk;1}=\{-e\}$.\ For the trivial gerbe $\,\cI(\om_\la)$,\ we
have the trivial surjective submersion $\,\id_{\cC_\la}:\cC_\la \to
\cC_\la$.\ The product of the two surjective submersions fibred over
the common base $\,\cC_\la\,$ of the two gerbes is given by
\qq
Y\cC_\la\x_{\cC_\la}\cC_\la=Y\cC_\la\,.
\qqq
Finally, for $\,YY_{1,2}\cC_\la\to Y\cC_\la$,\ we choose the trivial
surjective submersion of $\,Y\cC_\la$,\ i.e.\ $\,YY_{1,2}\cC_\la =
Y\cC_\la\,$ and $\,\pi_{YY_{1,2}\cC_\la} = \id_{Y\cC_\la}$.

In order to give the line bundle over $\,Y\cC_\la$,\ we shall,
again, work in the equivariant formulation via the pullback by
$\,c$.\ We therefore define
\qq
Y\widetilde\cC_\la=\widetilde\cC_{\la;0}\sqcup\widetilde\cC_{\la;1}\,,
\quad \text{with} \quad \widetilde \cC_{\la;l}= \widetilde
\cC_\la\cap \widetilde\cO_l\,.
\qqq
The difference of the pullback curvings in \eqref{eq:curv-E12} on
$\,\widetilde \cC_{\la;l} \subset Y\widetilde\cC_\la\,$ then reads
\be
\widetilde\omega_\la(h) - \widetilde B_l(\la,h) = - \widetilde F(\la
- \la_l,h) = \sfi\, (\la - \la_l)\, \tr\big(\Lambda\, \sfd
\th_L(h)\big)\,.
\labl{eq:Yla-curv}
The latter 2-form can be the curvature 2-form of a bundle only if
its periods over 2-cycles of $\,\widetilde \cC_{\la;l}\,$ take
values in $\,2\pi\bZ$.\ For $\,\la=0\,$ or $\,\la=\sfk$,\ the
manifold $\,Y\cC_\la\,$ is simply a point, and so the line bundle
exists -- it is the trivial bundle. For $\,\la \not\in \{0,\sfk\}$,\
the conjugacy classes $\,\cC_{\la; l}\,$ with $\,l=0,1\,$ provide a
choice of representatives of the generators of $\,H_2(Y\cC_\la)$.\
The conjugacy classes can be parametrised by mapping
$\,[0,\pi[\times [0,\pi/2]\to \widetilde \cC_{\la;l} \to
\cC_{\la;l}\,$ as $\,(\varphi,\theta)\mapsto\bigl(
\la,h(\varphi,\theta)\bigr)\mapsto h(\varphi,\theta)\cdot t_\la\cdot
h(\varphi,\theta)^{-1}$,\ where $\,h(\varphi,\theta) = \erm^{\sfi
\varphi\si_3}\cdot\erm^{\sfi\th\si_2}$.\ The pullback of $\,-F\,$ is
then given by
\be
-\widetilde F(\varphi,\theta)= -2(\la-\la_l)\,\sin
2\th\,\sfd\varphi\wedge\sfd\th\,,
\ee
so that
\qq\label{eq:Bdiff-int}
\int_{\cC_{\la;l}}\,(-F) = 2 (\la_l-\la)\, \int_0^{\pi}\, d\varphi\,
\int_0^{\pi/2}\, \hspace{-1em} d\theta \, \sin 2\th = 2\pi\, (\la_l
-\la)\,.
\qqq
The result lies in $\,2 \pi \bZ\,$ if and only if $\,\la\in\bZ$,\ or
-- equivalently -- if and only if $\,\la\,\La\in\cA_{\rm W}^\sfk\,$
is an integral weight.

\subsection*{Line bundle over $\,\boldsymbol{Y\cC_\la}$}

The calculation of periods shows that the line bundle over
$\,Y\cC_\la\,$ with the curvature \eqref{eq:Yla-curv} exists if and
only if $\,\la \in \faff\,$ (so that $\,\la\,\La\in\cA_{\rm
W}^\sfk$). In order to construct this bundle -- or rather the
equivariant line bundle $\,\widetilde E_\la \rightarrow
Y\widetilde\cC_\la\,$ -- we can use the KKS bundle again. Namely, we
set
\be
\widetilde E_{\la;l}=\pi_{\sug}^*\overline K_{\la-\la_l}
\,,\qquad\tx{where} \quad\pi_{\sug}\ :\
\widetilde\cC_{\la;l}\to\sug\ :\ (\la,h)\mapsto h\,.
\ee
For $\,\la\not\in\{0,\sfk\}$,\ it inherits equivariance with respect
to the action of $\,\sug_{\la-\la_l}=\uj \,$ from $\,\overline
K_{\la-\la_l}$,\ and hence -- as $\,\sug_{t_\la}=\sug_{\la-\la_l}\,$
-- it descends to a well-defined quotient bundle
$\,E_{\la;l}\to\cC_{\la;l}$.\ For $\,\la =0$,\ we have $\,l=0\,$ and
both the connection and the action of the isotropy group are
trivial, and so $\,\widetilde E_{0;0}\,$ induces the trivial bundle
over $\,\cC_{0;0} = \{ e \}$.\ Similarly for $\,\la = \sfk$,\ we get
the trivial bundle over $\cC_{\sfk;1} = \{-e\}$.\ Altogether, the
quotient bundles $\,E_{\la;l}\to\cC_{\la; l}\,$ compose a line
bundle $\,E_\la\to Y\cC_\la$.

\subsection*{Isomorphism of pullback bundles}

The last piece of data that we need is the isomorphism
\eqref{eq:alpha12-iso} of the pullback bundles. The bundles are
pulled back to $\,YY_{1,2}\cC_\la \times_{\cC_\la} YY_{1,2}\cC_\la =
Y\cC_\la \times_{\cC_\la} Y\cC_\la = Y^{[2]}\cC_\la$,\ where
\qq
Y^{[2]}\cC_\la=\bigsqcup_{i,j\in\{0,1\}}\,\cC_{\la;i,j}\,,
\qquad\qquad\cC_{\la;i,j}=\cC_\la\cap\cO_{i,j}\,.
\qqq
As before, we denote by $\,\pr_i:Y^{[2]}\cC_\la\to Y\cC_\la\,$ the
projections onto the two factors. In the present case, the
isomorphism \eqref{eq:alpha12-iso} boils down to
\qq
\a_\la\ :\ L\vert_{Y^{[2]}\cC_\la}\ox\pr_2^*E_\la\to\pr_1^*E_\la\,.
\qqq
In the equivariant formulation, we work with the subsets
$\,\widetilde\cC_{\la;i,j}=\widetilde\cC_\la\cap\widetilde\cO_{i,
j}$.\ We denote points in $\,\widetilde\cC_{\la;i,j}\,$ as
$\,(i,j,\la,h)$,\ with $\,h \in \sug$.\ The isomorphism takes the
simple form
\be
  \widetilde\a_\la|_{\widetilde\cC_{\la;i,j}}\ :\
  \widetilde L_{i,j} \otimes \widetilde E_{\la,j}\big\vert_{\widetilde\cC_{\la;i,j}}
  \to \widetilde E_{\la,i}\big\vert_{\widetilde\cC_{\la;i,j}}\ :\ (i,j,\la,h,z \otimes z') \mapsto
  (i,j,\la,h,z \cdot z')\,.
\ee
This map is clearly unitary and it is easy to see that it is
compatible with the connections. Equivariance with respect to the
maximal torus, common to all isotropy groups, amounts to the
statement that, for all $\,s\in\uj$,
\be
  \widetilde\a_\la\bigl(i,j,\la,h \cdot s,\chi_{\la_{i,j}}(s) \cdot z \otimes \chi_{\la-\la_j}(s)\cdot z'\bigr)
  = \bigl(i,j,\la,h \cdot s, \chi_{\la-\la_i}(s) \cdot z \cdot  z'\bigr)\,.
\ee
This holds true because of the equality
\qq
\chi_{\la_{i,j}}(s) \cdot \chi_{\la-\la_j}(s) =
\chi_{\la-\la_i}(s)\,.
\qqq
For the non-generic isotropy groups over $\,\la =0\,$ (in which case
$\,i=j=0$) and $\,\la=\sfk\,$ (in which case $\,i=j=1$), the
argument is the same, only that now $\,s \in \sug\,$ and
$\,\chi_0(s)=1$.\ We conclude that $\,\widetilde\a_\la\,$ yields a
bundle isomorphism $\,\a_\la$.

Finally, it is clear that $\,\a_\la\,$ is also compatible with the
groupoid structures on $\,\cG \vert_{\cC_\la}\,$ and
$\,\cI(\om_\la)$.

\section{Elements of the multiplicative structure}\label{sec:el-mult}

A multiplicative structure is a 1-isomorphism $\,\cM\,$ between
gerbes on $\,\Gx\x\Gx$,\ together with a 2-isomorphism subject to
some coherence properties \cite{Carey:2004xt,Waldorf:2008mult}.
Here, $\,\Gx\,$ is a Lie group, and the gerbes in question are
obtained by pulling back a power of the basic gerbe on $\,\Gx\,$
along the projection maps to each factor, and along the
multiplication map, respectively. If $\,\Gx\,$ is compact simple and
connected, a multiplicative structure exists under certain
conditions on the level \cite{Gawedzki:2009jj}, and it exists for
all levels if $\,\Gx\,$ is simply connected \cite{Waldorf:2008mult}.
Whenever a multiplicative structure exists, it is unique up to an
isomorphism \cite{Waldorf:2008mult,Gawedzki:2009jj}. In the present
paper, we shall only need the 1-isomorphism $\,\cM\,$ restricted to
certain submanifolds of $\,\sug \x \sug$.

\subsection*{2-isomorphisms}

The last element of the general theory necessary to fully understand
the subsequent discussion is the notion of a 2-isomorphism between
1-isomorphisms. We begin with
\bef\cite{Stevenson:2000wj,Waldorf:2007mm}\label{def:2-iso}
Let $\,\cG_i=(Y_i M,B_i,L_i,\mu_{L_i}),\ i\in\{1,2\}\,$ be a pair of
gerbes over a common base $\,M\,$ and of the same curvature, and let
$\,\Phi_{1,2}^A=(Y^A Y_{1,2}M,E_{1,2}^A,\a_{1, 2}^A),\
A\in\{1,2\}\,$ be a pair of 1-isomorphisms between $\,\cG_1\,$ and
$\,\cG_2$,
\qq
\Phi^A_{1,2}\ :\ \cG_1\xrightarrow{\sim}\cG_2\,.
\qqq
A {\it 2-isomorphism} $\,\varphi\,$ between $\,\Phi_{1,2}^1\,$ and
$\,\Phi_{1,2}^2$,\ denoted as
\qq
\varphi\ :\ \Phi_{1,2}^1\Longrightarrow\Phi_{1,2}^2\,,
\qqq
is a pair $\,(YY^{1,2}Y_{1,2}M,\varphi)\,$ determined by a choice of
a surjective submersion
\qq
YY^{1,2}Y_{1,2}M\xrightarrow{\pi_{YY^{1,2}Y_{1,2}M}}Y^{1,2}Y_{1,2}M
\,,\qquad\qquad Y^{1,2}Y_{1,2}M=Y^1 Y_{1,2}M\x_{Y_{1,2}M}Y^2 Y_{1,2}
M\,,
\qqq
together with a unitary connection-preserving isomorphism
\qq
\varphi\ :\ p_1^*E_{1,2}^1\xrightarrow{\sim}p_2^*E_{1,2}^2
\qqq
between the two line bundles over $\,YY^{1 ,2}Y_{1,2}M$,\ the latter
coming with the maps $\,p_A=\pr_A\circ
\pi_{YY^{1,2}Y_{1,2}M}:YY^{1,2}Y_{1,2}M\to Y^A Y_{1,2}M\,$ written
in terms of the canonical projections $\,\pr_A:Y^{1,2}Y_{1,2}M\to
Y^A Y_{1,2}M$.\ The isomorphism must be compatible with the two
isomorphisms $\,\a_{1,2}^A\,$ in the sense specified by the
commutativity of the diagram
\qq
\alxydim{@C=3cm@R=2cm}{
p_{1,3}^*L_1\ox\pi_3^*E_{1,2}^1
\ar[d]_{\id_{p_{1,3}^*L_1}\ox\pr_2^*\varphi}
\ar[r]^{\pi_{1,3}^*\a_{1,2}^1} & \pi_1^*E_{1,2}^1\ox p_{2,4}^*
L_2 \ar[d]^{\pr_1^*\varphi\ox\id_{p_{2,4}^*L_2}} \\
p_{1,3}^* L_1\ox\pi_4^*E_{1,2}^2 \ar[r]_{\pi_{2,4}^*\a_{1,2}^2}
& \pi_2^*E_{1,2}^2\ox p_{2,4}^*L_2}
\qqq
of bundle isomorphisms over $\,Y^{[2]}Y^{1,2}Y_{1,2}M=YY^{1,2}Y_{1,
2}M\x_M YY^{1,2}Y_{1,2}M$,\ the latter coming with the following
maps:
\begin{itemize}
\item[-] the canonical projections $\,\pr_i:Y^{[2]}Y^{1,2}Y_{1,2}M
\to YY^{1,2}Y_{1,2}M$;
\item[-] the natural maps
$\,p_{i,i+2}=\pr_{i,i+2}\circ(\pi_{Y^A Y_{1,2}M}\x\pi_{Y^A
Y_{1,2}M})\circ(p_A\x p_A)$,\ independent of $\,A\,$ and defined in
terms of $\,p_A\,$ as above, alongside the canonical projections
$\,\pr_{i,i+2}:Y_{1,2}M\x_M Y_{1,2}M\to Y_i^{[2]}M$;\
\item[-] the natural maps
$\,\pi_{2i-1}:Y^{[2]}Y^{1,2}Y_{1,2}M
\xrightarrow{\pr_i}YY^{1,2}Y_{1,2}M\xrightarrow{p_1}Y^1Y_{1,2}M\,$
and
$\,\pi_{2i}:Y^{[2]}Y^{1,2}Y_{1,2}M$ $\xrightarrow{\pr_i}YY^{1,2}Y_{1,
2}M\xrightarrow{p_2}Y^2Y_{1,2}M$,\ alongside
$\,\pi_{1,3}:Y^{[2]}Y^{1,2}Y_{1,2}M\to Y^1 Y_{1,2}M\x_M Y^1
Y_{1,2}M\,$ and $\,\pi_{2,4}:Y^{[2]}Y^{1,2}Y_{1,2}M \to Y^2
Y_{1,2}M\x_M Y^2 Y_{1,2}M$.
\end{itemize}
\eef

The notion of equivalence of 2-isomorphisms is specified in the
following definition.

\bef\cite{Waldorf:2007mm}\label{def:2-iso-equiv}
Let $\,(Y_i Y^{1,2}Y_{1,2}M,\varphi_i),\ i\in\{1,2\}\,$ be two
2-isomorphisms between a pair $\,\Phi_{1,2}^A=(Y^A Y_{1,2}M,E_{1,
2}^A,\a_{1, 2}^A),\ A\in\{1,2\}\,$ of 1-isomorphisms between two
gerbes $\,\cG_i= (Y_i M,B_i,L_i,\mu_{L_i}),\ i\in\{1,2\}\,$ over a
common base $\,M\,$ and of the same curvature. We call the
2-isomorphisms {\em equivalent} if there exists a smooth space
$\,Y_{1,2}Y^{1,2}Y_{1,2}M\,$ and a pair of surjective submersions
$\,\pi_{Y_{1,2}Y^{1,2}Y_{1,2}M}^i:Y_{1,2}Y^{1,2}Y_{1,2}M\to Y_i Y^{1
,2}Y_{1,2}M\,$ satisfying the identity
\qq
\pi_{Y_1 Y^{1,2}Y_{1,2}M} \circ \pi_{Y_{1,2}Y^{1,2}Y_{1,2}M}^1 =
\pi_{Y_2 Y^{1,2}Y_{1,2}M} \circ \pi_{Y_{1,2}Y^{1,2}Y_{1,2}M}^2
\qqq
and such that
\qq
\pi_{Y_{1,2}Y^{1,2}Y_{1,2}M}^{1*}\varphi_1=\pi_{Y_{1,2}Y^{1,2}Y_{1,
2}M}^{2*}\varphi_2\,.
\qqq
\eef
\bigskip

We have
\berop\label{prop:2-iso-triv-subm}
Every 2-isomorphism between a pair $\,\Phi_{1,2}^A=(Y^A Y_{1,2}M,
E_{1,2}^A,\a_{1, 2}^A),\ A\in\{1,2\}\,$ of 1-isomorphisms between
two gerbes $\,\cG_i= (Y_i M,B_i,L_i,\mu_{L_i}),\ i\in\{1,2\}\,$ over
a common base $\,M\,$ and of the same curvature is equivalent to a
2-isomorphism between these 1-isomorphisms with the trivial
surjective submersion $\,\id_{Y^{1,2}Y_{1,2}M}$.
\eerop
The proposition follows from the proof of Theorem 1 in
\cite{Waldorf:2007mm} (there, the statement is made for
2-isomorphisms between stable isomorphisms but it can be
generalised).

\begin{Rem}
Gerbes over a given base $\,M$,\ together with 1-isomorphisms and
2-isomorphisms, form a 2-category, introduced in
\cite{Stevenson:2000wj} and further explored in
\cite{Waldorf:2007mm}. In this sense, 2-isomorphisms are a logical
completion of the family of structures discussed in the preceding
sections. They are essential constituents of (twisted) equivariant
structures on bundle gerbes
\cite{Gawedzki:2002se,Gawedzki:2003pm,Schreiber:2005mi,Gawedzki:2007uz,Gawedzki:2008um}
and play a fundamental r\^ole in the construction of $\si$-models on
world-sheets with (intersecting) defect lines \cite{Runkel:2008gr}.
They are also part of the definition of the multiplicative structure
on gerbes over Lie groups \cite{Carey:2004xt,Waldorf:2008mult},
which we shall explain next to the extent strictly necessary for
understanding subsequent considerations.
\end{Rem}

There is a convenient cohomological classification of stable
isomorphisms between a given pair of gerbes; it can be proved
similarly to Proposition \ref{prop:stab-iso-cl}.

\berop\label{prop:2-iso-cl}
The set of 2-isomorphism classes of stable isomorphisms between two
given gerbes over a common base $\,M\,$ is a torsor under the
cohomology group $\,H^1\bigl(M,\uj\bigr)$.
\eerop

\subsection*{The submanifolds $\,\boldsymbol{\Tlmn}$}

Consider the triple of smooth maps
\qq\label{eq:pi-m-def}
\txp_i\ &:\ \sug\x\sug\to\sug\ &:\ (g_1,g_2)\mapsto g_i\,,\qquad
i\in\{1,2\}\,,\cr\cr \txm\ &:\ \sug\x\sug\to\sug\ &:\
(g_1,g_2)\mapsto g_1\cdot g_2\,.
\qqq
The submanifolds
\qq\label{eq:T_lmn-def}
\Tlmn=\txp_1^{-1}(\cC_\la)\cap\txp_2^{-1}(\cC_\mu)\cap
\txm^{-1}(\cC_\nu) \ \subset \ \sug \x \sug
\qqq
will play a central r\^ole in the following. The subset of the
parameter space for which these submanifolds are non-empty is called
the {\em fusion polytope} of $\,\sug$,
\be
  \xcF = \{\ (\la,\mu,\nu) \in [0,\sfk]^{\times 3} \ | \ \Tlmn
  \neq \emptyset \ \}\,.
\labl{eq:setF-maindef}

\begin{Rem}
Using the three maps $\,\txp_1,\txp_2,\txm:\sug\x\sug\to\sug$,\ we
may pull the various geometric objects introduced previously, such
as the canonical gerbe $\,\cG\,$ and its trivialisations
$\,\Phi_\la, \Phi_\mu\,$ and $\,\Phi_\nu$,\ back to $\,\Tlmn$.\ We
may also restrict the multiplicative structure on $\sug \x \sug$ to
it. The motivation to do so ultimately derives from the study of the
maximally symmetric bi-brane of the WZW model for $\,\sug\,$
\cite{Fuchs:2007fw} and will be expounded in full detail in the
companion paper \cite{Runkel:2009}. As already pointed out in
\cite{Fuchs:2007fw}, the space $\,\Tlmn\,$ is closely related to the
moduli space of flat connections on a principal $\sug$-bundle over a
Riemann sphere with three punctures, with the holonomy around the
punctures constrained to take values in the three conjugacy classes
$\,\cC_\la,\cC_\mu\,$ and $\,\cC_\nu$,\ as reconstructed in
\cite{Alekseev:1993rj}. The moduli space, in turn, reproduces, upon
quantisation \cite{Witten:1988hf}, the space of conformal blocks of
the WZW model for $\,\sug\,$ with insertions of primary fields from
the corresponding representations of highest weights $\,\la,\mu\,$
and $\,\nu\,$ of the affine Lie algebra $\,\suk$.\ Our point of view
is slightly different in that we are trying to recover the fusion
rules from the existence of three-fold junctions of defect lines in
the classical WZW model \cite{Runkel:2008gr} rather than from the
quantisation of a moduli space. In the quantum WZW model, the two
approaches to fusion rules are equivalent, as such defect junctions
exist if and only if the corresponding space of conformal
three-point blocks is nonzero \cite{Frohlich:2006ch,Runkel:2008gr}.
\end{Rem}

An equivalent way of writing the submanifolds $\,\Tlmn\,$ is
\qq\label{eq:T-param}
\Tlmn=\bigl\{\ \bigl(\Ad_x(t_\la),\Ad_{x\cdot a}(t_\mu)\bigr) \
\vert \ x, a \in\sug\quad \text{and}\quad t_\la\cdot \Ad_a (t_\mu)
\in \cC_\nu \ \bigr\}\,.
\qqq
Here, $\,a\,$ is taken to run over all of $\,\sug$,\ but it is easy
to see that it suffices to pick
a representative $\,a\,$ from each of the equivalence classes
\qq
[a] \in \sug_{t_\la}\backslash\sug/\sug_{t_\mu}\,.
\qqq
Recalling the parametrisation \eqref{eq:Euler-param} by the Euler
angles, we infer that we can choose this representative in the form
$\,a = \erm^{\sfi \th \si_2}\,$ and let $\,\th\,$ vary over those
values in $\,[0,\pi/2]\,$ which are allowed by the condition
$\,t_\la\cdot \Ad_a(t_\mu) \in \cC_\nu$.\ In fact, the set of such
$\,\th\,$ is either empty or contains a single element. In order to
see this, note, first, that
\be
\erm^{\sfi\th\si_2}\cdot t_\mu\cdot\erm^{-\sfi\th\si_2}=
\begin{pmatrix}
\cos\tfrac{\pi\mu}{\sfk}+\sfi\,\sin\tfrac{\pi\mu}{\sfk}\,\cos 2\th
& -\sfi\,\sin\tfrac{\pi\mu}{\sfk}\,\sin 2\th \\
-\sfi\,\sin\tfrac{\pi\mu}{\sfk}\,\sin 2\th & \cos\tfrac{\pi
\mu}{\sfk}-\sfi\,\sin\tfrac{\pi\mu}{\sfk}\,\cos 2\th
\end{pmatrix}\,.
\ee
For $\,\sug$,\ the condition $\,t_\la\cdot\Ad_a(t_\mu)\in\cC_\nu\,$
is equivalent to $\,\tr\bigl(t_\la\cdot\Ad_a(t_\mu)\bigr)=\tr(
t_\nu)$.\ Explicitly, this condition reads
\qq\label{eq:class-cos-range}
\cos\bigl(\tfrac{\pi\,(\la+\mu
)}{\sfk}\bigr)+(\sin\th)^2\,\bigl[\cos\bigl(\tfrac{\pi\,(\la-\mu
)}{\sfk}\bigr)-\cos\bigl(\tfrac{\pi\,(\la+\mu)}{\sfk}\bigr)\bigr] =
\cos\bigl(\tfrac{\pi\nu}{\sfk}\bigr)
\qqq
and has a solution for $\,\th\,$ if and only if
$\,\vert\la-\mu\vert\leq\nu \leq\min(\la+\mu,2\sfk-\la-\mu)$.\
Altogether, we have proved (cf., e.g.,
\cite[Prop.\,3.1]{Jeffrey:1991rp}):

\begin{Lem}
The fusion polytope of $\,\sug\,$ is given by
\qq\label{eq:F-alternate}
\xcF=\bigl\{\ (\la,\mu,\nu) \in [0,\sfk]^{\x3} \ \big|\
\vert\la-\mu\vert\leq\nu\leq\min(\la+\mu,2\sfk-\la-\mu) \ \big\}\,.
\qqq
\end{Lem}

For later use, we also define the following subsets of the fusion
polytope: the boundary $\,\p\xcF$,\ the corners $\,\partial_0
\xcF$,\ the boundary less the corners $\,\partial_{12} \xcF$,\ the
interior $\,\mathring \xcF$,\ and the fusion polytope less its
corners and edges $\,\dot\xcF$.\ Altogether,
\be\begin{array}{rl}\displaystyle
  \partial\xcF \etb=
    \bigl\{\ (\la,\mu,\nu) \in [0,\sfk]^{\x3} \ \big|\ \nu=\vert\la-\mu\vert
    \quad \text{or}\quad  \nu=\min(\la+\mu,2\sfk-\la-\mu) \ \big\}\,,\enl
  \partial_0 \xcF \etb=
    \big\{\ (\la,\mu,\nu) \in \xcF \ \big|\ \la,\mu,\nu \in \{0,\sfk\} \ \big\}\,,\enl
  \partial_{12} \xcF \etb=
    \partial\xcF \backslash \partial_0\xcF\,, \enl
  \mathring\xcF \etb=
    \xcF \backslash \partial \xcF\,,\enl
  \dot\xcF \etb=
    \xcF \cap \, ]0,\sfk[{}^{\x3}\,.
\eear\ee

For $\,(\la,\mu,\nu) \in \xcF\,$ and $\,\la,\mu\not\in\{0,\sfk\}$,\
equation \eqref{eq:class-cos-range} has a unique solution for
$\,\th\,$ in the range $\,[0,\tfrac\pi2]$,\ and we use this to
define
\be
  a_{\la,\mu}^{\ \nu} = \erm^{\sfi \th \si_2}\,,
  \ \text{with} \quad
  \th\in[0,\tfrac\pi2]\quad\tx{such that}\quad(\sin\th)^2 =
  \frac{ \cos\tfrac{\pi \nu}{\sfk} - \cos\tfrac{\pi(\la+\mu)}{\sfk} }{
  2 \sin\tfrac{\pi\la}{\sfk}\,\sin\tfrac{\pi\mu}{\sfk}}\,.
\ee
By construction, the $\,a_{\la,\mu}^{\ \nu}\,$ satisfy $\,t_\la \cdot
\Ad_{a_{\la,\mu}^{\ \nu}}(t_\mu) = \Ad_b(t_\nu)\,$ for some $\,b \in
\sug$.\ Clearly, only the class of $\,b\,$ in the coset
$\,[b]\in\sug/\sug_{t_\nu}\,$ is relevant, hence we are free to
choose $\,b\,$ in the form $\,b = \erm^{\sfi \varphi' \si_3}\cdot
\erm^{\sfi \th' \si_2}$.\ In order to fix $\,\th'$,\ we take the
trace of $\,\Ad_{a_{\la,\mu}^{\ \nu}}(t_\mu) = t_{-\la} \cdot \Ad_b
(t_\nu)$, whereby $\,a_{\la,\mu}^{\ \nu}\,$ and $\,\erm^{\sfi
\varphi' \si_3}\,$ drop out. The value of $\,\varphi'\,$ can then be
determined by matching the matrix entries of $\,t_\la\cdot
\Ad_{a_{\la,\mu}^{\ \nu}}(t_\mu)\,$ with those of
$\,\Ad_b(t_\nu)$.\ One finds that both $\,\th'\,$ and $\,\varphi'\,$ are defined
uniquely for $\,(\la,\mu,\nu)\in\xcF\,$ with $\la\not\in\{0,
\sfk\}\,$ and $\,\nu\not\in\{\vert\la-\mu\vert,\min(\la+\mu,2\sfk-
\la-\mu)\}$.\ In fact,
$\,\th'\,$ is unique in the larger range, namely, for all
$\,(\la,\mu,\nu)\in\xcF\,$ such that $\la,\nu\not\in\{0,\sfk\}$,\
and so we may use the arbitrariness of $\,\varphi'\,$ on the
difference of the two ranges to extend the formula for
$\,\varphi'\,$ to it, whereby we obtain
\be
b_{\la,\mu}^{\ \nu}=\erm^{\sfi\varphi'\si_3}\cdot\erm^{\sfi\th'
\si_2}\,,\ \text{with}\quad\varphi'=\frac{\pi\la}{2\sfk}\quad
\text{and}\quad\th'\in[0,\tfrac\pi2]\quad\tx{such
that}\quad(\cos\th')^2 =
  \frac{ \cos\tfrac{\pi \mu}{\sfk} - \cos\tfrac{\pi(\la+\nu)}{\sfk} }{
  2 \sin\tfrac{\pi\la}{\sfk}\,\sin\tfrac{\pi\nu}{\sfk}}
\ee
for all $\,(\la,\mu,\nu)\in\xcF\,$ such that $\la,\nu\not\in\{0,
\sfk\}$.\ In their respective ranges, $\,a_{\la,\mu}^{\ \nu}\,$ and
$\,b_{\la,\mu}^{\ \nu}\,$ are smooth functions of $\,(\la,\mu,\nu
)$.\ In particular, both are smooth functions on the subset
$\,\dot\xcF\,$ of $\,\xcF$.

It will be convenient to extend the range of $\,a_{\la,\mu}^{\
\nu}\,$ and $\,b_{\la,\mu}^{\ \nu}\,$ (non-continuously) to include
all triples from $\,\xcF\,$ as
\qq\label{eq:extend-ab}
a_{0,\la}^{\ \la}=a_{\la,0}^{\ \la}=a_{\sfk,\la}^{\ \sfk-\la}=a_{\la
,\sfk}^{\ \sfk-\la}=e \qquad \text{and} \qquad b_{\mu,\mu}^{\
0}=b_{\mu,\sfk-\mu}^{\ \sfk}=e=b_{0,\la}^{\ \la}\,,\qquad
b_{\sfk,\la}^{\ \sfk-\la}=\sfi\,\si_2
\qqq
for arbitrary $\,\la\in[0,\sfk]$ and $\,\mu\in]0,\sfk[$.\ With these
assignments, we have
\be
t_\la\cdot\Ad_{a_{\la,\mu}^{\ \nu}}(t_\mu)=\Ad_{b_{\la,\mu}^{\ \nu}}
(t_\nu)\quad \text{for all} \quad (\la,\mu,\nu) \in \xcF \,.
\labl{eq:tata=btb}
With the help of $\,a_{\la,\mu}^{\ \nu}$,\ we obtain a surjective
map \label{p:F-sug-redund-param}
\be
  \tau\ :\ \xcF \x \sug \to \sug\x\sug\ :\ (\la,\mu,\nu,h)\mapsto
  \big( \Ad_h( t_\la ), \Ad_{h \cdot a_{\la,\mu}^{\ \nu}}(t_\mu) \big) \,.
\ee
The map $\tau$ is {\em not} smooth on all of $\,\xcF \x \sug\,$ but
it is smooth on the subset $\,\dot\xcF \x \sug\,$.\ Furthermore, for
any fixed choice $\,(\la,\mu,\nu) \in \xcF$,\ the map
$\,\tau(\la,\mu,\nu,-)\,$ is a smooth surjection from $\,\sug\,$ to
$\,\Tlmn$,\ and provides a convenient (redundant) parametrisation of
the submanifolds $\,\Tlmn$,\ i.e.\ we have
\be
\Tlmn = \big\{ \ \tau(\la,\mu,\nu,h) \ \big|\ h \in \sug \ \big\}\,.
\labl{eq:Tlmn-as-tau-im}

\begin{Rem}
The map $\,\tau\,$ intertwines the left action $\,\sug \x \bigl(\xcF
\x \sug\bigr) \to \xcF \x \sug$,\ given by $\,g.(\la,\mu,\nu,h) =
(\la,\mu,\nu,g \cdot h)$,\ with the diagonal adjoint action on
$\,\sug\x\sug$.\ It follows that $\,\Tlmn\,$ is a {\em single} orbit
under the diagonal $\Ad$-action of $\,\sug\,$ on $\,\sug\x\sug$.\
This is peculiar to $\,\sug$;\ for other Lie groups, the
corresponding intersection of conjugacy classes as in
\eqref{eq:T_lmn-def} will typically decompose into a continuum of
diagonal $\Ad$-orbits.
\end{Rem}

Below, we shall also need lifts of the maps $\,\txp_1,\txp_2\,$ and
$\,\txm\,$ of \eqref{eq:pi-m-def} to maps $\,\txq_1,
\txq_2,\txq_\txm:\dot\xcF\x\sug\to{}]0,k[{} \x \sug\,$ satisfying
\be
\txp_1\circ\tau=c\circ\txq_1\,,\qquad\qquad\txp_2\circ\tau=c\circ
\txq_2\,,\qquad\qquad\txm\circ\tau=c\circ\txq_m
\labl{eq:q-lift-p}
on $\,\dot\xcF \x \sug$.\ Their definition uses both
$\,a_{\la,\mu}^{\ \nu}\,$ and $\,b_{\la,\mu}^{\ \nu}$,
\be
\txq_1(\la,\mu,\nu,h) = (\la,h)\,,\qquad \txq_2(\la,\mu,\nu,h) =
(\mu,h \cdot a_{\la,\mu}^{\ \nu})\,, \qquad \txq_\txm(\la,\mu,\nu,h)
= (\nu,h \cdot b_{\la,\mu}^{\ \nu})\,.
\ee
Due to the restriction to $\,\dot\xcF$,\ these maps are smooth, and
the property \eqref{eq:q-lift-p} follows from \eqref{eq:tata=btb}.
This allows us to express the pullback-gerbe data on
$\,\sug\x\sug\,$ in terms of the corresponding equivariant data on
$\,\dot\xcF \x \sug$,\ obtained -- in turn -- via pullback of
equivariant data from $\,]0,\sfk[\x\sug$.\ For fixed values of
$\,\la,\mu\,$ and $\,\nu$,\ when smoothness in the weight variables
ceases to play a r\^ole, we shall also consider maps
$\,\txq_{(\la),\mu}^{\ \nu},\txq_{\la,(\mu)}^{\ \nu}\,$ and
$\,\txq_{\la,\mu}^{\ (\nu)}\,$ defined just like $\,\txq_1,\txq_2\,$
and $\,\txq_\txm$,\ respectively, but for an arbitrary fixed triple
$\,(\la,\mu,\nu)\in \xcF$.\ These will be useful when pulling back
objects from fixed conjugacy classes in $\,\sug\,$ to the manifolds
$\,\{(\la,\mu,\nu) \}\x\sug$.

The adjoint action being transitive on $\,\Tlmn$,\ we can write
\qq\label{eq:T-as-coset}
\Tlmn\cong\sug/\cS_{\la,\mu}^{\ \nu} \,,\qquad\qquad
\cS_{\la,\mu}^{\ \nu}=\sug_{t_\la}\cap\Ad_{a_{\la,\mu}^{\
\nu}}\sug_{t_\mu}
\qqq
for $\,\cS_{\la,\mu}^{\ \nu} \subset \sug\,$ the isotropy subgroup
with respect to the diagonal adjoint action. As indicated by
\eqref{eq:Tlmn-as-tau-im}, whenever we describe gerbe data on
$\,\sug\x\sug\,$ in terms of its pullback to (subsets of) $\,\xcF \x
\sug$,\ we have to ensure equivariance with respect to the right
action of $\,\sug\,$ restricted to the isotropy subgroups
$\,\cS_{\la,\mu}^{\ \nu}$.\ One verifies that the latter are given
by
\qq
\cS_{\la,\mu}^{\ \nu} = \left\{\barr{ll}
  \sug \quad & \tx{if} \quad (\la,\mu,\nu) \in \partial_0 \xcF \\
  \uj  \quad & \tx{if} \quad (\la,\mu,\nu) \in \partial_{12} \xcF \\
  \bZ_2 \quad & \tx{if} \quad (\la,\mu,\nu) \in \mathring\xcF
\earr\right.\,.
\qqq
Thus, with the above three possibilities for $\,\cS_{\la,\mu}^{\
\nu}$,\ the manifold $\,\Tlmn\,$ is isomorphic to a point, to a
2-sphere, and to $\,\mathrm{SO(3)}$,\ respectively. The resulting
first cohomology groups are
\be
  H^1\bigl(\Tlmn,\uj\bigr)
  = \left\{\barr{ll}
  \triv \quad & \tx{if} \quad (\la,\mu,\nu) \in \partial \xcF \\
  \bZ_2 \quad & \tx{if} \quad (\la,\mu,\nu) \in \mathring\xcF
  \earr\right.\,.
\labl{eq:cohom-Tlmn}
It now follows from Proposition \ref{prop:2-iso-cl} that in the
interior of $\,\xcF$,\ there appears an obstruction to the existence
of the fusion 2-isomorphism to be constructed in Section
\ref{sec:fusmorph}, which will ultimately lead to the
parity-conservation rule that is the main point of this paper.

\medskip

Finally, note that since
$\,\tau(\la,\mu,\nu,h)=\tau(\la,\mu,\nu,\zeta\cdot h)\,$ for
$\,\zeta=\pm e$,\ the map $\,\tau\,$ factors through
$\,\mathrm{SO}(3)$,
\be
  \xcF \x \sug \xrightarrow{~\tau~} \sug\x\sug
  ~=~
  \xcF \x \sug \xrightarrow{~\pi~} \xcF \x \mathrm{SO}(3) \xrightarrow{~\bar\tau~} \sug\x\sug\,.
\labl{eq:tau-factor}
Here, $\,\pi\,$ is the projection from $\,\sug\,$ to
$\,\mathrm{SO}(3) = \sug/\{\pm e\}\,$ and $\,\bar\tau\,$ is the
induced map on the quotient, which is smooth when restricted to
$\,\dot\xcF\x\sug$.

\subsection*{Pullback gerbes over $\,\boldsymbol{\Tlmn}$}

Our next goal is to pull the canonical gerbe $\,\cG\,$ back to the
manifolds $\,\Tlmn\,$ along each of the three maps $\,\txp_1,\txp_2,
\txm\,$ of \eqref{eq:pi-m-def}. To this end, we should first specify
a surjective submersion over each $\,\Tlmn$,\ which we take in the
form
\qq
Y\Tlmn:=\bigsqcup_{\vec k=(l,m)\in\{0,1\}^2}\,\cT_{\la,\mu;\vec
k}^{\ \nu}\,,\qquad\qquad\cT_{\la,\mu;\vec k}^{\ \nu}= \txp_1^{-1}
(\cC_{\la;l})\cap\txp_2^{-1} (\cC_{\mu;m})\cap \txm^{-1}
(\cC_{\nu;l\dotplus m})\,,
\qqq
where $\,l\dotplus m=\min(l+m,2-l-m) = l + m~(\text{mod}\,2)$.\
Recall that either $\,\cC_{\la;l} = \cC_\la\,$ or $\,\cC_{\la;l} =
\emptyset$,\ etc., and so either $\,\cT_{\la,\mu;\vec k}^{\ \nu} =
\cT_{\la,\mu}^{\ \nu}\,$ or $\,\cT_{\la,\mu;\vec k}^{\ \nu} =
\emptyset$.\ In the above definition, it is understood that we keep
those index pairs $\,\vec k = (l,m)\,$ only for which
$\,\cT_{\la,\mu;\vec k}^{\ \nu}\,$ is non-empty. The advantage of
working with $\,Y\Tlmn\,$ is that it comes with a natural choice of
maps $\,\widehat\txp_i,\ i\in\{1,2 \}\,$ and $\,\widehat\txm\,$ that
cover the respective maps on its base $\,\Tlmn$,\ cf.\ Diagram
\eqref{diag:pullback}, to wit
\qq
\widehat\txp_1\bigl(\vec k,g,h\bigr)=\bigl(l,g\bigr)\,,\qquad\qquad
\alxydim{@C=5em@R=3em}{Y\Tlmn \ar[d]_{\pi_{Y\Tlmn}}
\ar[r]^{\widehat\txp_1} & Y\cC_\la \ar[d]^{\pi_{Y\cC_\la}} \\
\Tlmn \ar[r]_{\txp_1} & \cC_\la}\,,
\qqq
and similarly for $\,\widehat\txp_2\bigl(\vec
k,g,h\bigr)=\bigl(m,h\bigr)$,\ alongside
\qq
\widehat\txm\bigl(\vec k,g,h\bigr)=\bigl(l\dotplus m,g \cdot h \bigr)\,,
\qquad\qquad\alxydim{@C=5em@R=3em}{Y\Tlmn
\ar[d]_{\pi_{Y\Tlmn}} \ar[r]^{\widehat\txm} & Y\cC_\nu
\ar[d]^{\pi_{Y\cC_\nu}} \\ \Tlmn \ar[r]_{\txm} & \cC_\nu}\,.
\qqq
The definition of $\,\widehat\txm\,$ ensures that its image is
non-empty for every triple $\,(\la,\mu,\nu)\in\xcF$.\ Furthermore,
maps $\,\widehat \txp_i,\widehat\txm\,$ canonically induce the
corresponding maps on the fibred square $\,Y^{[2]}\Tlmn\,$ that we
shall also need later. These are, for $\,\vec k_i=(l_i,m_i),\ i\in
\{1,2\}$,\ given by the formul\ae
\qq
&&\widehat\txp_1^{[2]}\bigl(\vec k_1,\vec k_2,g,h\bigr)=
\bigl(l_1,l_2,g \bigr)\,,\qquad\qquad\widehat\txp_2^{[2]}\bigl(\vec
k_1,\vec k_2,g,h\bigr) = \bigl(m_1,m_2 ,h\bigr)\,, \cr\cr
&&\widehat\txm^{[2]}\bigl(\vec k_1,\vec k_2,g,h\bigr)=
\bigl(l_1\dotplus m_1,l_2 \dotplus m_2,g \cdot h\bigr)\,.
\qqq
With these choices, we may pull the canonical gerbe $\,\cG\,$ back
to $\,\Tlmn\,$ along the three maps $\,\txp_1,\txp_2\,$ and
$\,\txm$.\ As usual, it proves more convenient to work in the
equivariant formalism. For that purpose, we introduce the subsets
\be
  \widetilde\cT_{\la,\mu}^{\ \nu}
  = \{(\la,\mu,\nu)\} \x \sug \subset \xcF \x \sug \,.
\ee
Over these sets, we need surjective submersions
$\,\pi_{Y\widetilde\cT_{\la,\mu}^{\
\nu}}:Y\widetilde\cT_{\la,\mu}^{\ \nu} \to
\widetilde\cT_{\la,\mu}^{\ \nu}\,$ together with a map $\,\widehat
\tau\,$ that covers $\,\tau\,$ restricted to
$\,\widetilde\cT_{\la,\mu}^{\ \nu}$.\ In analogy with
\eqref{eq:tilde-O}, we define the sets
\qq
\widetilde\cO_{i,j}^{\ k}=([0,\sfk]\setminus\{\la_{1-i}\})\x([0,\sfk]\setminus\{\la_{1-j}
\})\x([0,\sfk]\setminus\{\la_{1-k}\})\x\sug\,.
\qqq
We now have the natural choices
\qq
Y\widetilde\cT_{\la,\mu}^{\ \nu}&=&\bigsqcup_{\vec k=(l,m)\in\{0,1\}^{\x 2}}\,
\widetilde\cT_{\la,\mu;\vec k}^{\ \nu}\,,\qquad\qquad\widetilde\cT_{\la,\mu;\vec k}^{\ \nu}
=\widetilde\cT_{\la,\mu}^{\ \nu}\cap\widetilde\cO_{l,m}^{\ l\dotplus m}\,,\cr\cr
Y^{[2]}\widetilde\cT_{\la,\mu}^{\ \nu}&=&\bigsqcup_{\vec k_1,\vec k_2\in\{0,1\}^{\x
2}}\,\widetilde\cT_{\la,\mu;\vec k_1,\vec k_2}^{\ \nu}\,,\qquad\qquad
\widetilde\cT_{\la,\mu;\vec k_1,\vec k_2}^{\ \nu}=\widetilde\cT_{\la,\mu}^{\ \nu}\cap\widetilde
\cO_{l_1,m_1}^{\ l_1\dotplus m_1}\cap\widetilde\cO_{l_2,m_2}^{\ l_2
\dotplus m_2}\,,
\qqq
alongside
\qq
\widehat\tau(\vec k,\la,\mu,\nu,h)=\bigl(\vec k,\tau(\la,\mu,\nu,h)
\bigr)\,,\qquad\qquad\widehat\tau^{[2]}(\vec k_1,\vec k_2,\la,\mu,
\nu,h)=\bigl(\vec k_1,\vec k_2,\tau(\la,\mu,\nu,h)\bigr)\,.
\qqq
We want to lift both sides of \eqref{eq:q-lift-p} to the respective
surjective submersions. To this end, we use the map $\,\widehat c\,$
from \eqref{eq:c-hat-def} and introduce the maps
\bea
\widehat\txq_{(\la),\mu}^{\ \nu}(\vec k,\la,\mu,\nu,h)=(l,\la,h)\,,
\qquad\widehat\txq_{\la,(\mu)}^{\ \nu}(\vec k,\la,\mu,\nu,h)=(m,\mu,
h\cdot a_{\la,\mu}^{\ \nu})\,,\enl
\widehat\txq_{\la,\mu}^{\
(\nu )}(\vec k,\la,\mu,\nu,h)=(l\dotplus m,\nu,h\cdot b_{\la,\mu}^{\
\nu} )\,,
\eear\ee
as well as their counterparts for fibred squares of the surjective
submersions,
\qq
\widehat c^{[2]}(i,j,\la,h)=\bigl(i,j,\Ad_h(t_\la)\bigr)
\qqq
and
\bea
\widehat\txq_{(\la),\mu}^{\ \nu\ [2]}(\vec k_1,\vec k_2,\la,\mu,\nu,
h)=(l_1,l_2,\la,h)\,,~~\widehat\txq_{\la,(\mu)}^{\ \nu\ [2]}(
\vec k_1,\vec k_2,\la,\mu,\nu,h)=(m_1,m_2,\mu,h\cdot a_{\la,\mu}^{\
\nu})\,,
\enl
\widehat\txq_{\la,\mu}^{\ (\nu)\ [2]}(\vec k_1,\vec
k_2,\la,\mu,\nu,h)=(l_1 \dotplus m_1,l_2\dotplus m_2,\nu,h\cdot
b_{\la,\mu}^{\ \nu})\,.
\eear\ee

We may next pull back restrictions of the bundle $\,L\,$ along the
maps $\,\widehat\txp_1^{[2]}\circ\widehat\tau^{[2]},\widehat
\txp_2^{[2]}\circ\widehat\tau^{[2]}\,$ and $\,\widehat\txm^{[2]}
\circ\widehat\tau^{[2]}$ to the common base
$\,Y^{[2]}\widetilde\cT_{\la,\mu}^{\ \nu}\,$.\ By construction, we
have that $\,\widetilde
L_{l_1,l_2}\vert_{\widetilde\cC_{\la;l_1,l_2}} \cong \widehat
c^{[2]*} L\vert_{\cC_{\la;l_1,l_2}}\,$ as equivariant line bundles.
Thus, since $\,\widehat\txp_1^{[2]}\circ\widehat\tau^{[2]} =
\widehat c^{[2]} \circ \widehat\txq_{(\la),\mu}^{\ \nu\ [2]}$,\ we
may equivalently pull back the corresponding equivariant bundles
\qq
\bigl( \widehat\txp_1^{[2]}\circ\widehat\tau^{[2]}\bigr)^*
L\vert_{\cC_{\la;l_1,l_2}} \cong \widehat\txq_{(\la),\mu}^{\ \nu\
[2]*}\widetilde L_{l_1,l_2}\vert_{\widetilde\cC_{\la;l_1,l_2}} \,,
\qqq
and similarly for $\,\widetilde L_{m_1,m_2}\,$ and $\,\widetilde
L_{l_1\dotplus m_1,l_2\dotplus m_2}$.\ For the reader's convenience,
and for later reference, we write out the action of the isotropy
group $\,\cS_{\la ,\mu}^{\ \nu}\,$ on the pullback bundles that
enables us to descend them to the quotient $\,\Tlmn$.\ The actions
take the form (with $\,g\in\cS_{\la ,\mu}^{\ \nu}\,$ and
$\,z\in\bC$)
\qq\label{eq:funny-act-1}
(\vec k_1,\vec k_2,\la,\mu,\nu,h,z).g = \bigl(\vec k_1,\vec
k_2,\la,\mu,\nu,h\cdot g,\chi_{\la_{l_1,l_2}}(g)\cdot z \bigr)
\qqq
for $\,\widehat\txq_{(\la),\mu}^{\ \nu\ [2]*}\widetilde L_{l_1,l_2}
\vert_{\widetilde\cC_{\la;l_1,l_2}}$,
\qq
(\vec k_1,\vec k_2,\la,\mu,\nu,h,z).g = \bigl(\vec
k_1,\vec k_2,\la,\mu,\nu,h\cdot g,\chi_{\la_{m_1,m_2}}\bigl(
\Ad_{a_{\la,\mu}^{\ \nu\ -1}}(g)\bigr)\cdot z\bigr)
\qqq
for $\,\widehat\txq_{\la,(\mu)}^{\ \nu\ [2]*}\widetilde L_{m_1,m_2}
\vert_{\widetilde\cC_{\mu;m_1,m_2}}$,\ and
\qq\label{eq:funny-act-3}
(\vec k_1,\vec k_2,\la,\mu,\nu,h,z).g = \bigl(\vec k_1,\vec
k_2,\la,\mu,\nu,h\cdot g,\chi_{\la_{l_1\dotplus m_1,l_2\dotplus
m_2}}\bigl( \Ad_{b_{\la,\mu}^{\ \nu\ -1}}(g)\bigr)\cdot z\bigr)
\qqq
for $\,\widehat\txq_{\la,\mu}^{\ (\nu)\ [2]*}\widetilde
L_{l_1\dotplus m_1,l_2 \dotplus
m_2}\vert_{\widetilde\cC_{\nu;l_1\dotplus m_1,l_2\dotplus m_2}}$.\
Define the map $\,\vep:\sug\to\{-1,1\}\,$ as
\qq\label{eq:epsil-on-sug}
\vep(g)=\left\{\barr{ll}
  -1 \quad & \tx{if} \quad g\in\uj\cdot \sfi\,\si_2 \\
  1 \quad & \tx{otherwise} \earr\right.\,.
\qqq
Using $\,\vep$,\ we may rewrite the action of the isotropy subgroup
on the fibre of $\,\widehat\txq_{\la,(\mu)}^{\ \nu\ [2]*}\widetilde
L_{m_1,m_2} \vert_{\widetilde\cC_{\mu;m_1,m_2}}\,$ and
$\,\widehat\txq_{\la, \mu}^{\ (\nu)\ [2]*}\widetilde L_{l_1\dotplus
m_1,l_2\dotplus m_2} \vert_{\widetilde\cC_{\nu;l_1\dotplus
m_1,l_2\dotplus m_2}}\,$ in the useful form
\qq\label{eq:pull-ab-out-chi}
z\mapsto\chi_{\la_{m_1,m_2}}(g)^{\vep(a_{\la,\mu}^{\ \nu})}\cdot z
\,,\qquad\qquad z\mapsto\chi_{\la_{l_1\dotplus m_1,l_2\dotplus
m_2}}(g)^{\vep(b_{\la,\mu}^{\ \nu})}\cdot z\,,
\qqq
respectively. To see this, we distinguish and check three cases:
\\
\nxt $(\la,\mu,\nu)\in \p_0\xcF$:\ The characters entering the above
expressions are the trivial ones, $\,\chi_0=1$,\ and so there
remains nothing to demonstrate.
\\
\nxt $(\la,\mu,\nu)\in\p_{12} \xcF$:\ Here, $\,\cS_{\la,\mu}^{\
\nu}=\uj$,\ and both $\,a_{\la,\mu}^{\ \nu}\,$ and $\,b_{\la,\mu}^{\
\nu}\,$ take values in the set $\,\{e\}\cup\uj\cdot\sfi\,\si_2$,\
and so \eqref{eq:pull-ab-out-chi} follows from the identities
$\,\Ad_{x^{-1}}(\si_3)=\vep(x)\,\si_3\,$ for
$\,x\in\{e,\sfi\,\si_2\}$.
\\
\nxt $(\la,\mu,\nu)\in \mathring\xcF$:\ In this case,
$\,\cS_{\la,\mu}^{\ \nu}=\bZ_2$,\ and neither $\,a_{\la,\mu}^{\
\nu}\,$ nor $\,b_{\la, \mu}^{\ \nu}\,$ belongs to
$\,\uj\cdot\sfi\,\si_2$.\ The resulting equalities
$\,\vep(a_{\la,\mu}^{\ \nu})=1=\vep(b_{\la, \mu}^{\ \nu})\,$ are
consistent with $\,\Ad_{a_{\la,\mu}^{\ \nu\
-1}}(g)=g=\Ad_{b_{\la,\mu}^{\ \nu\ -1}}( g)\,$ implied by
$\,g\in\bZ_2$.

\subsection*{$\boldsymbol{\bZ_2}$-equivariant line bundles}

In this section, we briefly collect our conventions and some
statements about $\bZ_2$-equivariant line bundles. Of course, the
definitions could equally well be given for other discrete groups,
but we shall only need the $\bZ_2$-case.

Here, and for the remainder of this paper, by a line bundle, we
shall mean a Hermitian line bundle with unitary connection, and by
its connection form, the unique 1-form $\,\widehat A\,$ on the total
space of the bundle with the properties listed, e.g., in
\cite[Def.\,2.2.4]{Brylinski:1993ab}. Furthermore, by a manifold
$\,M\,$ with $\bZ_2$-action $\,\rho$,\ we shall mean a smooth
manifold $\,M\,$ together with diffeomorphisms $\,\rho_\z : M\to
M\,$ for $\,\z \in \bZ_2\,$ such that $\,\rho_\zeta \circ \rho_\xi =
\rho_{\zeta \cdot \xi}\,$ for all $\,\zeta,\xi\in\bZ_2$.

\bef\label{def:z2-equiv-bundle}
Let $\,M\,$ be a manifold with $\bZ_2$-action $\,\rho$.\\[.3em]
(i) A $\bZ_2$-{\em equivariant line bundle over} $\,M\,$ is a pair
$\,(L,\widehat\rho)$,\ where $\,\pi : L \rightarrow M\,$ is a line
bundle, and $\,\widehat\rho: \bZ_2\x L\to L\,$ is an action of
$\,\bZ_2\,$ by diffeomorphisms, i.e.\ $\,\widehat\rho_\zeta \circ
\widehat\rho_\xi = \widehat\rho_{\zeta \cdot \xi}\,$ for all
$\,\zeta,\xi\in\bZ_2$.\ The action preserves fibres, is linear and
unitary on the fibres, and preserves the connection form $\,\widehat
A\in\Om^1(L)\,$ in the sense that $\,{\widehat\rho_\zeta}^* \widehat
A = \widehat A\,$ for all $\,\zeta \in \bZ_2$.\\[.3em]
(ii) Let $\,(L,\widehat\rho)\,$ and $\,(L,'\widehat\rho')\,$ be two
$\bZ_2$-equivariant line bundles over $\,M$.\ A $\bZ_2$-{\em
equivariant isomorphism} $\,L \rightarrow L'\,$ is an isomorphism
$\,\phi : L \rightarrow L'\,$ of Hermitian line bundles with unitary
connection that intertwines the $\bZ_2$-action, $\,\phi \circ
\widehat\rho_\zeta = \widehat\rho'_\zeta \circ \phi\,$ for all
$\,\zeta \in\bZ_2$.
\eef

Altogether, given a manifold with $\bZ_2$-action, we obtain the
category (or groupoid) of $\bZ_2$-equivariant line bundles over
$\,M$,\ with objects and morphisms as described in the above
definition. If the $\bZ_2$-action is free, the quotient space
$\,M/\bZ_2\,$ is again a manifold, and the category of line bundles
over $\,M/\bZ_2\,$ is equivalent to the category of
$\bZ_2$-equivariant line bundles over $\,M$.

Let $\,M\,$ be a manifold with $\bZ_2$-action $\,\rho$,\ and let
$\,A\,$ be a 1-form on $\,M\,$ that is invariant under $\,\rho$,\
i.e.\ $\,\rho_\zeta^* A = A\,$ for all $\,\zeta \in \bZ_2$.\
Consider the trivial bundle $\,M \times \mathbb{C}\,$ with
connection $\,\nabla = d+\frac{1}{\sfi}\,A\,$ and the maps
$\,\widehat\rho_\pm : \bZ_2 \times M \rightarrow M\,$ given by
\be
  \widehat\rho_\pm(m,z) = (\rho_\zeta(m) , (\pm1)^{\vep_\z} \cdot z )\ ,
\ee
written in terms of the function $\,\vep:\bZ_2\to\{0,1\}\,$ with
values
\qq
\vep_\z=\left\{\barr{ll} 0 \quad & \tx{if} \quad \z=e \\
1 \quad & \tx{if} \quad \z=-e \earr\right.\,.
\qqq
One checks that the $\,\widehat\rho_\pm\,$ each turn $\,M \times
\mathbb{C}\,$ into a $\bZ_2$-equivariant line bundle. In order to
keep track of the invariant 1-form $\,A$,\ we denote them as $\,(M
\times \mathbb{C}, A, \widehat\rho_\pm)$.\ We have

\berop\label{prop:Z2-equiv-bundle-2-conn}
Let $\,M\,$ be a 2-connected manifold with $\bZ_2$-action. Every
$\bZ_2$-equivariant line bundle over $\,M\,$ is
$\bZ_2$-equivariantly isomorphic to one of the form $\,(M \times
\mathbb{C}, A, \widehat\rho_\pm)$.
\eerop

\begin{proof}
Since $\,M\,$ is 2-connected, the line bundle $\,L\,$ with
connection $\,\nabla_L\,$ is isomorphic to the trivial bundle $\,M
\x \Cb \to M\,$ with connection $\,\nabla' = \sfd +
\tfrac{1}{\sfi}\,A'\,$ for some 1-form $\,A'\,$ on $\,M$.\ To obtain
a $\bZ_2$-invariant connection, we average the 1-form
with respect to the $\bZ_2$-action, i.e.\ we define $\,\nabla = \sfd
+ \tfrac{1}{\sfi}\,A$,\ where $\,A = \tfrac12\,( A' + \rho_{-e}^*
A')$.\ Since the curvature of $\,L\,$ is $\bZ_2$-invariant, we have
$\,\sfd A = \sfd A'$.\ The set of isomorphism classes of line
bundles with a given curvature is isomorphic to
$\,H^1\bigl(M,\uj\bigr)$,\ which is trivial as $\,M\,$ is
1-connected. Thus, also $\,M \x \Cb \to M\,$ with connection
$\,\nabla\,$ is isomorphic to $\,L$.\ Let $\,f : M \x \Cb \to L\,$
be this isomorphism. Denote by $\,\widehat A_L\,$ the connection
form of $\,\nabla_L$,\ and by $\,\widehat A\,$ that of $\,\nabla$,\
so that $\,\widehat A = f^* \widehat A_L$.

Define a $\bZ_2$-action $\,\sigma\,$ on $\,M \x \Cb \to M\,$ as
$\,\sigma_\zeta = f^{-1} \circ \widehat\rho_\zeta \circ f$.\ Since
$\,\widehat A_L\,$ obeys $\,\widehat\rho_\zeta^* \widehat A_L =
\widehat A_L$,\ we also have $\,\sigma_\zeta^* \widehat A=  \widehat
A$.\ Let us write the $\bZ_2$-action $\,\sigma\,$ as
$\,\sigma_\zeta(m,z) = \bigl(\rho_\zeta(m),s_\zeta(m) \cdot
z\bigr)\,$ for some $\uj$-valued map $\,s_\zeta$.\ For a trivial
bundle, the connection form is given by (cf.\
\cite[Sect.\,2.2]{Brylinski:1993ab})
$\,\widehat A(m,z) = \sfi\, z^{-1}\, \sfd z + A(m)$,\ and so the
condition $\,\sigma_\zeta^* \widehat A= \widehat A\,$ implies
$\,\rho_\zeta^* A = -\sfi\, s_\zeta^{-1}\, \sfd s_\zeta + A$,\
i.e.\ $\,\sfd s_\zeta = 0$.\ Thus, $\,s_\zeta\,$ is locally
constant, and since $\,M\,$ is connected, it is globally constant.
The relation $\,\sigma_\zeta \circ \sigma_\xi = \sigma_{\zeta \cdot
\xi}\,$ implies that $\,s_\zeta\,$ is a $\bZ_2$-character. Finally,
the isomorphism $\,f\,$ is $\bZ_2$-equivariant by construction.
\end{proof}

\subsection*{Stable isomorphisms between pullback gerbes
on $\,\boldsymbol{\Tlmn}$}

The pullback gerbes $\,\txp_1^*\cG,\txp_2^*\cG_2\,$ and
$\,\txm^*\cG\,$ over the product manifold $\,\sug\x\sug\,$ are
closely related, as indicated by the Polyakov--Wiegmann identity
\cite{Polyakov:1984et}
\qq\label{eq:Pol-Wieg}
\txp_1^*\txH+\txp_2^*\txH=\txm^*\txH+\sfd\rho\,,\qquad\qquad\rho=
\tfrac{\sfk}{4\pi}\,\tr\bigl(\txp_1^*\th_L\wedge\txp_2^*\th_R\bigr)
\,,\qquad\th_R=\Ad_\bullet\th_L\,.
\qqq
In fact, we have

\berop
There exists a stable isomorphism
\be
  \cM\ :\ \txp_1^*\cG\star\txp_2^*\cG\xrightarrow{\sim}\txm^*\cG\star\cI(\rho)
\labl{eq:mult-stable-iso}
between pull-back gerbes on $\,\sug\x\sug$,\ and this stable
isomorphism is unique up to a 2-isomorphism.
\eerop

\begin{proof}
By Proposition \ref{prop:stab-iso-cl}, stable-isomorphism classes of
gerbes with a given curvature are in a one-to-one correspondence
with elements of $\,H^2(\sug\x\sug,\uj)$,\ which is trivial. Thus,
\eqref{eq:Pol-Wieg} implies that $\,\cM\,$ exists. By Proposition
\ref{prop:2-iso-cl}, 2-isomorphism classes of stable isomorphisms,
in turn, form a torsor over $\,H^1(\sug\x\sug,\uj)$,\ which is also
trivial. Hence, $\,\Mc\,$ is unique up to a 2-isomorphism.
\end{proof}

The stable isomorphism $\,\cM\,$ constitutes part of the data of a
multiplicative structure on $\,\cG\,$
\cite{Carey:2004xt,Waldorf:2008mult}; here, we shall only require
$\,\cM$.\ An explicit expression for $\,\cM\,$ is currently not
known. In what follows, we shall give such an expression for the
restriction of $\,\cM\,$ to each of the subsets $\,\Tlmn\,$ of
$\,\sug\x\sug$.\ We shall do so in two steps: First, we determine
all 2-isomorphism classes of stable isomorphisms
$\,\txp_1^*\cG\star\txp_2^*\cG\xrightarrow{\sim}\txm^*\cG\star\cI(
\rho)\,$ restricted to $\,\Tlmn\,$ (Lemmas
\ref{lem:multipl-str-restr-dF} and
\ref{lem:multipl-str-restr-intF}), and, then, we identify those
which arise as a restriction of $\,\cM\,$ (Lemmas \ref{lem:chi-chi}
and \ref{lem:fix-int-equiv-str}).

By virtue of Proposition \ref{prop:2-iso-cl}, taken in conjunction
with \eqref{eq:cohom-Tlmn}, there are two 2-isomorphism classes of
stable isomorphisms $\,\txp_1^*\cG\star\txp_2^*\cG\xrightarrow{\sim}
\txm^*\cG\star\cI(\rho)\,$ over $\,\Tlmn\,$ for any triple
$\,(\la,\mu,\nu)\in \mathring\xcF$,\ and a unique class for
$\,(\la,\mu,\nu)\in\p\xcF$.\ Below, we give a representative of each
of these classes in terms of equivariant bundles over the simply
connected covers $\,\widetilde\cT_{\la,\mu}^{\ \nu}\,$ of $\,\Tlmn$.

\medskip

The surjective submersion of the gerbe
$\,\bigl(\txp_1^*\cG\star\txp_2^*\cG\bigr)\vert_{\Tlmn}\,$ is given
by $\,Y\Tlmn\x_{\Tlmn} Y\Tlmn = Y^{[2]}\Tlmn$,\ and that of
$\,\bigl(\txm^*\cG\star\cI(\rho)\bigr)\vert_{\Tlmn}\,$ by
$\,Y\Tlmn\x_{\Tlmn} \Tlmn = Y\Tlmn$.\ The surjective submersion of
the stable isomorphism is therefore $\,Y^{[2]}\Tlmn \x_{\Tlmn}
Y\Tlmn = Y^{[3]}\Tlmn$.\ In the equivariant formulation, the line
bundle on $\,Y^{[3]}\Tlmn\,$ and the isomorphism of pullback line
bundles on $\,Y^{[6]}\Tlmn\,$ are described as follows: Introduce
the canonical projections $\,\pr_i:Y^{[3]}\widetilde\cT_{\la,\mu}^{\
\nu}\to Y\widetilde\cT_{\la,\mu}^{\ \nu},\ i\in\{ 1,2,3\}\,$ and
$\,\pi_{\sug}:Y^{[3]}\widetilde\cT_{\la,\mu}^{\ \nu}\to\sug$.\ The
first datum that we need to give is a $\,\cS_{\la,\mu}^{\
\nu}$-equivariant line bundle $\,\widetilde E_{\la,\mu;\vec k_1,\vec
k_2,\vec k_3}^{\ \nu}\to\widetilde\cT_{\la,\mu;\vec k_1,\vec k_2,
\vec k_3}^{\ \nu}\,$ with connection $\,\nabla_{\widetilde E_{\la,
\mu;\vec k_1,\vec k_2,\vec k_3}^{\ \nu}}\,$ of curvature
\qq
\curv(\nabla_{\widetilde E_{\la,\mu;\vec k_1,\vec k_2,\vec k_3}^{\
\nu}})&=&\pr_3^*\widehat\txq_{\la,\mu}^{\ (\nu)*}\widetilde
B\vert_{\widetilde \cC_{\nu;l_3\dotplus
m_3}}+\pr_3^*\pi_{Y\widetilde\cT_{\la,\mu}^{\
\nu}}^*\tau^*\rho\vert_{\Tlmn}-\pr_1^*\widehat\txq_{(\la),\mu}^{\
\nu*}\widetilde
B\vert_{\widetilde\cC_{\la;l_1}}-\pr_2^*\widehat\txq_{\la,(\mu )}^{\
\nu*}\widetilde B\vert_{\widetilde\cC_{\mu;m_2}}\cr\cr
&=&\pr_1^*\pi_{Y\widetilde\cT_{\la,\mu}^{\ \nu}}^*\widetilde\Om_{\la,\mu}^{\ \nu}+\pi_{\sug}^*
\sfd\widetilde A_{\la,\mu;\vec k_1,\vec k_2,\vec k_3}^{\ \nu}\,.
\label{eq:curv-mult}
\qqq
Here,
\qq
\widetilde\Om_{\la,\mu}^{\ \nu}=\txq_{\la,\mu}^{\ (\nu)*}\widetilde Q\vert_{\widetilde
\cC_\nu}+\tau^*\rho\vert_{\Tlmn}-\txq_{(\la),\mu}^{\ \nu*}\widetilde Q
\vert_{\widetilde\cC_\la}-\txq_{\la,(\mu)}^{\ \nu*}\widetilde Q
\vert_{\widetilde\cC_\mu}\,,
\qqq
and
\qq\label{eq:conn-form-mult}
\widetilde A_{\la,\mu;\vec k_1,\vec k_2,\vec k_3}^{\ \nu}(h)=\sfi\,
\tr\bigl[\bigl((\la-\la_{l_1})\,\La+(\mu-\la_{m_2})\,\Ad_{a_{\la,\mu}^{\
\nu}}(\La)-( \nu-\la_{l_3\dotplus m_3})\,\Ad_{b_{\la,\mu}^{\
\nu}}(\La )\bigr)\,\th_L(h)\bigr]\,.
\qqq
In order to proceed further, we need
\belem\label{lem:Om-ker}
\qq
\widetilde\Om_{\la,\mu}^{\ \nu}\equiv 0\,.
\qqq
\elem

\begin{proof}
The tangent space of $\,\widetilde\cT_{\la,\mu}^{\ \nu}\,$ at a
point $\,(\la,\mu,\nu,h )\,$ is spanned by vectors
$\,\xcV_A(\la,\mu,\nu,h)=R_A(h),\ A\in\{1,2,3\}$,\ where $\,R_A\,$
are the standard right-invariant vector fields on $\,\sug$,\ dual to
the right-invariant Maurer--Cartan 1-forms,
$R_A\con\th_R=\sfi\,\si_A$.\ The vanishing of
$\,\widetilde\Om_{\la,\mu}^{\ \nu}\,$ can be rephrased equivalently
as the identity
\qq\label{eq:Om-annihil}
\xcV_A\con\widetilde\Om_{\la,\mu}^{\ \nu}(\la,\mu,\nu,h)=0\qquad\tx{for
all}\quad A\in\{1,2,3\}\,.
\qqq
Using \eqref{eq:om_la-def} and the shorthand notation $\,(g_1,
g_2)=\tau(\la,\mu,\nu,h)$,\ we readily compute
\qq
\xcV_A\con\txq_{(\la),\mu}^{\ \nu*}\widetilde
Q(\la,\mu,\nu,h)&=&\tfrac{\sfi\,\sfk}{4\pi}\,
\tr\bigl(\si_A\,(\Ad_{g_1}-\Ad_{g_1^{-1}})\th_R(h)\bigr)\,,
\cr\cr
\xcV_A\con\txq_{\la,(\mu)}^{\ \nu*}\widetilde
Q(\la,\mu,\nu,h)&=&\tfrac{\sfi\,\sfk}{4\pi}\,
\tr\bigl(\si_A\,(\Ad_{g_2}-\Ad_{g_2^{-1}})\th_R(h)\bigr)\,,
\cr\cr
\xcV_A\con\txq_{\la,\mu}^{\ (\nu)*}\widetilde
Q(\la,\mu,\nu,h)&=&\tfrac{\sfi\,\sfk}{4
\pi}\,\tr\bigl(\si_A\,(\Ad_{g_1\cdot g_2}-\Ad_{(g_1\cdot g_2)^{-1}})
\th_R(h)\bigr)\,.
\qqq
The 2-form $\,\rho\,$ pulls back to $\,\widetilde\cT_{\la,\mu}^{\
\nu}\,$ as
\qq
\tau^*\rho(\la,\mu,\nu,h)=-\tfrac{\sfk}{4\pi}\,\tr\bigl(\th_R(h)
\wedge(\id_{\sua}-\Ad_{g_1})\circ(\id_{\sua}-\Ad_{g_2})\th_R(h)
\bigr)\,,
\qqq
and so we also have
\qq
\xcV_A\con\tau^*\rho(\la,\mu,\nu,h)=\tfrac{\sfi\,\sfk}{4\pi}\,\tr
\bigl(\si_A\,(\Ad_{g_1}-\Ad_{g_1^{-1}}+\Ad_{g_2}-\Ad_{g_2^{-1}}-
\Ad_{g_1\cdot g_2}+\Ad_{(g_1\cdot g_2)^{-1}})\th_R (h)\bigr)\,.\cr
\qqq
Putting all the above formul\ae ~together, we obtain the desired
result, \eqref{eq:Om-annihil}, which concludes the proof.
\end{proof}

By the above lemma, the curvature has a global primitive, and since
$\,\widetilde\cT_{\la,\mu;\vec k_1,\vec k_2,\vec k_3}^{\ \nu}\,$ are
simply connected, the line bundle under consideration is trivial,
\qq
\widetilde E_{\la,\mu;\vec k_1,\vec k_2,\vec k_3}^{\ \nu}=\widetilde
\cT_{\la,\mu;\vec k_1,\vec k_2,\vec k_3}^{\ \nu}\x\bC\,.
\qqq
The connection is given by
\be
\nabla_{\widetilde E_{\la,\mu;\vec k_1,\vec k_2,\vec k_3}^{\ \nu}}=
\sfd+\tfrac{1}{\sfi}\,\pi_{\sug}^*\widetilde A_{\la,\mu;\vec k_1,
\vec k_2,\vec k_3}^{\ \nu}(h)\,.
\labl{eq:L-lmn-connection}

The second piece of data of a stable isomorphism between the pullback
gerbes is -- in the equivariant formulation -- a
family of $\cS_{\la,\mu}^{\ \nu}$-equivariant isomorphisms
\bea
\widetilde\a_{\la,\mu;\vec k_1,\vec k_2,\vec k_3,\vec k_4,\vec k_5,
\vec k_6}^{\ \nu}\  : \ \pr_{1,4}^*\widehat\txq_{(\la),\mu}^{\ \nu\
[2]*}\widetilde L_{l_1,l_4}\ox\pr_{2,5}^*\widehat\txq_{\la,(\mu
)}^{\ \nu\ [2]*}\widetilde L_{m_2,m_5}\ox\pr_{4,5,6}^* \widetilde
E_{\la,\mu;\vec k_4,\vec k_5,\vec k_6}^{\ \nu}\enl\hspace{10em}
\xrightarrow{~\sim~} ~ \pr_{1,2,3}^*\widetilde E_{\la,\mu;\vec k_1,\vec
k_2,\vec k_3}^{\ \nu}\ox\pr_{3,6}^*\widehat\txq_{\la,\mu}^{\ (\nu)\
[2]*}\widetilde L_{l_3\dotplus m_3,l_6\dotplus m_6}
\label{eq:mult-rest-stab-iso-iso}
\eear\ee
of line bundles over $\,\widetilde\cT_{\la,\mu;\vec k_1,\vec
k_2,\vec k_3,\vec k_4,\vec k_5,\vec k_6}^{\ \nu}$,\ the latter space
being equipped with the canonical projections $\,\pr_{i,j}(\vec
k_1,\dots,\vec k_6,\la,\mu,\nu,h)=(\vec k_i,\vec
k_j,\la,\mu,\nu,h)\,$ and $\,\pr_{i,j,k}(\vec k_1,\dots,\vec
k_6,\la,\mu,\nu,h)=(\vec k_i,\vec k_j,\vec k_k,$ $\la,\mu,\nu,h)$.\
As demonstrated below, the $\,\widetilde\a\,$ are of the form
\bea
\widetilde\a_{\la,\mu;\vec k_1,\vec k_2,\vec k_3,\vec k_4,\vec k_5,
\vec k_6}^{\ \nu}(\vec k_1,\vec k_2,\vec k_3,\vec k_4,\vec k_5,\vec
k_6,\la,\mu,\nu,h,z\ox z'\ox z'')
\enl \qquad
=(\vec k_1,\vec k_2,\vec
k_3,\vec k_4,\vec k_5,\vec k_6,\la,\mu,\nu,h,1\ox z\cdot z'\cdot z''
)\,.
\eear\labl{eq:alpha-for-Tbun}
It is manifest that the maps $\,\widetilde\a\,$ are unitary, associative and
compatible with the groupoid structures of the gerbes involved. It
is also easy to see that they preserve the bundle connections.

It remains to give the $\cS_{\la,\mu}^{\ \nu}$-action on
$\,\widetilde E_{\la,\mu;\vec k_1,\vec k_2,\vec k_3}^{\ \nu}\,$ and
verify equivariance of $\,\nabla_{\widetilde E_{\la,\mu;\vec
k_1,\vec k_2,\vec k_3}^{\ \nu}}\,$ and $\,\widetilde\a_{\la,\mu;\vec
k_1,\vec k_2,\vec k_3, \vec k_4,\vec k_5,\vec k_6}^{\ \nu}$.\ Once
this is done, we can descend the data to $\,\Tlmn$.
We shall denote the resulting stable isomorphism as
\be
\Phi_{\la,\mu}^{\ \nu}(\widetilde E_{\la,\mu;\vec k_1,\vec k_2,\vec
k_3}^{\ \nu})\ :\ \bigl(\txp_1^*\cG\star\txp_2^*\cG\bigr)|_{\Tlmn}
\xrightarrow{\sim}\bigl(\txm^*\cG\star\cI(\rho)\bigr)|_{\Tlmn}\,.
\labl{eq:Tlmn-1-iso}

For $\,(\la,\mu,\nu)\in\p\xcF$,\ one verifies that
\be
\Ad_{a_{\la,\mu}^{\ \nu}}(\si_3)=\vep(a_{\la,\mu}^{\ \nu})\,\si_3
\qquad \text{and} \qquad \Ad_{b_{\la,\mu}^{\
\nu}}(\si_3)=\vep(b_{\la,\mu}^{\ \nu})\, \si_3\,,
\labl{eq:Ad-ab-sig3}
with $\,\eps(\,\cdot\,)\,$ as defined in \eqref{eq:epsil-on-sug}.
Comparing the exponents in \eqref{eq:tata=btb},
we thus obtain the relation
\be
\la+\vep(a_{\la,\mu}^{\ \nu})\,\mu-\vep(b_{\la ,\mu}^{\ \nu})\,\nu=2
\sfk\,n_{\la,\mu}^{\ \nu}\,,\qquad\qquad n_{\la,\mu}^{\ \nu}\in\bZ\ .
\labl{eq:n-lmn-def}
The definition of $\,n_{\la,\mu}^{\ \nu}\,$ enters the formulation
of the following lemma.

\belem\label{lem:multipl-str-restr-dF}
For $\,(\la,\mu,\nu)\in\p\xcF$,\ any 1-isomorphism $\,\bigl(
\txp_1^*\cG\star\txp_2^*\cG \bigr)|_{\Tlmn} \xrightarrow{\sim}
\bigl( \txm^*\cG\star\cI(\rho) \bigr)|_{\Tlmn}\,$ is 2-isomorphic to
$\,\Phi_{\la,\mu}^{\ \nu}(\widetilde E_{\la,\mu;\vec k_1,\vec k_2,
\vec k_3}^{\ \nu})\,$ with
\qq
\widetilde E_{\la,\mu;\vec k_1,\vec k_2,\vec k_3}^{\ \nu}=
\pi_{\sug}^*\overline K_{2\sfk\,n_{\la,\mu}^{\ \nu}-\la_{l_1}-\vep(
a_{\la,\mu}^{\ \nu})\,\la_{m_2}+\vep(b_{\la,\mu}^{\ \nu})\,\la_{l_3
\dotplus m_3}}\,,
\qqq
and with the $\cS_{\la,\mu}^{\ \nu}$-equivariant structure inherited from
$\,\overline K_{2\sfk \,n_{\la,\mu}^{\ \nu}-\la_{l_1}-\vep(a_{\la,
\mu}^{\ \nu})\,\la_{m_2}+\vep(b_{\la,\mu}^{\ \nu})\,\la_{l_3\dotplus
m_3}}$,\ as defined in \eqref{eq:KKS-S-action}. \elem

\begin{proof}
Let us denote $\,\xi_{\vec k_1,\vec k_2,\vec k_3} =
2\sfk\,n_{\la,\mu}^{\ \nu}-\la_{l_1}-\vep(a_{\la, \mu}^{\
\nu})\,\la_{m_2}+\vep(b_{\la,\mu}^{\ \nu})\,\la_{l_3\dotplus m_3}$.\
When $\,(\la,\mu,\nu)\in\p\xcF$,\ the expression
\eqref{eq:conn-form-mult} for the 1-form defining the connection
$\,\nabla_{\widetilde E_{\la,\mu;\vec k_1,\vec k_2,\vec k_3}^{\
\nu}}=\sfd+\frac{1}{\sfi}\,\pi_{\sug}^*\widetilde A_{\la,\mu;\vec
k_1,\vec k_2,\vec k_3}^{\ \nu}(h)$\, on $\,\widetilde E_{\la,\mu;
\vec k_1,\vec k_2,\vec k_3}^{\ \nu}\,$ simplifies due to the
relations \eqref{eq:Ad-ab-sig3} and \eqref{eq:n-lmn-def}, giving
\be
\widetilde A_{\la,\mu;\vec k_1,\vec k_2,\vec k_3}^{\ \nu} = \sfi\,
\xi_{\vec k_1,\vec k_2,\vec k_3} \,\tr \bigl(\La\,\th_L(h)\bigr)\,.
\ee
As $\,\widetilde E_{\la,\mu;\vec k_1,\vec k_2,\vec k_3}^{\ \nu}\,$
is a trivial bundle, we infer, by comparing the respective
connection 1-forms, that it can be identified with the pullback of
the equivariant bundle $\,\overline K_{\xi_{\vec k_1,\vec k_2,\vec
k_3}}\to\sug\,$ along the (bijective) projection $\,\pi_{\sug} :
\widetilde\cT_{\la,\mu}^{\ \nu} \rightarrow \sug$.\ The
$\cS_{\la,\mu}^{\ \nu}$-equivariant structure on $\,\widetilde
E_{\la,\mu;\vec k_1,\vec k_2,\vec k_3}^{\ \nu}\,$ is obtained by
restricting that on $\,\overline K_{\xi_{\vec k_1,\vec k_2,\vec
k_3}}\,$ to $\,\cS_{\la,\mu}^{\ \nu}$.\ This is possible because the
isotropy subgroup $\,\sug_{\xi_{\vec k_1,\vec k_2,\vec k_3}}\,$ of
$\,\overline K_{\xi_{\vec k_1,\vec k_2,\vec k_3}}\,$ contains
$\,\cS_{\la,\mu}^{\ \nu}\,$ as it always contains $\,\uj$,\ and
equals $\,\sug\,$ for $\,\la,\mu\in\{0,\sfk \}$.\ According to
\eqref{eq:funny-act-1}--\eqref{eq:pull-ab-out-chi}, the
$\cS_{\la,\mu}^{\ \nu}$-equivariance of
$\,\widetilde\a_{\la,\mu;\vec k_1,\vec k_2,\vec k_3,\vec k_4,\vec
k_5,\vec k_6}^{\ \nu}\,$ amounts to
\bea
\widetilde\a_{\la,\mu;\vec k_1,\vec k_2,\vec k_3,\vec k_4,\vec k_5,
\vec k_6}^{\ \nu}\bigl(h\cdot s,\chi_{\la_{l_1,l_4}}(s)\cdot z\ox
\chi_{\la_{m_2,m_5}}(s)^{\vep( a_{\la,\mu}^{\ \nu})}\cdot z'\ox
\chi_{\xi_{\vec k_4,\vec k_5,\vec k_6}}(s)\cdot z''\bigr) \enl
\qquad \qquad = \bigl(h\cdot s,\chi_{\xi_{\vec k_1,\vec k_2,\vec
k_3}}(s)\otimes\chi_{\la_{l_3\dotplus m_3,l_6\dotplus m_6}}(s)^{\vep
(b_{\la,\mu}^{\ \nu})}z\cdot z'\cdot z'')
\eear\ee
for $\,s\in\cS_{\la,\mu}^{\ \nu}$.\ The above identity follows
immediately from
\qq
\la_{l_1,l_4}+\vep(a_{\la,\mu}^{\ \nu})\,\la_{m_2,m_5}+2\sfk\,n_{\la
,\mu}^{\ \nu}-\la_{l_4}-\vep(a_{\la,\mu}^{\ \nu})\,\la_{m_5}+\vep(
b_{\la,\mu}^{\ \nu})\,\la_{l_6\dotplus m_6}\cr\cr
=2\sfk\,n_{\la,\mu}^{\ \nu}-\la_{l_1}-\vep(a_{\la,\mu}^{\ \nu})\,
\la_{m_2}+\vep(b_{\la,\mu}^{\ \nu})\,\la_{l_3\dotplus m_3}+\vep(
b_{\la,\mu}^{\ \nu})\,\la_{l_3\dotplus m_3,l_6\dotplus m_6}\,.
\qqq
\end{proof}

For $\,(\la,\mu,\nu)\,$ from the interior of $\,\xcF$,\ there are
two 2-isomorphism classes of stable isomorphisms, as described by
the following lemma.

\belem\label{lem:multipl-str-restr-intF}
For $\,(\la,\mu,\nu)\in\mathring\xcF$,\ any 1-isomorphism $\,\bigl(
\txp_1^*\cG\star\txp_2^*\cG\bigr)|_{\Tlmn} \xrightarrow{\sim}
\bigl( \txm^*\cG\star\cI(\rho) \bigr)|_{\Tlmn}\,$ is 2-isomorphic to
$\,\Phi_{\la,\mu}^{\ \nu}(\widetilde E_{\la,\mu;\vec k_1,\vec k_2,
\vec k_3}^{\ \nu})\,$ with
\qq
\widetilde E_{\la,\mu;\vec k_1,\vec k_2,\vec k_3}^{\ \nu}=
\pi_{\sug}^*(\sug\x\bC)
\qqq
the trivial line bundle with
the connection 1-form given by \eqref{eq:conn-form-mult}, and a
$\bZ_2$-equivariant structure inherited from that on
$\,\sug\x\bC\to\sug\,$ which lifts the action of $\,\cS_{\la,\mu}^{\
\nu}\cong \bZ_2\,$ on the base $\,\sug\,$ to the total space as
\qq
(\sug\x\bC)\x\bZ_2\to\sug\x\bC\ &:&\ (h,z,\z)\mapsto\bigl(h\cdot\z,
\chi_{\la,\mu;\vec k_1,\vec k_2,\vec k_3}^{\ \nu}(\z)\cdot z\bigr)
\,,\cr\cr \chi_{\la,\mu;\vec k_1,\vec k_2,\vec k_3}^{\ \nu}\ :\
\bZ_2\to\uj\ &:&\ \z \mapsto (-1)^{\vep_\z\,(\la_{l_3\dotplus
m_3}-\la_{l_1}-\la_{m_2}+\vep_{\la,\mu}^{\ \nu})}\,,
\qqq
and where either $\,\vep_{\la,\mu}^{\ \nu} = 0\,$ or $\,\vep_{\la,\mu}^{\ \nu} = 1$.
\elem

\begin{proof}
For notational brevity, we shall, in this proof, omit the
isomorphism
$\,\pi_{\sug} : \widetilde\cT_{\la,\mu;\vec k_1,\vec k_2,\vec
k_3}^{\ \nu} = \{(\vec k_1,\vec k_2,\vec k_3,\la,\mu,\nu)\} \x \sug
\overset{\sim}{\rightarrow} \sug\,$ and identify
$\,\widetilde\cT_{\la,\mu;\vec k_1,\vec k_2,\vec k_3}^{\ \nu}\,$
with $\,\sug$.\ Each of the spaces
$\,\widetilde\cT_{\la,\mu;\vec k_1,\vec k_2,\vec k_3}^{\ \nu}\,$
carries a $\bZ_2$-action given by multiplication by elements of the
centre of $\,\sug$.\ The 1-forms
$\,\widetilde A_{\la,\mu;\vec k_1,\vec k_2,\vec k_3}^{\ \nu}\,$ are
$\bZ_2$-invariant.

We need to give a $\bZ_2$-equivariant line bundle on each of the
spaces $\,\widetilde\cT_{\la,\mu;\vec k_1,\vec k_2,\vec k_3}^{\
\nu}\,$ with curvature $\,\sfd \widetilde  A_{\la,\mu;\vec k_1,\vec
k_2,\vec k_3}^{\ \nu}$.\ Since $\,\widetilde\cT_{\la,\mu;\vec
k_1,\vec k_2,\vec k_3}^{\ \nu}\,$ is 2-connected, each such bundle
is, by Proposition \ref{prop:Z2-equiv-bundle-2-conn},
$\bZ_2$-equivariantly isomorphic to a trivial bundle with
$\bZ_2$-action as given in the proposition. Specifically, up to a
$\bZ_2$-equivariant isomorphism, the most general
$\bZ_2$-equivariant line bundle is
\qq
\widetilde E_{\la,\mu;\vec k_1,\vec k_2,\vec k_3}^{\ \nu}=\sug\x\bC
\,,
\qqq
and the $\bZ_2$-action, given by
\qq\label{eq:gen-Z2-equiv-bdle}
\bigl( \zeta , (h,z) \bigr) \mapsto \bigl(h\cdot\z,\chi_{\la
,\mu;\vec k_1,\vec k_2,\vec k_3}^{\ \nu}(\z)\cdot z\bigr)\,,
\qqq
is determined by a family of $\bZ_2$-characters
$\,\chi_{\la,\mu;\vec k_1,\vec k_2,\vec k_3}^{\ \nu}$.\ These
characters are constrained by the requirement that
$\,\widetilde\a_{\la,\mu;\vec k_1,\vec k_2,\vec k_3, \vec k_4,\vec
k_5,\vec k_6}^{\ \nu}\,$ be a $\bZ_2$-equivariant bundle
isomorphism. As a bundle isomorphism, it is necessarily of the form
\qq\label{eq:al-gen}
\widetilde\a_{\la,\mu;\vec k_1,\vec k_2,\vec k_3,\vec k_4,\vec k_5,
\vec k_6}^{\ \nu}(h,z\ox z'\ox z'')=(h,1\ox\psi_{\la,\mu;\vec k_1,
\vec k_2,\vec k_3,\vec k_4,\vec k_5,\vec k_6}^{\ \nu}\cdot z\cdot
z'\cdot z'')
\qqq
for some constant phases $\,\psi_{\la,\mu;\vec k_1,\vec k_2,\vec k_3
,\vec k_4,\vec k_5,\vec k_6}^{\ \nu}\in\uj$.\ This follows from the
comparison of the respective connection 1-forms on both sides of
\eqref{eq:mult-rest-stab-iso-iso} (implying local constancy of the
phases), in conjunction with connectedness of their common base
(implying their global constancy). Equivariance of
$\,\widetilde\a_{\la,\mu;\vec k_1,\vec k_2,\vec k_3,\vec k_4,\vec
k_5,\vec k_6}^{\ \nu}\,$ with respect to the action of the isotropy
subgroup $\,\cS_{\la,\mu}^{\ \nu}\,$ on the bundles involved amounts
to the requirement that
\qq
\widetilde\a_{\la,\mu;\vec k_1,\vec k_2,\vec k_3,\vec k_4,\vec k_5,
\vec k_6}^{\ \nu} \bigl(h\cdot\z,\chi_{\la_{l_1,l_4}}(\z)\cdot z\ox
\chi_{\la_{m_2,m_5}}\bigl(\Ad_{a_{\la,\mu}^{\ \nu\
-1}}(\z)\bigr)\cdot z'\ox\chi_{\la,\mu;\vec k_4,\vec k_5,\vec
k_6}^{\ \nu}(\z)\cdot z'' \bigr)\cr\cr
=\bigl(h\cdot\z,\chi_{\la,\mu;\vec k_1,\vec k_2,\vec k_3}^{\
\nu}(\z)\ox\chi_{\la_{l_3\dotplus m_3,l_6\dotplus m_6}}
\bigl(\Ad_{b_{\la,\mu}^{\ \nu\ -1}}(\z)\bigr)\cdot\psi_{\la,\mu;\vec
k_1 ,\vec k_2,\vec k_3,\vec k_4,\vec k_5,\vec k_6}^{\ \nu}\cdot
z\cdot z'\cdot z''\bigr)
\qqq
hold true for arbitrary $\,\z\in\bZ_2$.\ This is equivalent to the
equality
\qq
\chi_{\la,\mu;\vec k_1,\vec k_2,\vec k_3}^{\ \nu}\cdot
\chi_{\la_{l_1}+\la_{m_2}-\la_{l_3\dotplus m_3}}=\chi_{\la,\mu;\vec
k_4,\vec k_5,\vec k_6}^{\ \nu}\cdot\chi_{\la_{l_4}+\la_{m_5}-
\la_{l_6\dotplus m_6}}
\qqq
of functions on $\,\bZ_2$.\ It follows that both sides have to be
independent of
$\,\vec k_1,\vec k_2,\vec k_3$,\ so that
\qq
\chi_{\la,\mu;\vec k_1,\vec k_2,\vec k_3}^{\ \nu}(\z)=
\chi_{\la_{l_3\dotplus m_3}-\la_{l_1}-\la_{m_2}}(\z) \cdot
(-1)^{\vep_\z\,\vep_{\la,\mu}^{\ \nu}}
\qqq
for some constants $\,\vep_{\la,\mu}^{\ \nu} \in \{ 0,1 \}$.

We may now proceed to constrain $\,\psi_{\la,\mu; \vec k_1,\vec k_2,
\vec k_3,\vec k_4,\vec k_5,\vec k_6}^{\ \nu}\,$ by demanding
compatibility of the isomorphism $\,\widetilde\a_{\la,\mu;\vec
k_1,\vec k_2, \vec k_3,\vec k_4,\vec k_5,\vec k_6}^{\ \nu}\,$ with
the (trivial) groupoid structures on
$\,\txp_1^*\cG\star\txp_2^*\cG\,$ and $\,\txm^*\cG\star\cI(\rho )$.\
When rewritten in terms of the data $\,\psi_{\la,\mu;\vec k_1, \vec
k_2,\vec k_3,\vec k_4,\vec k_5,\vec k_6}^{\ \nu}$,\ Diagram
\eqref{diag:stab-iso-comp-mu} yields the relation
\qq
\psi_{\la,\mu;\vec k_1,\vec k_2,\vec k_3,\vec k_7,\vec k_8,\vec
k_9}^{\ \nu}=\psi_{\la,\mu;\vec k_1,\vec k_2,\vec k_3,\vec k_4,\vec
k_5,\vec k_6}^{\ \nu}\cdot \psi_{\la,\mu;\vec k_4,\vec k_5,\vec
k_6,\vec k_7,\vec k_8,\vec k_9}^{\ \nu}\,.
\qqq
Setting all $\,\vec k_i\,$ equal to $\,\vec 0\,$ shows that
$\,\psi_{\la,\mu;\vec 0,\vec 0,\vec 0,\vec 0,\vec 0,\vec 0}^{\ \nu}=
1$;\ setting $\,\vec k_1 = \vec k_2 = \vec k_3 = \vec k_7 = \vec k_8
= \vec k_9 = \vec 0\,$ gives $\,\psi_{\la,\mu;\vec 0,\vec 0,\vec
0,\vec k_4,\vec k_5,\vec k_6}^{\ \nu} = \psi_{\la,\mu;\vec k_4,\vec
k_5,\vec k_6,\vec 0,\vec 0,\vec 0}^{\ \nu}{}^{-1}$;\ and, finally,
setting $\,\vec k_4 = \vec k_5 = \vec k_6 = \vec 0\,$ implies,
in conjunction with the above results, that
\be
\psi_{\la,\mu;\vec k_1,\vec k_2,\vec k_3,\vec k_4,\vec k_5,\vec
k_6}^{\ \nu}= \frac{ \psi_{\la,\mu;\vec k_1,\vec k_2,\vec k_3,\vec
0,\vec 0,\vec 0}^{\ \nu} }{\psi_{\la,\mu;\vec k_4,\vec k_5,\vec k_6,
\vec 0,\vec 0,\vec 0}^{\ \nu} }\,.
\ee

The line bundles $\,\widetilde E_{\la,\mu;\vec k_1,\vec k_2,\vec
k_3}^{\ \nu}\,$ together with the isomorphisms
$\,\widetilde\a_{\la,\mu;\vec k_1,\vec k_2,\vec k_3,\vec k_4,\vec
k_5,\vec k_6}^{\ \nu}\,$ provide the $\bZ_2$-equivariant formulation
of the data of a stable isomorphism $\,\bigl(
\txp_1^*\cG\star\txp_2^*\cG \bigr)|_{\Tlmn} \xrightarrow{\sim}
\bigl( \txm^*\cG\star\cI(\rho) \bigr)|_{\Tlmn}$.\ The latter is
2-isomorphic to the stable isomorphism obtained by replacing
\be
\psi_{\la,\mu;\vec k_1,\vec k_2,\vec k_3,\vec k_4,\vec k_5,\vec
k_6}^{\ \nu} \leadsto \psi_{\la,\mu;\vec k_1,\vec k_2,\vec k_3,\vec
k_4,\vec k_5,\vec k_6}^{\ \nu}  \cdot \frac{\gamma_{\la,\mu;\vec
k_1,\vec k_2,\vec k_3}^{\ \nu}}{\gamma_{\la,\mu;\vec k_4,\vec
k_5,\vec k_6}^{\ \nu}}
\ee
for some $\uj$-valued constants $\,\gamma_{\la,\mu;\vec k_1,\vec
k_2,\vec k_3}^{\ \nu}$.\ In order to see this, just take the trivial
submersion in Definition \ref{def:2-iso} and the locally constant
maps $\,(\vec k_1,\vec k_2,\vec k_3,h,z)\mapsto(\vec k_1,\vec
k_2,\vec k_3,h,\gamma_{\la,\mu;\vec k_1,\vec k_2,\vec k_3}^{\ \nu}
\cdot z)\,$ as bundle isomorphisms. This allows to choose
$\,\psi_{\la,\mu;\vec k_1,\vec k_2,\vec k_3,\vec 0,\vec 0,\vec 0}^{\
\nu} = 1$,\ and hence also $\,\psi_{\la,\mu;\vec k_1,\vec k_2,\vec
k_3,\vec k_4,\vec k_5,\vec k_6}^{\ \nu}=1$.\ Thus,
$\,\widetilde\a_{\la,\mu;\vec k_1,\vec k_2,\vec k_3,\vec k_4,\vec
k_5,\vec k_6}^{\ \nu}\,$ are of the form claimed in
\eqref{eq:alpha-for-Tbun}.
\end{proof}

\subsection*{The stable isomorphism of the multiplicative structure
on $\,\boldsymbol{\Tlmn}$}

By Proposition \ref{eq:mult-stable-iso}, there is a unique stable
isomorphism $\,\cM\,$ between the pullback gerbes $\,\txp_1^*\cG
\star\txp_2^*\cG\,$ and $\,\txm^*\cG\star\cI( \rho)\,$ on
$\,\sug\x\sug$.\ As proved in Lemma \ref{lem:multipl-str-restr-dF},
the restriction of $\,\cM\,$ to $\,\Tlmn\,$ with
$\,(\la,\mu,\nu)\in\partial\xcF\,$ has to be 2-isomorphic to the
stable isomorphism described there. On the other hand, according to
Lemma \ref{lem:multipl-str-restr-intF}, there are two choices of the
stable isomorphism for $\,(\la,\mu,\nu)\in\mathring\xcF$.\ In this
section, we describe which one of these is 2-isomorphic to the
restriction of $\,\cM\,$ to $\,\Tlmn$.

Let us commence our study by fixing indices $\,(\vec k_1,\vec
k_2,\vec k_3)\,$ and defining the 2-connected space
\qq
\widetilde\cT_{\vec k_1,\vec k_2,\vec k_3}=\{(\vec k_1,\vec k_2,\vec
k_3)\}\x\dot\xcF\x\sug\,.
\qqq
Note, in particular, that the factor $\,\dot\xcF\,$ contains points
from $\,\p_{12}\xcF$.\ This will be of prime significance to our
subsequent considerations. As discussed on
p.\,\pageref{p:F-sug-redund-param}, the space $\,\widetilde\cT_{\vec
k_1,\vec k_2,\vec k_3}\,$ parametrises a connected submanifold
$\,\tau\bigl(\dot\xcF\x\sug\bigr)\subset\sug\x\sug$,\ and the
mapping $\,\tau\,$ factorises through the $\bZ_2$-orbifold
$\,\widetilde\cT_{\vec k_1,\vec k_2,\vec k_3}\x\sug/\bZ_2$.\ Hence,
$\,\cM\,$ defines a $\bZ_2$-equivariant line bundle
\qq
\widetilde\cM_{\vec k_1,\vec k_2,\vec k_3}\to\widetilde\cT_{\vec
k_1,\vec k_2,\vec k_3}\,.
\qqq
As $\,\widetilde\cT_{\vec k_1,\vec k_2,\vec k_3}\,$ is 2-connected,
$\,\widetilde\cM_{\vec k_1,\vec k_2,\vec k_3}\,$ is, by Proposition
\ref{prop:Z2-equiv-bundle-2-conn}, $\bZ_2$-equivariantly isomorphic
to the trivial bundle $\,\widetilde\cT_{\vec k_1,\vec k_2 ,\vec k_3}
\x \Cb\,$ with a connection determined by a $\bZ_2$-invariant 1-form
$\,\widetilde P_{\vec k_1,\vec k_2,\vec k_3}\,$ on
$\,\widetilde\cT_{\vec k_1,\vec k_2,\vec k_3}$.\ To avoid
introducing yet another symbol, we denote this trivial (and
explicitly trivialised) bundle by $\,\widetilde\cM_{\vec k_1,\vec
k_2,\vec k_3}\,$ from now on. The $\bZ_2$-equivariance is
implemented on $\,\widetilde\cM_{\vec k_1,\vec k_2,\vec k_3}\,$ as
\qq
&&\bZ_2\x \widetilde\cM_{\vec k_1,\vec k_2,\vec k_3}\to\widetilde
\cM_{\vec k_1,\vec k_2,\vec k_3}\cr\cr \ &:&\ \bigl(\z,(\vec k_1,
\vec k_2,\vec k_3,\la,\mu,\nu,h,z)\bigr)\mapsto\bigl(\vec k_1,\vec
k_2,\vec k_3,\la,\mu,\nu,\z\cdot h,\eta_{\vec k_1,\vec k_2,\vec k_3}
(\z)\cdot z\bigr)
\qqq
by characters $\,\eta_{\vec k_1,\vec k_2,\vec k_3}:\bZ_2\to\uj$.\ In
particular, the $\,\eta_{\vec k_1,\vec k_2,\vec k_3}\,$
are independent of $\,\la,\mu,\nu\,$ and $\,h$.\

As already mentioned above, the stable isomorphism $\,\cM\,$
restricted to $\,\Tlmn\,$ is 2-isomorphic to one of the stable
isomorphisms constructed in Lemmas \ref{lem:multipl-str-restr-dF}
and \ref{lem:multipl-str-restr-intF}. By Proposition
\ref{prop:2-iso-triv-subm}, the surjective submersion of the
2-isomorphism can be chosen trivial, and so $\,\cM\,$ yields an
isomorphism of line bundles over $\,\cT_{\la,\mu;\vec k_1,\vec
k_2,\vec k_3}^{\ \nu}$.\ In the $\bZ_2$-equivariant formulation,
this provides an isomorphism of $\bZ_2$-equivariant line bundles.
Namely, let $\,\unl{\chi}_{\la,\mu;\vec k_1,\vec k_2,\vec k_3}^{\
\nu}\,$ be the data of the $\cS_{\la,\mu}^{\ \nu}$-equivariant
structure on that bundle $\,\widetilde E_{\la,\mu;\vec k_1,\vec
k_2,\vec k_3}^{\ \nu}\,$ which is $\bZ_2$-equivariantly isomorphic
to $\,\widetilde\cM_{\vec k_1,\vec k_2,\vec
k_3}\vert_{\widetilde\cT_{\la,\mu;\vec k_1,\vec k_2,\vec k_3}^{\
\nu}}$.\ We now have the following

\belem\label{lem:chi-chi}
$\eta_{\vec k_1,\vec k_2,\vec k_3}=\unl{\chi}_{\la,\mu;\vec k_1,\vec
k_2,\vec k_3}^{\ \nu}\big\vert_{\bZ_2}$. \elem

\begin{proof}
Let us identify $\,\cT_{\la,\mu;\vec k_1,\vec k_2,\vec k_3}^{\ \nu}
\equiv \sug\,$ for brevity. By Proposition
\ref{prop:Z2-equiv-bundle-2-conn}, there exists a
$\bZ_2$-equivariant isomorphism $\,f :  \widetilde\cM_{\vec k_1,
\vec k_2,\vec k_3}\vert_{\widetilde\cT_{\la,\mu;\vec k_1,\vec k_2,
\vec k_3}^{\ \nu}} \to \widetilde E_{\la,\mu;\vec k_1,\vec k_2,\vec
k_3}^{\ \nu}$.\ Since both bundles are just trivial bundles over
$\,\sug$,\ we can write $\,f(h,z) = \bigl(h,b(h)\cdot z\bigr)\,$ for
some $\uj$-valued map $\,b\,$ on $\,\sug$.\ The compatibility of
$\,f\,$ with the $\bZ_2$-action implies the following identity for
$\,b$:
\qq\label{eq:chi-chi-ff}
\frac{\unl{\chi}_{\la,\mu;\vec k_1,\vec k_2,\vec k_3}^{\ \nu}(\z
)}{\eta_{\vec k_1,\vec k_2,\vec k_3}(\z)}=\frac{b(\zeta\cdot h)}{b(h)}
\qquad \text{for all $\,h \in \sug\,$ and $\,\z\in\{e,-e\}$}\,.
\qqq
Its left-hand side takes values in the set $\,\{-1,+1\}$.\ Suppose
that $\,b(-h) = - b(h)\,$ for all $\,h \in \sug$.\ Restricting
$\,b\,$ to the maximal torus $\,\uj \subset \sug$,\ we obtain a map
$\,\uj \to \uj$.\ This map necessarily has a non-zero winding
number, as illustrated by the simple calculation
\be
\tfrac{1}{2\pi\sfi}\,\int_{\uj}\,\hspace{-.5em}\sfd\log b =
\tfrac{1}{2\pi\sfi}\,
\int_0^{2\pi}\,\hspace{-.5em}\sfd\phi\,\tfrac{\sfd\ }{\sfd\phi}\log
b\bigl(e^{\sfi \sigma_3 \phi}\bigr) =
\tfrac{1}{\pi\sfi}\,\int_0^\pi\,
\hspace{-.5em}\sfd\phi\,\tfrac{\sfd\ }{\sfd\phi}\log b\bigl(e^{\sfi
\sigma_3 \phi} \bigr) = \tfrac{1}{\pi\sfi}\,\log\frac{b(-e)}{b(e)}
\in 2\bZ+1\,,
\ee
where we used $\,\tfrac{\sfd\ }{\sfd\phi}\log b\bigl(e^{\sfi\sigma_3
(\phi+\pi)}\bigr) = \tfrac{\sfd\ }{\sfd\phi}\log b\bigl(e^{\sfi
\sigma_3 \phi}\bigr)$.\ However, the maximal torus in $\,\sug\,$ is
the boundary of a disc, and $\,b\,$ restricts to a smooth map on
this disc, and so cannot have non-zero winding. Hence,
$\,b(-h)/b(h)=1$,\ which proves the lemma.
\end{proof}

Taking into account the last lemma, we can now use the known form of
$\,\unl{\chi}_{\la,\mu;\vec k_1,\vec k_2,\vec k_3}^{\ \nu}\,$ for
$\,(\la,\mu,\nu)\,$ from the boundary $\,\dot\xcF\,$ to determine
$\,\eta_{\vec k_1,\vec k_2,\vec k_3}$,\ and subsequently employ this
result to deduce the characters $\,\unl{\chi}_{\la,\mu;\vec k_1,\vec
k_2,\vec k_3}^{\ \nu}\,$ for $\,(\la,\mu,\nu)\,$ from the interior
of $\,\xcF$.\ This is done in the following lemma.

\belem\label{lem:fix-int-equiv-str}
Let $\,(\la,\mu,\nu)\in\xcF\,$ and take an arbitrary element $\,\z=
(-1)^{\vep_\z}\,e\in\bZ_2$.\ Then,
\qq
\unl{\chi}_{\la,\mu;\vec k_1,\vec k_2,\vec k_3}^{\ \nu}(\z)=(-1
)^{\vep_\z(\la_{l_3\dotplus m_3}-\la_{l_1}-\la_{m_2})}\,.
\qqq
\elem

\begin{proof}
By virtue of Lemma \ref{lem:chi-chi}, it suffices to verify the
thesis for $\,(\la,\mu,\nu)\in\p_{12}\xcF$.\ Lemma
\ref{lem:multipl-str-restr-dF} shows that the relevant
$\bZ_2$-equivariant structure is simply a restriction of the
$\uj$-equivariant structure of the KKS bundle
$\,\overline K_{2\sfk \,n_{\la,\mu}^{\
\nu}-\la_{l_1}-\vep(a_{\la,\mu}^{\ \nu})\,
\la_{m_2}+\vep(b_{\la,\mu}^{\ \nu})\,\la_{l_3\dotplus m_3}}$.\ Upon
recalling \eqref{eq:char-la} and using integrality of $\,n_{\la,
\mu}^{\ \nu}\,$ alongside \eqref{eq:epsil-on-sug}, we obtain the
desired result.
\end{proof}

\section{The fusion 2-isomorphism}\label{sec:fusmorph}

The preceding sections have equipped us with all the tools necessary
to address the issue of the existence of the fusion 2-isomorphism
$\,\filmn\,$ and -- in so doing -- prove the main result of our
paper, as expressed in Theorem \ref{thm:V=V}.

We begin by noting that, for any triple $\,(\la,\mu,\nu)\,$ of
weights from the discrete subset
\qq
F_\bZ=\xcF\cap\bigl(\faff\bigr)^{\x 3}
\qqq
of the fusion polytope, the product gerbe
$\,\txp_1^*\cG\star\txp_2^*\cG\,$ discussed previously admits two
different trivialisations over the corresponding manifold
$\,\Tlmn$,\ namely,
\qq
\txp_1^*\Phi_\la\star\txp_2^*\Phi_\mu\ :\ \bigl(\txp_1^*\cG\star
\txp_2^*\cG\bigr)\big\vert_{\Tlmn}\xrightarrow{\sim}\cI\bigl(
\txp_1^*\om_\la+\txp_2^*\om_\mu\bigr)\big\vert_{\Tlmn}
\qqq
and
\qq
\bigl(\txm^*\Phi_\nu\star\id_{\cI(\rho)}\bigr)\circ\cM\ :\ \bigl(
\txp_1^*\cG\star\txp_2^*\cG\bigr)\big\vert_{\Tlmn}\xrightarrow{\sim}
\cI\bigl(\txm^*\om_\nu+\rho\bigr)\big\vert_{\Tlmn} \,.
\qqq
Upon invoking Lemma \ref{lem:Om-ker} in conjunction with
\eqref{eq:om_la-def}, we conclude that both trivialisations yield
the same trivial gerbe over $\,\Tlmn$,
\qq
\cI\bigl(\txp_1^*\om_\la+\txp_2^*\om_\mu\bigr)\equiv\cI\bigl(\txm^*
\om_\nu+\rho\bigr)\,,
\qqq
and so it is natural to enquire when there exists a 2-isomorphism
$\,\filmn\,$ between the 1-isomorphisms
$\,\txp_1^*\Phi_\la\star\txp_2^*\Phi_\mu\,$ and
$\,\bigl(\txm^*\Phi_\nu\star\id_{\cI(\rho)}\bigr)\circ\cM$,\ as
depicted in the standard 2-categorial diagram
\qq
\alxydim{@C=6em@R=5em}{\bigl(\txp_1^*\cG\star\txp_2^*\cG\bigr)\big
\vert_{\Tlmn} \ar[d]_{\txp_1^*\Phi_\la\star\txp_2^*\Phi_\mu}
\ar[r]^{\cM\vert_{\Tlmn}} & \bigl(\txm^*\cG\star\cI(\rho)\bigr)\big
\vert_{\Tlmn} \ar[d]^{\txm^*\Phi_\nu\star\id_{\cI(\rho)}}
 \\ \cI\bigl(\txp_1^*\om_\la+\txp_2^*\om_\mu
\bigr)\big \vert_{\Tlmn} \ar@{=>}[ur]|{\filmn} \ar@{=}[r] &
\cI\bigl(\txm^*\om_\nu+\rho\bigr)\big \vert_{\Tlmn}}\,.
\qqq
This 2-isomorphism was dubbed the fusion 2-isomorphism in the
introduction, with reference to its underlying physical
interpretation detailed in \cite{Runkel:2009}.

\belem\label{lem:2-iso-exists}
A 2-isomorphism $\,\filmn\,$ exists if and only if $\,\la+\mu-\nu\in
2\bZ$,\ in which case it is unique up to a globally defined
$\uj$-valued constant. \elem

\begin{proof}
Let us give the fusion 2-isomorphism in full detail by specialising
Definition \ref{def:2-iso} to the setting at hand. To this end, we
first write out the line bundles and the bundle isomorphisms of the
two 1-isomorphisms that appear in its definition, in keeping with
Definition \ref{def:stab-iso}. Starting with
$\,\txp_1^*\Phi_\la\star \txp_2^*\Phi_\mu$,\ we find -- for a fixed
pair $\,(\vec k_1,\vec k_2)=(l_1,m_1,l_2,m_2)\in\{0,1\}^{\x 4}\,$ --
the surjective submersion
\qq
\bigl(\cT_{\la,\mu;\vec k_1}^{\ \nu}\x_{\Tlmn}\Tlmn\bigr)\x_{\Tlmn}
\bigl(\cT_{\la,\mu;\vec k_2}^{\ \nu}\x_{\Tlmn}\Tlmn\bigr)\cong
\cT_{\la,\mu;\vec k_1,\vec k_2}^{\ \nu}
\qqq
as the base of the product line bundle
\qq
E_{\la\ox\mu}\vert_{\cT_{\la,\mu;\vec k_1,\vec k_2}^{\ \nu}}=
\pr_1^*\widehat\txp_1^*E_{\la;l_1}\ox\pr_2^*\widehat\txp_2^*E_{\mu;
m_2}\to\cT_{\la,\mu;\vec k_1,\vec k_2}^{\ \nu}\,,
\qqq
written in terms of the canonical projections $\,\pr_i:\cT_{\la,\mu;
\vec k_1,\vec k_2}^{\ \nu}\to\cT_{\la,\mu;\vec k_i}^{\ \nu}$,\
alongside the product bundle isomorphism
\qq
\a_{\la\ox\mu}=\pr_{1,3}^*\widehat\txp_1^{[2]*}\a_\la\ox\pr_{2,4}^*
\widehat\txp_2^{[2]*} \a_\mu\ :\ \pr_{1,3}^*\widehat\txp_1^{[2]*}L
\ox\pr_{2,4}^*\widehat\txp_2^{[2]*}L\ox\pr_{3,4}^*E_{\la\ox\mu}
\xrightarrow{\sim}\pr_{1,2}^*E_{\la\ox\mu}\,,
\qqq
of line bundles over a disjoint union of spaces $\,\cT_{\la,\mu;\vec
k_1,\vec k_2,\vec k_3,\vec k_4}^{\ \nu}$,\ each equipped with the
canonical projections $\,\pr_{i,j}:\cT_{\la,\mu;\vec k_1,\vec k_2,
\vec k_3,\vec k_4}^{\ \nu}\to\cT_{\la,\mu;\vec k_i,\vec k_j}^{\
\nu}$.\ In the case of the composite stable isomorphism
$\,\bigl(\txm^*\Phi_\nu\star\id_{\cI(\rho)}\bigr)\circ\cM$,\ we
obtain -- for fixed $\,\vec k_i\in\{0,1\}^{\x 2}\,$ -- the
surjective submersion
\qq
\cT_{\la,\mu;\vec k_1,\vec k_2,\vec k_3}^{\ \nu}\x_{\cT_{\la,\mu;
\vec k_3}^{\ \nu}}\bigl(\cT_{\la,\mu;\vec k_3}^{\ \nu}\x_{\Tlmn}
\Tlmn\bigr)\cong\cT_{\la,\mu;\vec k_1,\vec k_2,\vec k_3}^{\ \nu}
\qqq
as the base of the product line bundle
\qq
E_{\nu(\la\cdot\mu)}\vert_{\cT_{\la,\mu;\vec k_1,\vec k_2,\vec
k_3}^{\ \nu}}=E_{\la,\mu;\vec k_1,\vec k_2,\vec k_3}^{\ \nu}\ox
\pr_3^*\widehat\txm^*E_{\nu;l_3\dotplus m_3}\to\cT_{\la,\mu;\vec
k_1,\vec k_2,\vec k_3}^{\ \nu}\,,
\qqq
written in terms of the canonical projections $\,\pr_i:\cT_{\la,\mu;
\vec k_1,\vec k_2,\vec k_3}^{\ \nu}\to\cT_{\la,\mu;\vec k_i}^{\
\nu}$,\ together with the bundle isomorphism
\qq
\a_{\nu(\la\cdot\mu)}\vert_{\cT_{\la,\mu;\vec k_1,\vec k_2,\vec k_3,
\vec k_4,\vec k_5,\vec k_6}^{\ \nu}}&=&\bigl(\id_{\pr_{1,2,3}^*
E_{\la,\mu;\vec k_1,\vec k_2,\vec k_3}^{\ \nu}}\ox\pr_{3,6}^*
\widehat\txm^{[2]*}\a_\nu\vert_{\cC_{\nu;l_3\dotplus m_3,l_6\dotplus
m_6}}\bigr)\cr\cr
&&\circ\bigl(\a_{\la,\mu;\vec k_1,\vec k_2,\vec k_3,\vec k_4,\vec
k_5,\vec k_6}^{\ \nu}\ox\id_{\pr_6^*\widehat\txm^*E_{\nu;l_6
\dotplus m_6}}\bigr)\,,\cr\cr\cr \a_{\nu(\la\cdot\mu)}\ &:&\
\pr_{1,4}^*\widehat\txp_1^{[2]*}L\ox\pr_{2,5}^*\widehat\txp_2^{[2]
*}L\ox\pr_{4,5,6}^*E_{\nu(\la\cdot\mu)}\xrightarrow{\sim}\pr_{1,2,
3}^*E_{\nu(\la\cdot\mu)}
\qqq
of line bundles over a disjoint union of spaces $\,\cT_{\la,\mu;\vec
k_1,\vec k_2,\vec k_3,\vec k_4,\vec k_5,\vec k_6}^{\ \nu}$,\ each
equipped with the canonical projections $\,\pr_i:\cT_{\la,\mu;\vec
k_1,\vec k_2,\vec k_3,\vec k_4,\vec k_5,\vec k_6}^{\ \nu}\to\cT_{\la
,\mu;\vec k_i}^{\ \nu}$,\ $\,\pr_{i,j}: \cT_{\la,\mu;\vec k_1,\vec
k_2,\vec k_3,\vec k_4,\vec k_5,\vec k_6}^{\ \nu}\to\cT_{\la,\mu;
\vec k_i,\vec k_j}^{\ \nu}\,$ and $\,\pr_{i,j,k}: \cT_{\la,\mu;\vec
k_1,\vec k_2,\vec k_3,\vec k_4,\vec k_5,\vec k_6}^{\ \nu}\to\cT_{\la
,\mu;\vec k_i,\vec k_j,\vec k_k}^{\ \nu}$.\ We may now write the
fusion 2-isomorphism as an isomorphism
\qq
\filmn\ :\ \pr_{1,2}^*E_{\la\ox\mu}\xrightarrow{\sim}E_{\nu(\la
\cdot\mu)}
\qqq
of Hermitian line bundles with connection over a disjoint union of
bases $\,\cT_{\la,\mu;\vec k_1,\vec k_2}^{\ \nu}\x_{\cT_{\la,\mu;
\vec k_1,\vec k_2}^{\ \nu}}\cT_{\la,\mu;\vec k_1,\vec k_2,\vec
k_3}^{\ \nu}=\cT_{\la,\mu;\vec k_1,\vec k_2,\vec k_3}^{\ \nu}$,\ the
latter coming with the canonical projection $\,\pr_{1,2}:\cT_{\la,
\mu;\vec k_1,\vec k_2,\vec k_3}^{\ \nu}\to\cT_{\la,\mu;\vec k_1,\vec
k_2}^{\ \nu}$.\ Note that we have implicitly chosen the trivial
surjective submersion $\,\id_{\cT_{\la,\mu;\vec k_1,\vec k_2,\vec
k_3}^{\ \nu}}\,$ for the fusion 2-isomorphism, which can be done by
virtue of Proposition \ref{prop:2-iso-triv-subm}. Imposing the
requirement that $\,\filmn\,$ be compatible with
$\,\a_{\nu(\la\cdot\mu)}\,$ and $\,\a_{\la\otimes\mu}\,$ becomes
equivalent to demanding the commutativity of the diagram
\qq\label{diag:fus-iso-alal}
\alxydim{@C=8em@R=5em}{\pr_{1,4}^*\widehat\txp_1^{[2]*}L\ox\pr_{2,
5}^*\widehat\txp_2^{[2]*}L\ox\pr_{4,5}^*E_{\la\ox\mu}
\ar[r]^{\hspace{2.2cm}\pr_{1,2,4,5}^*\a_{\la\ox\mu}}
\ar[d]_{\id_{\pr_{1,4}^*\widehat\txp_1^{[2]*}L\ox\pr_{2,5}^*
\widehat\txp_2^{[2]*}L}\ox\pr_{4,5,6}^*\filmn}
& \pr_{1,2}^*E_{\la\ox\mu} \ar[d]^{\pr_{1,2,3}^*\filmn
}
\\ \pr_{1,4}^*\widehat\txp_1^{[2]*}L\ox\pr_{2,5}^*\widehat\txp_2^{[
2]*}L\ox\pr_{4,5,6}^*E_{\nu(\la\cdot\mu)} \ar[r]_{\hspace{2.2cm}
\a_{\nu(\la\cdot\mu)}} & \pr_{1,2,3}^*E_{\nu(\la\cdot\mu)}}
\qqq
of isomorphisms of Hermitian line bundles with connection over a
disjoint union of bases $\,\cT_{\la,\mu;\vec k_1,\vec k_2,\vec k_3,
\vec k_4,\vec k_5,\vec k_6}^{\ \nu}$,\ the latter taken with
canonical projections $\,\pr_{i,j}$,\ $\,\pr_{i,j,k}\,$ as above and
with $\,\pr_{i,j,k,l}:\cT_{\la,\mu; \vec k_1,\vec k_2,\vec k_3,\vec
k_4,\vec k_5,\vec k_6}^{\ \nu}\to \cT_{\la,\mu;\vec k_i,\vec
k_j,\vec k_k,\vec k_l}^{\ \nu}$.

In what follows, we work with the $\cS_{\la,\mu}^{\
\nu}$-equivariant counterparts $\,\widetilde E_{\la\ox\mu},
\widetilde\a_{\la\ox\mu},\widetilde E_{\nu(\la\cdot\mu)},\widetilde
\a_{\nu(\la\cdot\mu)}\,$ and $\,\widetilde\varphi_{\la,\mu}^{\
\nu}\,$ of the untilded objects, defined on the respective
surjective submersions over $\,\widetilde\cT_{\la,\mu}^{\ \nu}$.\
The above abstract definition of the fusion 2-isomorphism then
yields, on $\,\cT_{\la,\mu;\vec k_1,\vec k_2,\vec k_3}^{\ \nu}$,
\qq\label{eq:filmn-equiv-loc}
\widetilde\varphi_{\la,\mu}^{\ \nu}(\vec k_1,\vec k_2,\vec k_3,\la,
\mu,\nu,h,z\ox z')=(\vec k_1,\vec k_2,\vec k_3,\la,\mu,\nu,h,f_{\la,
\mu;\vec k_1,\vec k_2,\vec k_3}^{\ \nu}\cdot z\ox z')
\qqq
for some constant phases $\,f_{\la,\mu;\vec k_1,\vec k_2,\vec
k_3}^{\ \nu}\in\uj$.\ Constancy of $\,f_{\la,\mu;\vec k_1,\vec k_2,
\vec k_3}^{\ \nu}\,$ is an immediate consequence of the equality of
the connection 1-forms of the line bundles $\,\widetilde
E_{\nu(\la\cdot\mu)} \vert_{\widetilde\cT_{\la,\mu;\vec k_1,\vec
k_2,\vec k_3}^{\ \nu}}\,$ and $\,\widetilde E_{\la\ox
\mu}\vert_{\widetilde\cT_{\la, \mu;\vec k_1,\vec k_2}^{\ \nu}}$,\
and of connectedness of $\,\widetilde\cT_{\la,\mu;\vec k_1,\vec
k_2,\vec k_3}^{\ \nu}$.\ The compatibility condition encoded in
Diagram \ref{diag:fus-iso-alal} is rewritten as
\qq
f_{\la,\mu;\vec k_4,\vec k_5,\vec k_6}^{\ \nu}=f_{\la,\mu;\vec k_1,
\vec k_2,\vec k_3}^{\ \nu}
\qqq
and implies independence of $\,f_{\la,\mu;\vec k_1,\vec k_2,\vec
k_3}^{\ \nu}\,$ of the indices $\,\vec k_1,\vec k_2,\vec k_3$,
\qq
f_{\la,\mu;\vec k_1,\vec k_2,\vec k_3}^{\ \nu}=f_{\la,\mu}^{\ \nu}
\in\uj\,.
\qqq
The remaining freedom in the definition of $\,\varphi_{\la,\mu}^{\
\nu}\,$ is hence a $\,\uj$-valued constant. Clearly, for the fusion
2-isomorphism to be well-defined, \eqref{eq:filmn-equiv-loc} has to
be consistent with the $\cS_{\la,\mu}^{\ \nu}$-equivalences entering
the definitions of the bundles involved, as expressed by the
identity
\qq
\widetilde\varphi_{\la,\mu}^{\ \nu}(\vec k_1,\vec k_2,\vec k_3,\la,
\mu,\nu,h\cdot g,\chi_{\la-\la_{l_1}}(g)\cdot z\ox\chi_{\mu-
\la_{m_2}}(g)^{\vep(a_{\la,\mu}^{\ \nu})}\cdot z')\cr\cr =(\vec
k_1,\vec k_2,\vec k_3,\la,\mu,\nu,h\cdot g,\unl{\chi}_{\la,\mu;\vec
k_1,\vec k_2,\vec k_3}^{\ \nu}(g)\cdot f_{\la,\mu}^{\ \nu}\cdot z\ox
\chi_{\nu-\la_{l_3\dotplus m_3}}(g)^{\vep(b_{\la,\mu}^{\ \nu})}\cdot
z')\,,
\qqq
to be imposed for an arbitrary element $\,g\in\cS_{\la,\mu}^{\
\nu}$,\ cf.\ \eqref{eq:funny-act-1}--\eqref{eq:pull-ab-out-chi}.
Thus, we obtain the algebraic condition
\qq\label{eq:D=1}
\D_{\la,\mu;\vec k_1,\vec k_2,\vec k_3}^{\ \nu}(g)=\unl{\chi}_{\la,
\mu;\vec k_1,\vec k_2,\vec k_3}^{\ \nu}(g)\cdot\chi_{\nu-\la_{l_3
\dotplus m_3}}(g)^{\vep(b_{\la,\mu}^{\ \nu})}\cdot\chi_{\la-
\la_{l_1}}(g)^{-1}\cdot\chi_{\mu-\la_{m_2}}(g)^{-\vep(a_{\la,
\mu}^{\ \nu})}=1\,.
\qqq
The fusion 2-isomorphism $\,\filmn\,$ exists if and only if this
condition is satisfied for all $\,\vec k_1,\vec k_2,\vec k_3\in\{0,
1\}\,$ such that $\,\cT_{\la,\mu;\vec k_1,\vec k_2,\vec k_3}^{\ \nu}
\neq\emptyset\,$ and for all $\,g\in\cS_{\la,\mu}^{\ \nu}$.\ We
shall solve this condition by distinguishing two cases:
\begin{itemize}
\item $(\la,\mu,\nu)\in\p\xcF\cap\bigl(\faff \bigr)^{\x 3}$:\
For $\,(\la,\mu,\nu)\in\p_0\xcF\cap\bigl(\faff \bigr)^{\x 3}$,\ the
identity $\,\D_{\la,\mu;\vec k_1,\vec k_2,\vec k_3}^{\ \nu}(g)=1\,$
is automatically implied by triviality of the characters involved.
For $\,(\la,\mu,\nu)\in\p_{12}\xcF\cap\bigl(\faff\bigr)^{\x 3}$,\ on
the other hand, it is readily verified by inspection, using Lemma
\ref{lem:multipl-str-restr-dF} and \eqref{eq:char-la}.
\item $(\la,\mu,\nu)\in\mathring\xcF\cap\bigl(\faff\bigr)^{\x
3}$:\ In this case, $\,\D_{\la,\mu;\vec k_1,\vec k_2,\vec k_3}^{\
\nu}(\z)=(-1)^{\vep_\z\,(\la+\mu-\nu)}\,$ for all $\z\in\cS_{\la,
\mu}^{\ \nu}\cong\bZ_2$,\ as follows directly from Lemma
\ref{lem:fix-int-equiv-str} taken in conjunction with the identities
$\,\vep(a_{\la,\mu}^{\ \nu})=1=\vep(b_{\la,\mu}^{\ \nu})\,$ and
\eqref{eq:char-la}.
\end{itemize}
\end{proof}

Lemma \ref{lem:2-iso-exists} now implies the main result of the
paper,

\medskip\noindent
{\bf Theorem \ref{thm:V=V}.}
$V=V_\cG$,\ where $\,V\,$ is given in \eqref{eq:setV-def} and
$\,V_\cG\,$ in \eqref{eq:setV'-def}.


\small

\begin{thebibliography}{XXX}
\setlength{\itemsep}{-0em}

\bibitem[Al]{Alvarez:1984es}
O.~Alvarez, {\it Topological quantization and cohomology},
\doi{10.1007/BF01212452}{Comm.\ Math.\ Phys.\ {\bf 100} (1985)
279--309}.

\bibitem[AM]{Alekseev:1993rj}
A.Yu.~Alekseev and A.Z.~Malkin, {\it Symplectic structure of the
moduli space of flat connection on a Riemann surface},
\doi{10.1007/BF02101598}{Comm.\ Math.\ Phys.\ {\bf 169} (1995)
445--495} \arxiv{dg-ga/9707021}{[dg-ga/9707021]}.

\bibitem[AMW]{Alekseev:2000}
A.~Alekseev, E.~Meinrenken and C.~Woodward, {\it Formulas of Verlinde type for non-simply connected groups},
\arxiv{math/0005047}{math.SG/0005047}.

\bibitem[Be]{Beauville:1994}
A.~Beauville, {\it Conformal blocks, fusion rules and the Verlinde
formula}, in: Proceedings of the Hirzebruch 65 Conference on
Algebraic Geometry, Israel Math.\ Conf.\ Proc.\ {\bf 9} (1996)
75--96 \arxiv{alg-geom/9405001}{[alg-geom/9405001]}.

\bibitem[BL]{Bismut:1999}
J.-M.~Bismut, F.~ Labourie, {\it Symplectic geometry and Verlinde
formulas}, in:
\httpurl{www.intlpress.com/books/SDG/SDG-V.php}{Surveys in
Differential Geometry, Vol.\,V, International Press, 1999}.

\bibitem[Bo]{Bouwknegt:2001vu}
P.~Bouwknegt and A.L.~Carey and V.~Mathai and M.K.~Murray and
D.~Stevenson, {\it Twisted K-theory and K-theory of bundle gerbes},
\doi{doi:10.1007/s002200200646}{Comm.\ Math.\ Phys.\ {\bf 228}
(2002) 17--49}, \arxiv{hep-th/0106194}{[hep-th/0106194]}.

\bibitem[Br]{Brylinski:1993ab}
J.-L.~Brylinski, {\it Loop Spaces, Characteristic Classes and
Geometric Quantization}, in: Progress in Mathematics, vol.\,107,
Birkh{\"a}user, 1993.

\bibitem[Ca]{Carey:2004xt}
A.L.~Carey and S.~Johnson and M.K.~Murray and D.~Stevenson and
B.-L.~Wang, {\it Bundle gerbes for Chern--Simons and
Wess--Zumino--Witten theories},
\doi{10.1007/s00220-005-1376-8}{Comm.\ Math.\ Phys.\ {\bf 259}
(2005) 577--613} \arxiv{math/0410013}{[math/0410013]}.

\bibitem[CJM]{Carey:2002}
A.L.~Carey, S.~Johnson and M.K.~Murray, {\it Holonomy on D-branes},
\doi{doi:10.1016/j.geomphys.2004.02.008}{J.\ Geom.\ Phys.\ {\bf 52}
(2004) 186--216} \arxiv{hep-th/0204199}{[hep-th/0204199]}.

\bibitem[FGK]{Felder:1988sd}
G.~Felder, K.~Gaw{\c{e}}dzki and A.~Kupiainen, {\it Spectra of
Wess-Zumino-Witten models with arbitrary simple groups},
\doi{10.1007/BF01228414}{Comm.\ Math.\ Phys.\ {\bf 117} (1988)
127--158.}

\bibitem[Fr]{Frohlich:2006ch}
J.~Fr\"ohlich, J.~Fuchs, I.~Runkel and C.~Schweigert, {\it Duality
and defects in rational conformal field theory},
\doi{doi:10.1016/j.nuclphysb.2006.11.017}{Nucl.\ Phys.\ {\bf B763}
(2007) 354--430} \arxiv{hep-th/0607247}{[hep-th/0607247].}

\bibitem[FSW]{Fuchs:2007fw}
J.~Fuchs, C.~Schweigert and K.~Waldorf, {\it Bi-branes: Target space
geometry for world sheet topological defects},
\doi{10.1016/j.geomphys.2007.12.009}{J.\ Geom.\ Phys.\ {\bf 58}
(2008) 576--598} \arxiv{hep-th/0703145}{[hep-th/0703145].}

\bibitem[Ga1]{Gawedzki:1987ak}
K.~Gaw{\c{e}}dzki, {\it Topological actions in two-dimensional
quantum field theory}, in: Nonperturbative Quantum Field Theory
(G.~{'t} Hooft, A.~Jaffe, G.~Mack, P.~Mitter and R.~Stora, eds.),
Plenum Press, 1988, pp. 101--141.

\bibitem[Ga2]{Gawedzki:2004tu}
K.~Gaw{\c{e}}dzki, {\it Abelian and non-Abelian branes in WZW models
and gerbes}, \doi{10.1007/s00220-005-1301-1}{Comm.\ Math.\ Phys.\
{\bf 258} (2005) 23--73} \arxiv{hep-th/0406072}{[hep-th/0406072]}.

\bibitem[Gi]{Giraud:1971}
J.~Giraud, {\it Cohomologie non ab\'elienne}, in: Grundlehren der
Mathematischen Wissenschaften, vol.\,179, Springer, (1971).

\bibitem[Go]{Gomi:2003}
K.~Gomi, {\it Equivariant smooth Deligne cohomology},
\httpurl{projecteuclid.org/euclid.ojm/1153494380}{Osaka J.\ Math.\ {\bf 42} (2005) 309--337}
\arxiv{math/0307373}{[math.DG/0307373]}.

\bibitem[GR1]{Gawedzki:2002se}
K.~Gaw{\c{e}}dzki and N.~Reis, {\it WZW branes and gerbes},
\doi{10.1142/S0129055X02001557}{Rev.\ Math.\ Phys.\  {\bf 14} (2002)
1281--1334} \arxiv{hep-th/0205233}{[hep-th/0205233].}

\bibitem[GR2]{Gawedzki:2003pm}
K.~Gaw{\c{e}}dzki and N.~Reis, {\it Basic gerbe over non simply
connected compact groups},
\doi{doi:10.1016/j.geomphys.2003.11.004}{J.\ Geom.\ Phys.\ {\bf 50}
(2003) 28--55} \arxiv{math.DG/0307010}{[math.DG/0307010]}.

\bibitem[GSW1]{Gawedzki:2007uz}
K.~Gaw{\c{e}}dzki, R.R.~Suszek and K.~Waldorf, {\it WZW orientifolds
and finite group cohomology}, \doi{10.1007/s00220-008-0525-2}{Comm.\
Math.\ Phys.\ {\bf 284} (2008) 1--49}
\arxiv{hep-th/0701071}{[hep-th/0701071]}.

\bibitem[GSW2]{Gawedzki:2008um}
K.~Gaw{\c{e}}dzki, R.R.~Suszek and K.~Waldorf, {\it Bundle gerbes
for orientifold sigma models}, submitted to Adv.\ Theor.\ Math.\
Phys., \arxiv{0809.5125}{0809.5125 [math-ph]}.

\bibitem[GW]{Gawedzki:2009jj}
K.~Gaw{\c{e}}dzki and K.~Waldorf, {\it Polyakov--Wiegmann formula
and multiplicative gerbes}, \arxiv{0908.1130}{0908.1130 [hep-th]}.

\bibitem[Ha]{Hayashi:1999}
M.~Hayashi, {\it The moduli space of ${\rm SU}(3)$ flat connections
and the fusion rules}, \httpurl{www.jstor.org/stable/119328}{Proc.\
AMS {\bf 127} (1999) 1545--1555}.

\bibitem[JW]{Jeffrey:1991rp}
L.C.~Jeffrey and J.~Weitsman, {\it Bohr-Sommerfeld orbits in the
moduli space of flat connections and the Verlinde dimension
formula}, \doi{10.1007/BF02096964}{Comm.\ Math.\ Phys.\  {\bf 150}
(1992) 593--630}.

\bibitem[Ki]{Kirillov:1975}
A.A.~Kirillov, {\it Elements of the Theory of Representations}, in:
Grundlehren der mathematischen Wissenschaften, vol.\,220, Springer,
1975.

\bibitem[Ko]{Kostant:1970}
B.~Kostant, {\it Quantization and Unitary Representations}, in:
Lecture Notes in Mathematics, vol.\,170, Springer, 1970.

\bibitem[Me]{Meinrenken:2002}
E.~Meinrenken, {\it The basic gerbe over a compact simple Lie
group},
\httpurl{retro.seals.ch/cntmng?type=pdf&rid=ensmat-001:2003:49::234&subp=hires}{Enseign.
Math. {\bf 49} (2003) 307--333}
\arxiv{math.DG/0209194}{[math.DG/0209194]}.

\bibitem[MS]{Murray:1999ew}
M.K.~Murray and D.~Stevenson, {\it Bundle gerbes: stable isomorphism
and local theory},
\httpurl{jlms.oxfordjournals.org/cgi/content/abstract/62/3/925}{J.\
Lond.\ Math.\ Soc.\ {\bf 62} (2000) 925--937}
\arxiv{math.DG/9908135}{[math.DG/9908135]}.

\bibitem[Mu]{Murray:1994db}
M.K.~Murray, {\it Bundle gerbes},
\httpurl{jlms.oxfordjournals.org/cgi/content/abstract/54/2/403}{J.\
Lond.\ Math.\ Soc.\ {\bf 54} (1996) 403--416}
\arxiv{dg-ga/9407015}{[dg-ga/9407015]}.

\bibitem[PW]{Polyakov:1984et}
A.M.~Polyakov and P.B.~Wiegmann, {\it Goldstone fields in two
dimensions with multivalued actions},
\doi{doi:10.1016/0370-2693(84)90206-5}{Phys.\ Lett.\ {\bf B141}
(1984) 223--228}.

\bibitem[RS1]{Runkel:2008gr}
I.~Runkel and R.R.~Suszek, {\it Gerbe-holonomy for surfaces with
defect networks},
\httpurl{www.intlpress.com/ATMP/ATMP-issue_13_4.php}{Adv.\ Theor.\
Math.\ Phys.\ {\bf 13} (2009) 1137--1219}
\arxiv{0808.1419}{[0808.1419 [hep-th]]}.

\bibitem[RS2]{Runkel:2009}
I.~Runkel and R.R.~Suszek, {\it Maximally symmetric defects with
junctions in the classical WZW model}, in preparation.

\bibitem[SS]{Sarkissian:2008dq}
G.~Sarkissian and C.~Schweigert, {\it Some remarks on defects and
T-duality}, \doi{doi:10.1016/j.nuclphysb.2009.04.016}{Nucl.\ Phys.\
{\bf B819} (2009) 478--490} \arxiv{0810.3159}{[0810.3159 [hep-th]]}.

\bibitem[So]{Souriau:1970}
J.-M.~Souriau, {\it Structure des syst\`emes dynamiques}, Dunod,
1970.

\bibitem[SSW]{Schreiber:2005mi}
U.~Schreiber and C.~Schweigert and K.~Waldorf, {\it Unoriented WZW
models and holonomy of bundle gerbes},
\doi{10.1007/s00220-007-0271-x}{Comm.\ Math.\ Phys.\ {\bf 274}
(2007) 31--64} \arxiv{hep-th/0512283}{[hep-th/0512283]}.

\bibitem[St]{Stevenson:2000wj}
D.~Stevenson, {\it The Geometry of Bundle Gerbes}, Ph.D.\ Thesis,
University of Adelaide, 2000,
\arxiv{math.DG/0004117}{math.DG/0004117}.

\bibitem[TW]{Teleman:2000}
C.~Teleman and C.~Woodward, {\it Parabolic bundles, products of
conjugacy classes, and Gromov--Witten invariants}
\httpurl{www.numdam.org/item?id=AIF_2003__53_3_713_0}{Ann.\ Inst.\
Fourier {\bf 53} (2003) 713--748}
\arxiv{math/0012241}{[math.AG/0012241]}.

\bibitem[Ve]{Verlinde:1988sn}
E.P.~Verlinde, {\it Fusion rules and modular transformations in 2D
conformal field theory}, \doi{10.1016/0550-3213(88)90603-7}{Nucl.\
Phys.\ {\bf B300} (1988) 360--376}.

\bibitem[Wa1]{Waldorf:2007mm}
K.~Waldorf, {\it More morphisms between bundle gerbes},
\httpurl{www.tac.mta.ca/tac/volumes/18/9/18-09abs.html}{Theory\
Appl.\ Categories\ {\bf 18} (2007) 240--273}
\arxiv{math.CT/0702652}{[math.CT/0702652]}.

\bibitem[Wa2]{Waldorf:2008mult}
K.~Waldorf, {\it Multiplicative bundle gerbes with connection},
\doi{10.1016/j.difgeo.2009.10.006}{Differential Geom.\ Appl.\ {\bf
28} (2010) 313--340} \arxiv{0804.4835}{[0804.4835 [math.DG]]}.

\bibitem[Wi1]{Witten:1983ar}
E.~Witten, {\it Non-abelian bosonization in two dimensions},
\doi{10.1007/BF01215276}{Comm.\ Math.\ Phys.\ {\bf 92} (1984)
455--472}.

\bibitem[Wi2]{Witten:1988hf}
E.~Witten, {\it Quantum field theory and the Jones polynomial},
\doi{10.1007/BF01217730}{Comm.\ Math.\ Phys.\  {\bf 121} (1989)
351--399}.

\bibitem[Wo]{Woodhouse:1992de}
N.M.J.~Woodhouse, {\it Geometric Quantization}, in: Oxford
Mathematical Monographs, Oxford University Press, 1992.

\end{thebibliography}
\end{document}